\documentclass[conference, onecolumn]{IEEEtran}

\usepackage{verbatim}   
\usepackage{algorithm}
\usepackage{algorithmic}
\usepackage{times}
\usepackage{amsmath}\usepackage{amssymb}
\usepackage{latexsym}
\usepackage{graphicx, color}
\usepackage{epsfig}
\usepackage{url}
\usepackage{caption}
\usepackage{psfrag}
\usepackage{cite}
\usepackage{subfigure}

\renewcommand{\baselinestretch}{1}

 \newcommand{\bnum}{\begin{enumerate}}
 \newcommand{\enum}{\end{enumerate}}
 \newcommand{\ba}{\begin{array}}
 \newcommand{\ea}{\end{array}}
 \newcommand{\bitm}{\begin{itemize}}
 \newcommand{\eitm}{\end{itemize}}
 \newcommand{\beq}{\begin{equation}}
 \newcommand{\eeq}{\end{equation}}
 \newcommand{\beqn}{\begin{eqnarray}}
 \newcommand{\eeqn}{\end{eqnarray}}
 \newcommand{\beqno}{\begin{eqnarray*}}
 \newcommand{\eeqno}{\end{eqnarray*}}
 \newcommand{\bma}{\begin{displaymath}}
 \newcommand{\ema}{\end{displaymath}}
 \newcommand{\bnu}{\begin{enumerate}}
 \newcommand{\enu}{\end{enumerate}}
 \newcommand{\bce}{\begin{center}}
 \newcommand{\ece}{\end{center}}
 \newcommand{\btb}{\begin{tabular}}
 \newcommand{\etb}{\end{tabular}}
 \newcommand{\bmat}{\begin{pmatrix}}
 \newcommand{\emat}{\end{pmatrix}}

\newtheorem{theorem}{\textbf{Theorem}}
\newtheorem{lemma}{\textbf{Lemma}}
\newtheorem{corollary}{\textbf{Corollary}}

\newtheorem{observation}{\textbf{Observation}}
\newtheorem{remark}{\textbf{Remark}}
\newtheorem{definition}{Definition}

\begin{document}

\title{Dynamic Server Allocation 
over Time Varying Channels with Switchover Delay}

\author{G\"{u}ner D. \c{C}elik,~\IEEEmembership{Student Member,~IEEE,}
        Long B. Le,~\IEEEmembership{Member,~IEEE,}
        Eytan~Modiano,~\IEEEmembership{Senior Member,~IEEE}
\thanks{G\"{u}ner D. \c{C}elik and Eytan Modiano are with the Laboratory
for Information and Decision Systems, Massachusetts Institute
of Technology, Cambridge, MA 02139, gcelik@mit.edu,
modiano@mit.edu}
\thanks{Long B. Le is with INRS - EMT, University of Quebec, Montréal (Québec), Canada,
H5A 1K6, long.le@emt.inrs.ca.}
\thanks{A preliminary version of this paper was
presented in IEEE~INFOCOM'11, Apr.\ 2011 \cite{Celik_Infocom11}.}}

\markboth{IEEE TRANSACTIONS ON INFORMATION THEORY, Submitted May~2011}%
{Celik \MakeLowercase{\textit{et al.}}}

%
%


\maketitle\thispagestyle{empty}

\begin{abstract}
We consider a dynamic server allocation problem over parallel queues
with \emph{randomly varying connectivity} and \emph{server
switchover delay} between the queues. At each time slot the server
decides either to stay with the current queue 
or switch to another queue based on the current connectivity and
the queue length information. 
Switchover delay occurs in many telecommunications applications
and is a new modeling component of this problem that has not been previously addressed.
%
We show that the simultaneous presence of
randomly varying
connectivity
and switchover delay 
changes the system stability region 
and the structure of optimal policies. 
%
%
%
In the first part of the paper, we consider a system of two parallel queues, and develop a novel approach to explicitly characterize the stability region of the
system using \emph{state-action frequencies} which are
stationary solutions to a Markov Decision Process (MDP) formulation.
%
We then develop a frame-based
dynamic control (FBDC) policy, 
based on the state-action frequencies, 
and show that it is throughput-optimal asymptotically in the frame length.
The FBDC policy 
is applicable to a broad class of network control systems 
and provides \emph{a new framework for developing
throughput-optimal network control policies} using state-action
frequencies.
Furthermore, we develop simple \emph{Myopic policies} that provably achieve more
than $90\%$ of the stability region. 

In the second part of the paper we extend our results to 
systems with an arbitrary number of queues.
In particular, we show that the stability region characterization in terms of state-action frequencies and the throughput-optimality of the FBDC policy follow for the general case.
Furthermore, we characterize an outer bound on the stability region and an upper bound on sum-throughput
and show that a simple Myopic policy can achieve this sum-throughput upper-bound in the corresponding saturated system.
Finally, simulation results show that the Myopic policies may
achieve the full stability region and are more delay efficient than
the FBDC policy in most cases.
\end{abstract}



\section{Introduction}

\IEEEPARstart{S}{cheduling} a dynamic{\renewcommand{\thefootnote}{}\footnote{
This work was supported by NSF grants CNS-0626781 and CNS-0915988, and by ARO Muri grant number W911NF-08-1-0238.}}
\addtocounter{footnote}{-1} server over time varying wireless channels is
an important and well-studied research problem
which provides useful mathematical
modeling for many practical applications \cite{AltKush00,EryilOzMod07,LinShr05,ModShahZuss06,neely03,neely05,Stolyar04,tass92,tass93,WuSri06,shakk08}.
However, to the best of our knowledge, the joint
effects of \emph{randomly varying connectivity} and \emph{server
switchover delay} 
have not been considered previously.
In fact, switchover delay is a widespread phenomenon that can be observed in many practical
network systems. In satellite systems where a mechanically steered antenna
is providing service to 
ground stations, the time to
switch from one station to another can be around 10ms \cite{BlakeLong09}, \cite{TolShu96}.
Similarly, the delay for electronic beamforming can
be more than $300\mu s$ in wireless radio systems \cite{AkyilWang09}, \cite{BlakeLong09}, \cite{TolShu96},
and in optical
communication systems tuning delay for transceivers can take
significant time ($\mu$s-ms) \cite{Berz}, \cite{ModBarry00}.

\begin{figure}
\centering \psfrag{1\r}[l][][1]{$\!\!\lambda_1$}
\psfrag{2\r}[l][][1]{$\lambda_2$}
\psfrag{3\r}[l][][1]{$\lambda_3$}
\psfrag{4\r}[l][][1]{$\lambda_N$}
\psfrag{0\r}[l][][0.85]{\!\!\!\!\!\!\!\!Server}
\psfrag{5\r}[l][][1]{$\!\!\!\!\!\!\!C_1\!(t)$}
\psfrag{6\r}[l][][1]{$\!\!\! \!\!\!\!\!\!C_2\!(t)$}
\psfrag{7\r}[l][][1]{$\!\!\!C_3\!(t)$}
\psfrag{8\r}[l][][1]{$\!\! C_N\!(t)$}
\psfrag{a\r}[l][][1]{$\!t_s$}
\includegraphics[width=0.40\textwidth]{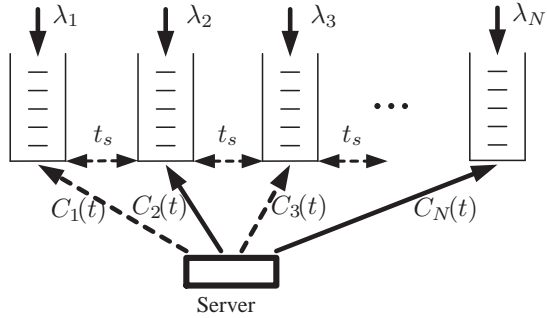}
\captionsetup{font=footnotesize}\caption{System model.
Parallel queues with randomly varying connectivity processes $C_1(t),C_2(t),...,C_N(t)$
and $t_s=1\; \textrm{slot}$ switchover time. 
}\label{Fig:Parallel_queues_servers_switchover_Rconnect}
\end{figure}

We consider the dynamic server control problem for parallel
queues with time varying channels and server switchover
delay as shown in Fig.~\ref{Fig:Parallel_queues_servers_switchover_Rconnect}. 
We consider a slotted system where the slot length is equal to a
single packet transmission time and it takes one slot for the server
to switch from one queue to another\footnote{In a slotted system,
even a minimal switchover delay will lead to a loss of a slot due to
synchronization issues.}. One packet is successfully received from
queue $i$ if the server is currently at queue $i$, it decides to
stay at queue $i$, and queue $i$ is connected, where $i\in\{1,2,...,N\}$ and $N$ is the total number of queues. The server
dynamically decides to stay with the current queue or switch to another queue based on the connectivity and the queue
length information. 
Our goal is to study 
the impact of the simultaneous presence of switchover delays and randomly varying connectivity
%
%
on system stability and to develop optimal control
algorithms.
%
%
We show that the stability region changes as a function of the memory in the channel processes, and it
is significantly reduced as compared to systems without
switchover delay.
Furthermore, we show that throughput-optimal policies take a very different structure from the celebrated Max-Weight algorithm or its variants.
%
%
%
%

\subsection{Main Results}

In the first part of the paper we consider a two-queue system and develop fundamental insights for the problem. We first consider the case of memoryless (i.i.d.) channels where we characterize the stability region explicitly and show that simple Exhaustive type policies that ignore the current queue size and channel state information are throughput-optimal. Next we consider the Gilbert-Elliot channel model \cite{AhmadTaraKrish09}, \cite{Gilbert60} which is a commonly used model to abstract physical channels with memory.
We develop a novel methodology to characterize
the stability region of the 
system using \emph{state-action frequencies} 
which are steady-state solutions to an MDP 
formulation for the corresponding saturated system, and 
characterize the stability region explicitly in terms of the
connectivity parameters. 
%
Using this state-action frequency approach, we develop a novel 
frame-based dynamic control (FBDC)
policy and show that 
it is throughput-optimal asymptotically in the frame length. Our FBDC policy is the only known policy to stabilize systems with randomly varying connectivity and switchover delay and it
is novel in that it 
utilizes the state-action frequencies of the MDP formulation in a dynamic queuing system. 
%
%
%
%
%
%
Moreover, we develop a simple 1-Lookahead Myopic policy that provably achieves 
at least $90\%$ of the stability region, and myopic policies
with
2 and 3 lookahead that achieve more than $94\%$ and $96\%$ of the
stability region respectively.
%
Finally,
we present 
simulation results suggesting that the myopic policies may be
throughput-optimal and more delay efficient than the
FBDC policy. 

In the second part of the paper we consider the general model with arbitrary number of parallel queues. 
For memoryless (i.i.d.) channel processes we explicitly characterize the stability region and the throughput-optimal policy.
For channels with memory, we show that the stability region characterization in terms of state-action frequencies extend to the general case and establish a tight outer bound on the stability region and an upper bound on the sum-throughput explicitly in terms of the connectivity parameters. We quantify the \emph{switching loss} in sum-throughput as compared to the system with no switchover delays and show that simple myopic policies achieve the sum-throughput upper bound in the corresponding saturated system.
We also show that the throughput-optimality of the
FBDC policy extend to the general case.
In fact, the FBDC policy provides \emph{a new framework for achieving throughput-optimal network control} by applying the state-action frequencies of the corresponding saturated system over frames in the dynamic queueing system. The FBDC policy is applicable to a broad class of systems 
whose corresponding saturated model is Markovian with a weakly communicating and finite state space, for example, 
systems with 
arbitrary switchover delays (i.e., systems that take any finite number of time slots for switching)
and general
Markov modulated channel processes. 
%
Moreover, the framework of the FBDC policy can be utilized to achieve throughput-optimality in systems without switchover delay,
%
for instance, in classical network control problems such as 
those considered in \cite{neely03}, \cite{ShahWis06}, \cite{tass93}, \cite{shakk08}.
%


\subsection{Related Work}
Optimal control of queuing systems and communication networks has
been a very active research topic over the past two decades \cite{EryilOzMod07,LinShr05,ModShahZuss06,neely03,neely05,shakk08,Stolyar04,tass92,tass93,WuSri06}. In
the seminal paper \cite{tass92}, Tassiulas and Ephremides
characterized the stability region of multihop wireless networks and proposed 
the throughput-optimal 
Max-Weight scheduling algorithm. 
In \cite{tass93}, the same authors considered a parallel queuing system
with randomly varying connectivity where they 
characterized the stability region of the system explicitly and
proved the throughput-optimality of
the Longest-Connected-Queue scheduling
policy. These results were later extended to joint power allocation and
routing in wireless networks in \cite{neely03,neely05} and optimal scheduling for 
switches in \cite{ShahWis06,Stolyar04}. More recently, suboptimal distributed
scheduling algorithms with throughput guarantees were studied
in \cite{ChapKarSar05,Long10,LinShr05,WuSri06}, while \cite{EryilOzMod07,ModShahZuss06}
developed distributed algorithms that achieve throughput-optimality
(see \cite{Geor_Neely_Tass06}, \cite{Neely10} for a detailed review). The effect of delayed channel state information was considered in \cite{KarLuo07,AnnaEphr09,shakk08} which
showed that the stability region is reduced and that a policy similar to the Max-Weight algorithm is throughput-optimal.

Perhaps the closest problem to ours is that of dynamic server allocation over parallel channels with randomly varying connectivity and limited channel sensing that
has been investigated in \cite{AhmadTaraKrish09,LiNeely10,ZhaoKrishLiu08} under the Gilbert-Elliot channel model.
The saturated system was considered and the optimality of a myopic policy was established for
a single server and two channels in \cite{ZhaoKrishLiu08}, for arbitrary number of channels in \cite{AhmadTaraKrish09}, and for arbitrary number of channels and servers in \cite{AhmadLiu}.
The problem of maximizing the throughput in the network while meeting average delay constraints for a small subset of users was considered in \cite{neely09}. The average delay constraints were turned into penalty functions in \cite{neely09}
and the
the theory of Stochastic Shortest Path problems, which is used for solving Dynamic Programs with certain special structures, was utilized to minimize the resulting drift+penalty terms.
%
%
%
Finally, a partially observable Markov decision process (POMDP) model was used in \cite{ChenZhaoSwami08} to analyze dynamic multichannel access in cognitive radio systems. However, none of these existing
works consider the server switchover delays.

Switchover delay has been considered in Polling models in queuing theory
community (e.g., \cite{AltKonsLiu92}, \cite{Levy},
\cite{LiuNainTow}, \cite{vish06}), however, randomly varying connectivity was not
considered since it may not arise in classical Polling applications.
Similarly, scheduling in optical
networks under reconfiguration delay was considered in
\cite{Berz,Celik_ISIT11}, again in the absence of randomly varying connectivity, where the transmitters and receivers were assumed to be
unavailable during reconfiguration. A detailed survey of the works in this field can be found in
\cite{vish06}.
To the best of our knowledge, this paper is the first to
simultaneously consider random connectivity and server switchover
times.

\subsection{Main Contribution and Organization}
The main contribution of this paper is solving the scheduling
problem in parallel queues with \emph{randomly varying connectivity}
and \emph{server switchover delays} for the first time. 
For this, the paper provides \emph{a novel framework for solving network
control problems} via
characterizing the stability region 
in terms of state-action frequencies 
and achieving throughput-optimality by utilizing the state-action
frequencies over frames.

This paper is organized as follows.
We consider the two-queue system in Section \ref{Sec:2queues} where we
characterize the stability region together with the throughput-optimal policy for memoryless channels. 
We develop the state-action frequency framework in Section~\ref{Sec:Stab_Reg} for channels with memory and use it to explicitly characterize the system stability region. We prove the throughput-optimality of the
FBDC policy in Section
\ref{Sec:Opt_Pol} 
and analyze simple myopic policies 
in Section \ref{Sec:Myopic_Pol}.
We extend our results to the general case in Section \ref{Sec:General_sys}
where we also develop outer bounds on the stability region and an upper bound on the sum-throughput achieved by a simple Myopic policy.
We present simulation results in Section \ref{Sec:Sim} and conclude
the paper in Section \ref{Sec:Conc}.

\section{Two-Queue System}\label{Sec:2queues}
\subsection{System Model}\label{Sec:Model}

Consider two parallel queues with time varying channels and
one server receiving data packets from the queues. 
Time is slotted into unit-length time slots 
equal to one
packet transmission time; $t \in \{0,1,2,...\}$. 
It takes one slot for the server to switch from one queue to the
other, and $m(t)$ denotes
the queue at which the server is present at 
slot $t$. 
Let the i.i.d. stochastic process $A_i(t)$ with average
arrival rate $\lambda_i$ denote the number of packets arriving to
queue $i$ at time slot $t$, where
$\mathbb{E}[A_i^2(t)]\le
A_{\max}^2$, $i \in \{1,2\}$. 
Let $\mathbf{C}(t)=(C_1(t),C_2(t))$
be the channel (connectivity) process 
at time slot $t$, where $C_i(t)= 0$ for the OFF state (disconnected)
and $C_i(t)= 1$ for the ON state (connected).  We assume that the
processes $A_1(t),A_2(t),C_1(t)$ and $C_2(t)$ are independent. 

\begin{figure}
\centering \psfrag{2\r}[l][][1]{$\!\!p_{10}$}
\psfrag{6\r}[l][][1]{$p_{01}$} \psfrag{4\r}[l][][0.85]{\!\!\!\!ON}
\psfrag{5\r}[l][][0.85]{\!\!\!\!\!OFF}
\psfrag{1\r}[l][][1]{\!\!\!\!\!$\!\!1-p_{10}$}
\psfrag{3\r}[l][][1]{\!\!\!\!\!$\!\!\!\!1-p_{01}$}
\includegraphics[width=0.28\textwidth]{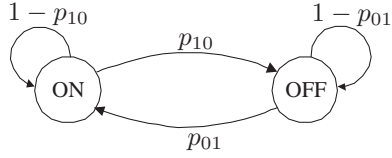}
\captionsetup{font=footnotesize}\caption{Markov
modulated ON/OFF channel process. We have $p_{10}+p_{01}<1$ ($\epsilon < 0.5$) for
positive correlation.}\label{Fig:channel_MC}
\end{figure}

The process $C_i(t),i   \in   \{1,2\}$, is assumed to form the two-state
Markov chain with transition probabilities $p_{10}$ and $p_{01}$  as
shown in Fig.~\ref{Fig:channel_MC}, i.e., the
Gilbert-Elliot channel model 
\cite{AhmadTaraKrish09}, \cite{Gilbert60}, \cite{LiNeely10}, \cite{ZhaoKrishLiu08}, \cite{ZorziRao95}. The
Gilbert-Elliot Channel model 
has been commonly used in modeling and
analysis
of wireless channels with memory \cite{AhmadTaraKrish09}, \cite{LiNeely10}, 
\cite{WangChang96}, \cite{ZhaoKrishLiu08}, \cite{ZorziRao95}. For ease of exposition, we present the analysis in this section for the symmetric Gilbert-Elliot channel model, i.e., $p_{10}=p_{01}=\epsilon$, and we state the corresponding results for the non-symmetric case in Appendix B.
The steady state
probability of each channel state is equal to 0.5 in the symmetric
Gilbert-Elliot channel model. Moreover, for $ \epsilon
 =   0.5$, $C_i(t)   =   1,\textrm{w.p. } 0.5$, 
independently and identically distributed (i.i.d.) at each time slot. 
We refer to this case as the \emph{memoryless channels} case.

Let $\mathbf{Q}(t) = (Q_1(t), Q_2(t))$ be the queue lengths at time
slot $t$. We assume that $\mathbf{Q}(t)$ and $\mathbf{C}(t)$ are
known to the
server at the beginning of each time slot.  
Let $ a_t   \in   \{0,1\}$ denote the action taken at the beginning of slot $t$, where 
$a_t = 1$ if the server stays with the current queue and $a_t = 0$ if it
switches to the other queue. One packet is successfully
received
from queue $i$ at time slot $t$, if $m(t) = i,a_t  =  1$ and $C_i(t)  =  1$.

\begin{definition}[Stability \cite{Neely10},\hspace{-0.5mm} \cite{neely05}]\label{def:stab}
\emph{The system is 
stable if}
\begin{equation*}
\displaystyle \limsup_{t\rightarrow \infty} \frac{1}{t}
\sum_{\tau=0}^{t-1} \sum_{i \in \{1,2\}}\mathbb{E}[Q_i(\tau)] < \infty.
\end{equation*}
\end{definition}
For the case of integer valued arrival processes, this stability criterion 
implies the existence of a long
run stationary distribution for the queue sizes 
with bounded first moments \cite{Neely10}. 

\vspace{0.1cm}
\begin{definition}[Stability Region \cite{Neely10},\hspace{-0.5mm} \cite{neely05}]\label{def:throughput_region}
\emph{The stability region $\mathbf{\Lambda}$ is the set of all
arrival rate vectors $\boldsymbol \lambda=(\lambda_1,\lambda_2)$ such that there exists
a control algorithm that stabilizes the system.}
\end{definition}
\vspace{0.1cm}
The $\mathbf{\delta}$-stripped stability region is defined for some
$\delta >0$ as $\mathbf{\Lambda}^{\delta} \doteq \Big\{
(\lambda_1,\lambda_2)| (\lambda_1+\delta, \lambda_2+\delta) \in
\mathbf{\Lambda} \Big\}.$ A policy is said to achieve
$\gamma$-fraction of $\mathbf{\Lambda}$, if it stabilizes the system
for all input rates inside $\gamma \mathbf{\Lambda}$, 
where $\gamma=1$ for a throughput-optimal policy. 

%
%
%

%
%

In the following, we start by explicitly characterizing the stability region for both memoryless channels and channels with memory and
show that channel memory can be exploited to enlarge the
stability region 
significantly.

\subsection{Motivation: Channels Without Memory}\label{Sec:Motiv} 

In this section we assume 
that $\epsilon=0.5$ so that the channel processes are i.i.d. over time. The stability region of the corresponding system with no-switchover time 
was established in \cite{tass93}:
$\lambda_1,\lambda_2 \in [0,0.5]$ 
and
$\lambda_1 + \lambda_2 \le 0.75$. Note that when the switchover time is zero, the stability region is the same
for both i.i.d. and Markovian channels, which is a special case of
the results in \cite{neely05}. However, when the
switchover time is non-zero, the stability region is reduced
considerably: 
%
%

\vspace{0.5mm}
\begin{theorem}\label{thm:stab_iid}
\emph{The stability region of the system with i.i.d. channels and
one-slot switchover delay is given by,
\begin{equation}\label{eq:Stab_iid}
\mathbf{\Lambda}=\{ (\lambda_1,\lambda_2) \big | \lambda_1 +
\lambda_2 \le 0.5, \lambda_1,\lambda_2 \ge 0 \}.
\end{equation}
In addition, the simple Exhaustive (Gated) policy is
throughput-optimal.}
\end{theorem}
The proof is given in Appendix A for a more general system. 
%
%
The basic idea behind the proof is that as soon as the
server switches to queue $i$ under some policy, 
the time to the ON state is a geometric random variable with mean 2 slots, independent of the policy.
Therefore, a necessary condition for stability is given by the
stability condition for a system without switchover times 
and i.i.d. service times with geometric distribution of mean 2 slots as given by (\ref{eq:Stab_iid}). 
The fact that the simple Gated policy is throughput-optimal 
follows from the observation that as the arrival rates are close to
the boundary of the stability region, the fraction of time the
server spends receiving packets dominates the fraction of time spent
on switching \cite{vish06}.

\begin{figure}
\centering \psfrag{1\r}[l][][0.95]{$\!\!\!\!0.5$}
\psfrag{2\r}[l][][0.95]{$\!\!\!\!\!\!0.5$}
\psfrag{3\r}[l][][1]{$\!\!\lambda_2$} \psfrag{4\r}[l][][1]{$\!\!\lambda_1$}
\psfrag{6\r}[l][][0.85]{$\!\!\!\!\!\!\!\!\!\!\!\!\!\!\!\!\!$\textrm{no-switchover
time}} \psfrag{9\r}[l][][0.85]{\textrm{i.i.d.}}
\psfrag{5\r}[l][][0.85]{$\!\!\!\!\!$\textrm{channels}}
\psfrag{7\r}[l][][.85]{$\!$\textrm{Markovian}}
\psfrag{8\r}[l][][.85]{$\!$\textrm{channels}}
\includegraphics[width=0.23\textwidth]{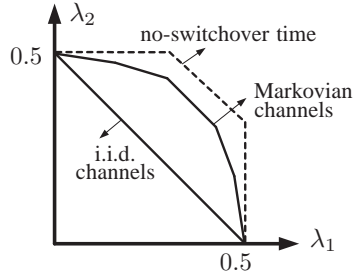}
\captionsetup{font=footnotesize}\caption{Stability
region under memoryless (i.i.d.) channels and channels with memory (Markovian with $\epsilon < 0.5$) with and without switchover delay.}\label{Fig:stab_region_two_links}
\end{figure}

As depicted in Fig.~\ref{Fig:stab_region_two_links}, 
the stability region of the system is considerably reduced for nonzero
switchover delay. 
Note that 
for systems in which channels are always connected, the stability
region is given by $\lambda_1+\lambda_2 \le 1, \lambda_1,\lambda_2 \ge
0$ and is \emph{not} affected by the switchover delay \cite{vish06}. 
%
Therefore, \emph{it is the combination of switchover delay and
random connectivity that result in fundamental changes in system
stability}.
%
%

\vspace{2mm}
\begin{remark}\label{rem:iid_stab}
As shown in Appendix A, the results in this subsection can easily be generalized to the case of non-symmetric Gilbert-Elliot channels with arbitrary switchover delays. \emph{For a system of $2$ queues with 
arbitrary switchover delays 
and i.i.d. channels with probabilities
$p_i, i \in \{1,2 \}$, $\mathbf{\Lambda}$ is the set of all $\boldsymbol \lambda \boldsymbol\ge0$ such that
$\lambda_1/p_1 + \lambda_2/p_2 \le 1$.
Moreover, simple Exhaustive (Gated) policy is throughput-optimal}. 
\end{remark}

When channel processes have memory, it is clear that one can achieve better
throughput region than the i.i.d. channels case if the channels are
positively correlated over time. This is because we can exploit the
channel diversity when the channel states stay the same with high
probability. 
In the following, we show that indeed the throughput region
approaches the throughput region of no switchover time case in in
\cite{tass93} as the channels become more correlated over time. Note
that the throughput region in \cite{tass93} is the same for both
i.i.d. and Markovian channels under the condition that probability
of ON state for the $i.i.d.$ channels is the same as the steady
state probability of ON state for the two state Markovian channels.
This fact can be derived as a special case of the seminal work of
Neely in \cite{neely05}.

\subsection{Channels With Memory - Stability Region}\label{Sec:Stab_Reg}

When switchover times are non-zero, the memory in the channel can be exploited to improve
the stability region considerably.
Moreover, as $\epsilon \rightarrow 0$, the
stability region tends to that achieved by the system with
no-switchover time and for $0<\epsilon<0.5$ it lies between the
stability regions corresponding to the two extreme cases $\epsilon
=0.5$ and $\epsilon \rightarrow 0$ as  shown in
Fig.~\ref{Fig:stab_region_two_links}.
%
%

We start by analyzing the corresponding system with saturated
queues, i.e., both queues are always non-empty. Let
$\mathbf{\Lambda}_s$ denote the set of all time average expected
departure rates that can be obtained from the two queues in the
saturated system under all possible policies that are possibly
history dependent, randomized and non-stationary. We will show that
$\mathbf{\Lambda}=\mathbf{\Lambda}_s$. We prove the necessary
stability conditions 
in the following Lemma and establish sufficiency in the next section.    

\vspace{0.1cm}
\begin{lemma}\label{lem:Lam_under_Lam_s}
\emph{We have that}
\begin{equation*}
\mathbf{\Lambda} \subseteq \mathbf{\Lambda}_s.
\end{equation*}
\end{lemma}
\begin{proof}
Given a policy $\pi$ for the dynamic queueing system specifying the switch
and stay actions based possibly on observed channel and queue state
information, consider the saturated system with \emph{the same
sample path of channel realizations} for $t\in \{0,1,2,... \}$ and
\emph{the same set of actions} as policy $\pi$ at each time slot
$t\in \{0,1,2,... \}$. Let this policy for the saturated system be
$\pi'$, and note that the policy $\pi'$ can be non-stationary\footnote{The policy $\pi'$ can be based on keeping virtual queue length information, i.e., the queue lenghts in the dynamic queueing system.} 
Let $D_i(t),i \in \{1,2\}$ be total number departures by time $t$ from queue-$i$ in the original system under policy $\pi$
and let $D_i'(t),i \in \{1,2\}$ be the corresponding quantity for
the saturated system under policy $\pi'$. It is clear that $\lim_{t
\rightarrow \infty} (D_1(t)+D_2(t))/t \le 1$, where the same
statement also holds for the limit of $D_i'(t),i \in \{1,2\}$. Since
some of the ON channel states are wasted in the original system due
to empty queues, we have
\begin{equation}\label{eq:sat_dept_bound}
D_1(t) \le D_1'(t),\;\;\;\textrm{and},\;\;\;D_2(t)\le D_2'(t).
\end{equation}
Therefore, the time average expectation of $D_i(t),i \in \{1,2\}$ is
also less than or equal to the time average expectation of
$D_i'(t),i \in \{1,2\}$. This completes the proof since
(\ref{eq:sat_dept_bound}) holds under any policy $\pi$ for the
original system.
\end{proof}

Now, we establish the region $\mathbf{\Lambda}_s$ by formulating the
system dynamics as a Markov Decision Process (MDP). Let
$\mathbf{s}_t=(m(t),C_1(t),C_2(t)) \in S$ denote the system state at
time $t$ where $\mathcal{S}$ is the set of all states. Also, let $a_t \in
\mathcal{A}=\{0,1\}$ denote the action taken at time slot $t$ where $\mathcal{A}$ is the
set of all actions at each state. 
Let $\mathbb{H}(t)=[\mathbf{s}_{\tau}]|_{\tau=0}^{t} \cup [a_{\tau}]|_{\tau=0}^{t-1}$ denote the full
history of the system state until time $t$ and let $\Upsilon(\mathcal{A})$ denote
the set of all probability distributions on $\mathcal{A}$.
For the saturated
system,
a policy is a mapping from the set of all possible past histories to $\Upsilon(\mathcal{A})$ \cite{AltSch91}, \cite{shie05}. 
A policy is said to be \emph{stationary} if, given a particular state, it applies the same decision rule in all stages and under a stationary policy, the process $\{\mathbf{s}_t;t \in \mathbb{N}\cup \{0\} \}$ forms a Markov chain.
%
In each time slot
$t$, the server observes the current state $\mathbf{s}_t$ and
chooses an action $a_t$. Then the next state $j$ is realized
according to the transition probabilities $\mathbf{P}(j|\mathbf{s},a)$, which
depend on the random channel processes. Now, we define the reward
functions as follows: \beqn
\!\!\!\!\!\!\overline{r}_1(\mathbf{s},a) \!\!\!\!\!&\doteq& \!\!\!\!\!1 \textrm{ if } \mathbf{s}=(1,1,1) \textrm{ or } \mathbf{s}=(1,1,0) \textrm{, and } a\!=\!1 \label{eq:dept_rate1} \\
\!\!\!\!\!\!\overline{r}_2(\mathbf{s},a) \!\!\!\!\!&\doteq&\!\!\!\!\! 1 \textrm{ if }
\mathbf{s}=(2,1,1) \textrm{ or } \mathbf{s}=(2,0,1) \textrm{, and
} a\!=\!1, \label{eq:dept_rate} \eeqn and
$\overline{r}_1(\mathbf{s},a)=\overline{r}_2(\mathbf{s},a)\doteq0$ otherwise. That is, a reward is
obtained when the server stays at an ON channel. We are interested
in the
set of all possible time average expected departure rates, 
therefore, given some
$\alpha_1,\alpha_2  \ge 0$, 
$\alpha_1+\alpha_2=1$,
we define the system reward at time $t$ by
$\overline{r}(\mathbf{s},a) \doteq \alpha_1\overline{r}_1(\mathbf{s},a) + \alpha_2
\overline{r}_2(\mathbf{s},a)$. 
The average reward of policy $\boldsymbol{\pi}$ is defined as follows:
\begin{equation*}
r^{\boldsymbol{\pi}} \doteq \displaystyle \limsup_{K\rightarrow \infty} \frac{1}{K} E
\Big \{ \sum_{t=1}^K \overline{r}(\mathbf{s}_t,a_t^{\boldsymbol{\pi}}) \Big \}.
\end{equation*}
Given some $\alpha_1,\alpha_2 \ge 0$, we are interested in the policy that
achieves the maximum time average expected reward
$r^{*} \doteq \max_{\boldsymbol{\pi}} r^{\boldsymbol{\pi}}$. 
This optimization problem is a discrete time 
MDP characterized by the state transition probabilities
$\mathbf{P}(j|\mathbf{s},a)$ with 8 states and 2 actions per state.
Furthermore, any given pair of states are accessible from each other (i.e., there exists a positive probability path between the states) under some stationary-deterministic policy. Therefore this MDP belongs to the class of \emph{Weakly Communicating}
MDPs\footnote{In fact, other than the trivial suboptimal policy $\boldsymbol{\pi}_s$ that decides to stay with the current queue in all states, all stationary deterministic policies are unichain, namely, they have a single recurrent class regardless of the initial state. Hence, when $\boldsymbol{\pi}_s$ is excluded, we have a Unichain MDP.}\cite{puterman05}.

\subsubsection{The State-Action Frequency Approach}

For Weakly Communicating MDPs with finite state and action spaces and bounded rewards, there exists
an optimal stationary-deterministic policy, given as a solution to standard Bellman's equation, with optimal average reward independent of the initial state
\cite[Theorem 8.4.5]{puterman05}.
%
This is because if a stationary policy has a nonconstant gain over initial states, one can construct another stationary policy with constant gain which dominates the former policy, which is possible since there exists a positive probability path between any two recurrent states under some stationary policy \cite{shie05}. 
%
The \emph{state-action frequency} approach, or the \emph{Dual Linear Program (LP)} approach, given below provides a systematic and intuitive framework to solve such average cost MDPs, and it can be derived using Bellman's equation and the monotonicity property of Dynamic Programs
[Section 8.8]\cite{puterman05}: 
%
%
%
\begin{equation}\label{eq:dual_LP_def}
\mathrm{Maximize}\;\;\;\; \sum_{\mathbf{s} \in \mathcal{S}} \sum_{a \in \mathcal{A}} \overline{r}(\mathbf{s},a)\mathbf{x}(\mathbf{s},a)
\end{equation}
subject to the balance
equations
\beqn \mathbf{x}(\mathbf{s};1) + \mathbf{x}(\mathbf{s};0) = \!\!\sum_{\mathbf{s}'\in\mathcal{S}}
\sum_{a\in\mathcal{A}}\!\mathbf{P}\big(\mathbf{s}|\mathbf{s}',a\big)\mathbf{x}(\mathbf{s}',a) ,
\;\forall\; \mathbf{s} \in S, \label{eq:saf_balance} \eeqn the normalization
condition \beqn \sum_{\mathbf{s} \in \mathcal{S}}\mathbf{x}(\mathbf{s};1) + \mathbf{x}(\mathbf{s};0)=1,  \label{eq:saf_norm}
\eeqn and the nonnegativity constraints \beqn \mathbf{x}(\mathbf{s},a)\ge 0,  \:
\mathbf{s} \in \mathcal{S}, a\in \mathcal{A}. \label{eq:saf_nonneg}
\eeqn
The feasible region of this LP constitutes a polytope called the \emph{state-action polytope} $\mathbf{X}$ and the elements of this polytope $\mathbf{x} \in \mathbf{X}$ are called state-action frequency vectors.
Clearly, $\mathbf{X}$ is a convex, bounded and closed set. 
Note that
$\mathbf{x}(\mathbf{s};1)$
can be interpreted as the stationary
probability that action
\emph{stay} 
is taken at state $\mathbf{s}$. 
More precisely, a point $\mathbf{x}\in\mathbf{X}$ 
corresponds to a stationary randomized policy 
that takes action $a \in \{0,1\}$ 
at state $\mathbf{s}$ w.p.
\begin{equation}\label{eq:rand_stat}
\,\!\mathbf{P}(\textrm{\emph{action $a$} at state$\,\mathbf{s}$})\!=\!
\frac{\mathbf{x}(\mathbf{s},a)}{\mathbf{x}(\mathbf{s};1)+\mathbf{x}(\mathbf{s};0)}, a \in A,\mathbf{s}\in S_x, \!\!\!\!\!\!\!\!\!
\end{equation}
where $S_x$ is the set of 
recurrent states  given by
$S_x\equiv \displaystyle \{ \mathbf{s}\in S: \mathbf{x}(\mathbf{s};1)+\mathbf{x}(\mathbf{s};0)>0\},$
and actions are arbitrary for transient states $\mathbf{s} \in S/S_x$ 
\cite{shie05}, \cite{puterman05}.

%

Next we argue that the empirical state-action frequencies corresponding to any given policy (possibly randomized, non-stationary, or non-Markovian) lies in the state-action polytope $\mathbf{X}$. This ensures us that the optimal solution to the dual LP in (\ref{eq:dual_LP_def}) is over possibly non-stationary and history-dependent policies. In the following we give the precise definition and the properties of the set of empirical state-action frequencies.
We define the \emph{empirical} state-action frequencies $\hat{x}^T(\mathbf{s},a)$ as
\begin{equation}\label{eq:empr_saf}
\displaystyle \hat{x}^T(\mathbf{s},a) \doteq \frac{1}{T}\sum_{t=1}^T I_{\{\mathbf{s}_t=\mathbf{s},a_t=a\}},
\end{equation}
where $I_E$ is the indicator function of an event $E$, i.e., $I_E=1$ if $E$ occurs and $I_E=0$ otherwise. 
Given a policy $\boldsymbol{\pi}$, let $P^{\boldsymbol{\pi}}$ be the state-transition probabilities under the policy $\boldsymbol{\pi}$ and $\boldsymbol{\phi} = (\phi_s)$ an initial state distribution with $\sum_{\mathbf{s} \in \mathcal{S}} \phi_{\mathbf{s}}=1$. We let $x^T_{\boldsymbol{\pi},\phi}(\mathbf{s},a)$ be the \emph{expected empirical} state-action frequencies under policy $\boldsymbol{\pi}$ and initial state distribution $\phi$:
\begin{eqnarray*}
\displaystyle x^T_{\boldsymbol{\pi},\phi}(\mathbf{s},a) \!\!\!&\doteq& \!\!\! \mathbb{E}^{\boldsymbol{\pi},\phi}\left[ \hat{x}^T(\mathbf{s},a) \right] \nonumber\\
&=& \!\!\!\frac{1}{T}\sum_{t=1}^T \sum_{\mathbf{s}'\in \mathcal{S}} \!\!\phi_{\mathbf{s}'} P^{\boldsymbol{\pi}}\!(\mathbf{s}_t=\mathbf{s},a_t=a|\mathbf{s}_0=\mathbf{s}' ).
\end{eqnarray*}
We let $\mathbf{x}_{\boldsymbol{\pi},\phi} \in \Upsilon(\mathcal{S} \times \mathcal{A})$ (as in \cite{shie05}, \cite{puterman05}) be the limiting expected state-action frequency
vector, if it exists, starting from an initial state distribution $\phi$, under a general policy $\boldsymbol{\pi}$ (possibly randomized, non-stationary, or non-Markovian): 
\begin{equation}\label{eq:exp_saf}
x_{\boldsymbol{\pi},\phi}(\mathbf{s},a) = \lim_{T\rightarrow \infty} x^T_{\boldsymbol{\pi},\phi}(\mathbf{s},a).
\end{equation}
Let the set of all limit points be defined by 
\begin{eqnarray*}
\mathbf{X}^{\phi}_{\Pi} &\doteq \{x \in \Upsilon(\mathcal{S} \times \mathcal{A}): \textrm{ there exists a policy $\boldsymbol{\pi}$ s.t.} \\
&\textrm{the limit in (\ref{eq:exp_saf}) exists and $\mathbf{x}=\mathbf{x}_{\boldsymbol{\pi},\phi}$ } \}.
\end{eqnarray*}
Similarly let $X^{\phi}_{\Pi'}$ denote the set of all limit points of a particular class of policies $\Pi'$, starting from an initial state distribution $\phi$. We let $\Pi_{SD}$ denote the set of 
all stationary-deterministic policies and we let $co(\Xi)$ denote the closed convex hull of set $\Xi$. The following theorem establishes the equivalency between the set of all achievable limiting state-action frequencies and the state-action polytope:
\vspace{2mm}
\begin{theorem} \cite[Theorem 8.9.3]{puterman05}, \cite[Theorem 3.1]{shie05}. \label{thm:empr_saf_wc}
\emph{For any initial state distribution $\phi$}
\begin{equation*}
co(\mathbf{X}^{\phi}_{\Pi_{SD}}) = \mathbf{X}^{\phi}_{\Pi} = \mathbf{X}.
\end{equation*}
\end{theorem}
We have $co(\mathbf{X}^{\phi}_{\Pi_{SD}}) \subseteq \mathbf{X}^{\phi}_{\Pi}$
since convex combinations of vectors in $\mathbf{X}^{\phi}_{\Pi_{SD}}$ correspond to limiting expected state-action frequencies for stationary-randomized policies, which can also be obtained by time-sharing between stationary-deterministic policies.
The inverse relation $co(\mathbf{X}^{\phi}_{\Pi_{SD}}) \supseteq \mathbf{X}^{\phi}_{\Pi}$ holds since for weakly communicating MDPs, there exists a stationary-deterministic optimal policy independent of the initial state distribution. Next, for any stationary-deterministic policy, the underlying Markov chain is stationary and therefore the limits $\mathbf{x}_{\boldsymbol{\pi},\phi}$ exists and satisfies the constraints (\ref{eq:saf_balance}), (\ref{eq:saf_norm}) and (\ref{eq:saf_nonneg}) of the polytope $\mathbf{X}$. Using $co(\mathbf{X}^{\phi}_{\Pi_{SD}}) = \mathbf{X}^{\phi}_{\Pi}$ and the convexity of $\mathbf{X}$ establishes
$\mathbf{X}^{\phi}_{\Pi} \subseteq \mathbf{X}$. Furthermore, via (\ref{eq:rand_stat}), every $\mathbf{x} \in \mathbf{X}$ corresponds to a stationary-randomized policy for which the limits $\mathbf{x}_{\boldsymbol{\pi},\phi}$ exists, establishing
$\mathbf{X}^{\phi}_{\Pi} \supseteq \mathbf{X}$.

%
Letting $\textrm{ext}(\mathbf{X})$ denote the set of extreme (corner) points of $\mathbf{X}$, an immediate corollary to Theorem \ref{thm:empr_saf_wc} is as follows:
\vspace{1mm}
\begin{corollary}\cite{shie05},\!\cite{puterman05}.\label{cor:X_vortex}
\emph{For any initial state distribution $\phi$}
\begin{equation*}
\textrm{ext}(\mathbf{X})=\mathbf{X}^{\phi}_{\Pi_{SD}}. 
\end{equation*}
\end{corollary}
\vspace{0.1cm}
The intuition behind this corollary is that if $\mathbf{x}$ is a corner point of $\mathbf{X}$, it cannot be expressed as a convex combination of any two
other elements in $\mathbf{X}$, therefore, for each state $\mathbf{s}$ only one action has a nonzero probability. 

Finally, we have that under any policy the probability of a large distance between the empirical expected state-action
frequency vectors and the state-action polytope $\mathbf{X}$ decays exponentially fast in time. This result is similar to
the mixing time of an underlying Markov chain to its steady state and we utilize such convergence results within the Lyapunov drift analysis for the dynamic queuing system in Section \ref{Sec:Opt_Pol}.

\subsubsection{The Rate Polytope $\boldsymbol{\Lambda}_s$}

Using the theory on state-action polytopes in the previous section, we characterize the set of all achievable time-average expected rates in the saturated system, $\boldsymbol{\Lambda}_s$.
The following linear transformation of the state-action polytope
$\mathbf{X}$ defines 
the 2 dimensional \emph{rate
polytope}\cite{shie05}:
%
%
\begin{eqnarray*}
\mathbf{\Lambda}_s=
\Big\{(r_1,r_2)\big| &r_1&=
\sum_{\mathbf{s}\in \mathcal{S}}\sum_{a\in\mathcal{A}}\mathbf{x}(\mathbf{s},a)\overline{r}_1(\mathbf{s},a)\nonumber \\
&r_2&=
\sum_{\mathbf{s}\in \mathcal{S}}\sum_{a\in\mathcal{A}}\mathbf{x}(\mathbf{s},a)\overline{r}_2(\mathbf{s},a),\;
\mathbf{x} \in \mathbf{X} \Big\},
\end{eqnarray*}
where
$\overline{r}_1(\mathbf{s},a)$ and $\overline{r}_2(\mathbf{s},a)$ are the 
reward functions defined in (\ref{eq:dept_rate1}) and
(\ref{eq:dept_rate}). This polytope is the set of all time average
expected departure rate pairs that can be obtained in the saturated
system, i.e., it is the rate region $\mathbf{\Lambda}_s$. An
explicit way of characterizing
$\mathbf{\Lambda}_s$ is given in Algorithm \ref{alg:LP_for_Lambda_s}. 

\begin{algorithm}[h]
\caption{\emph{Stability Region Characterization}}
\label{alg:LP_for_Lambda_s}
\begin{algorithmic}[1]

\STATE Given $\alpha_1,\alpha_2\!\! \ge\! 0$, $\alpha_1+\alpha_2=1$ solve the following Linear
Program (LP)
\begin{eqnarray}
\max.  \quad  
\alpha_1 r_1(\mathbf{x}) +  \alpha_2r_2(\mathbf{x})  \nonumber\\
\mbox{subject to} \quad  \mathbf{x} \in \mathbf{X}. \label{eq:LP}
\end{eqnarray}

\STATE For a given $\alpha_2/\alpha_1$ ratio, there exists an optimal solution
$(r_1^{*}, r_2^{*})$ of the LP in
(\ref{eq:LP}) at a corner point of 
$\mathbf{\Lambda_s}$. Find all possible corner points and take their
convex combination.
\end{algorithmic}
\end{algorithm}
The fundamental theorem of Linear Programming guarantees that an optimal solution of the LP in (\ref{eq:LP}) lies at a corner (extreme) point of the polytope $\mathbf{X}$ \cite{BertTsit97}.
%
%
%
%
%
Furthermore, the one-to-one correspondence between the extreme points of the polytope $\mathbf{X}$ and stationary-deterministic polices stated in Corollary~\ref{cor:X_vortex} is useful for finding the solutions of the above LP for all possible $\alpha_2/\alpha_1$ ratios. 
%
Namely, there are a total of $2^8$ stationary-deterministic
policies since we have $8$ states and $2$ actions per state and
finding the rate pairs corresponding to these 256 stationary-deterministic
policies and taking their convex combination gives
$\mathbf{\Lambda}_s$. Fortunately, we do not have to go through this
tedious procedure. The fact that at a vertex of (\ref{eq:LP}) either
$\mathbf{x}(\mathbf{s};1)$ or $\mathbf{x}(\mathbf{s};0)$ has to be zero for each $\mathbf{s} \in
\mathbf{S}$ provides a useful guideline for analytically solving
this LP. 
%
The following theorem 
characterizes the stability region explicitly. It shows that the stability region enlarges as the channel has more memory and that 
there is a critical value of the channel correlation parameter given by $\epsilon_c \doteq 1-\sqrt{2}/2$ at which the structure of the stability region changes.
\begin{theorem}\label{thm:stab}
\emph{The rate region $\mathbf{\Lambda}_s$ is the set of all
rates $r_1\ge 0$, $r_2 \ge 0$ that for $\epsilon <
\epsilon_c $ satisfy} \beqn
\epsilon r_1 + (1-\epsilon)^2r_2 &\le& \frac{(1-\epsilon)^2}{2}\nonumber\\
(1-\epsilon)r_1 + (1+\epsilon-\epsilon^2)r_2  &\le& \frac{3}{4}-\frac{\epsilon}{2}\nonumber\\
r_1 + r_2  &\le& \frac{3}{4}-\frac{\epsilon}{2}\nonumber\\
(1+\epsilon-\epsilon^2)r_1 + (1-\epsilon)r_2  &\le& \frac{3}{4}-\frac{\epsilon}{2}\nonumber\\
(1-\epsilon)^2r_1 + \epsilon r_2 &\le&
\frac{(1-\epsilon)^2}{2},\nonumber 
\eeqn
 \emph{and for $\epsilon \ge \epsilon_c $ satisfy}
\beqn
r_1 + (1-\epsilon)(3-2\epsilon)r_2 &\le& \frac{(1-\epsilon)(3-2\epsilon)}{2}\nonumber\\
r_1 +r_2  &\le& \frac{3}{4}-\frac{\epsilon}{2}\nonumber\\
(1-\epsilon)(3-2\epsilon)r_1 + r_2  &\le&
\frac{(1-\epsilon)(3-2\epsilon)}{2}.\nonumber
 \eeqn
\end{theorem}
The proof of the theorem is given in Appendix B and it 
is based on solving the LP in (\ref{eq:LP}) for all weights $\alpha_1$ and $\alpha_2$ 
to find the corner points of $\mathbf{\Lambda_s}$,
and then applying
Algorithm~\ref{alg:LP_for_Lambda_s}. The following observation follows from Theorem~\ref{thm:stab}.
\vspace{0.1cm}
\begin{observation}\label{obs:sum_tp_2queues}
\emph{The maximum achievable sum-rate in the saturated system is given by}
\begin{equation*}
r_1 + r_2 = \frac{3}{4} - \frac{\epsilon}{2}.
\end{equation*}
\end{observation}
\vspace{0.1cm}
%
Note that $r_1 + r_2 \le \frac{3}{4}$ is the boundary of the stability region for the system without switchover delay analyzed in \cite{tass93}, where the probability that at least 1 channel is in ON state is $3/4$. Therefore, $\epsilon/2$ is the \emph{throughput loss due to the 1 slot switchover delay}. 
This throughput loss corresponds to the probability that the server is at a queue with an OFF state when the other queue is in an ON state.

\begin{figure}
\begin{center}
\includegraphics[width=0.45\textwidth]{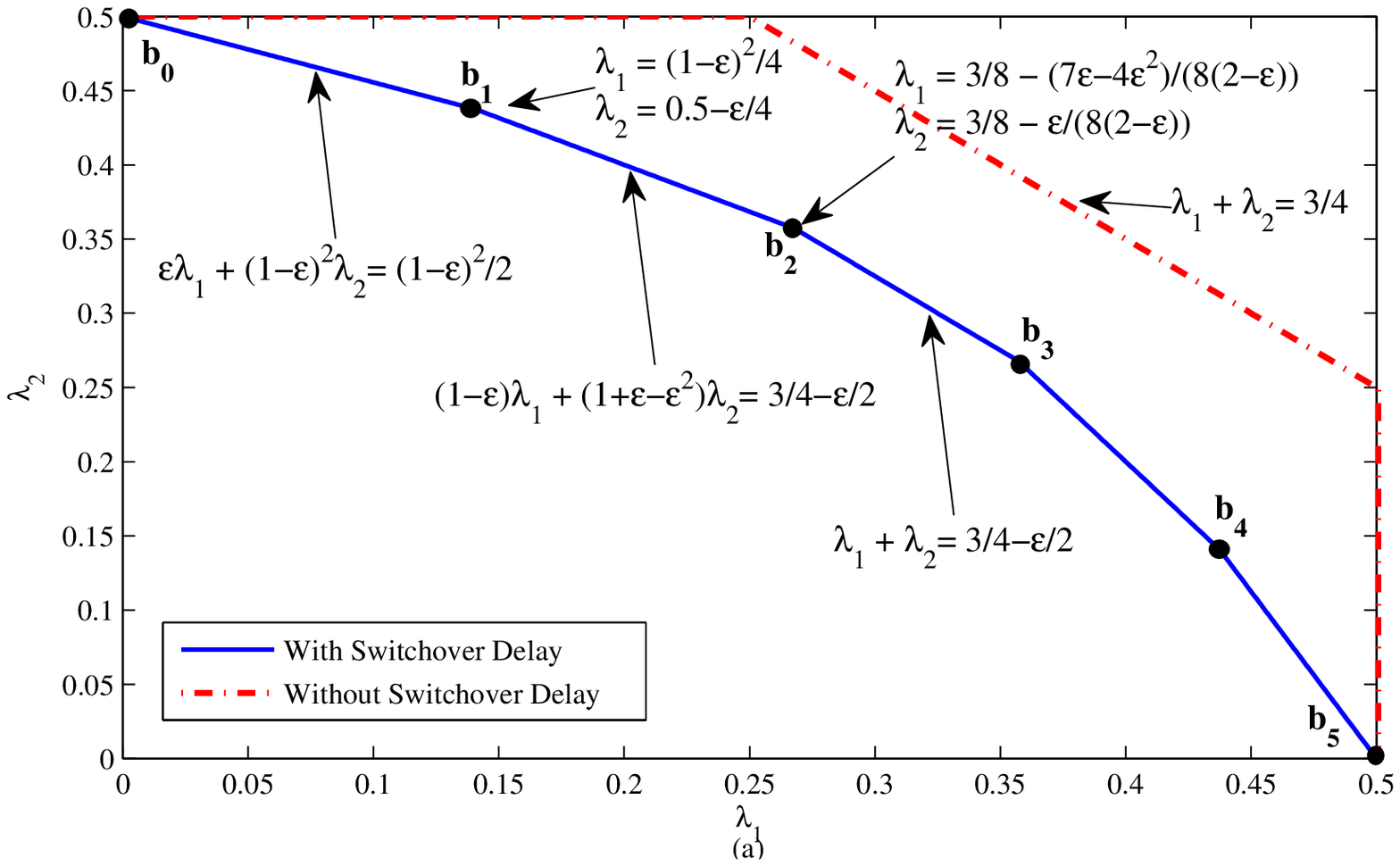}
\includegraphics[width=0.45\textwidth]{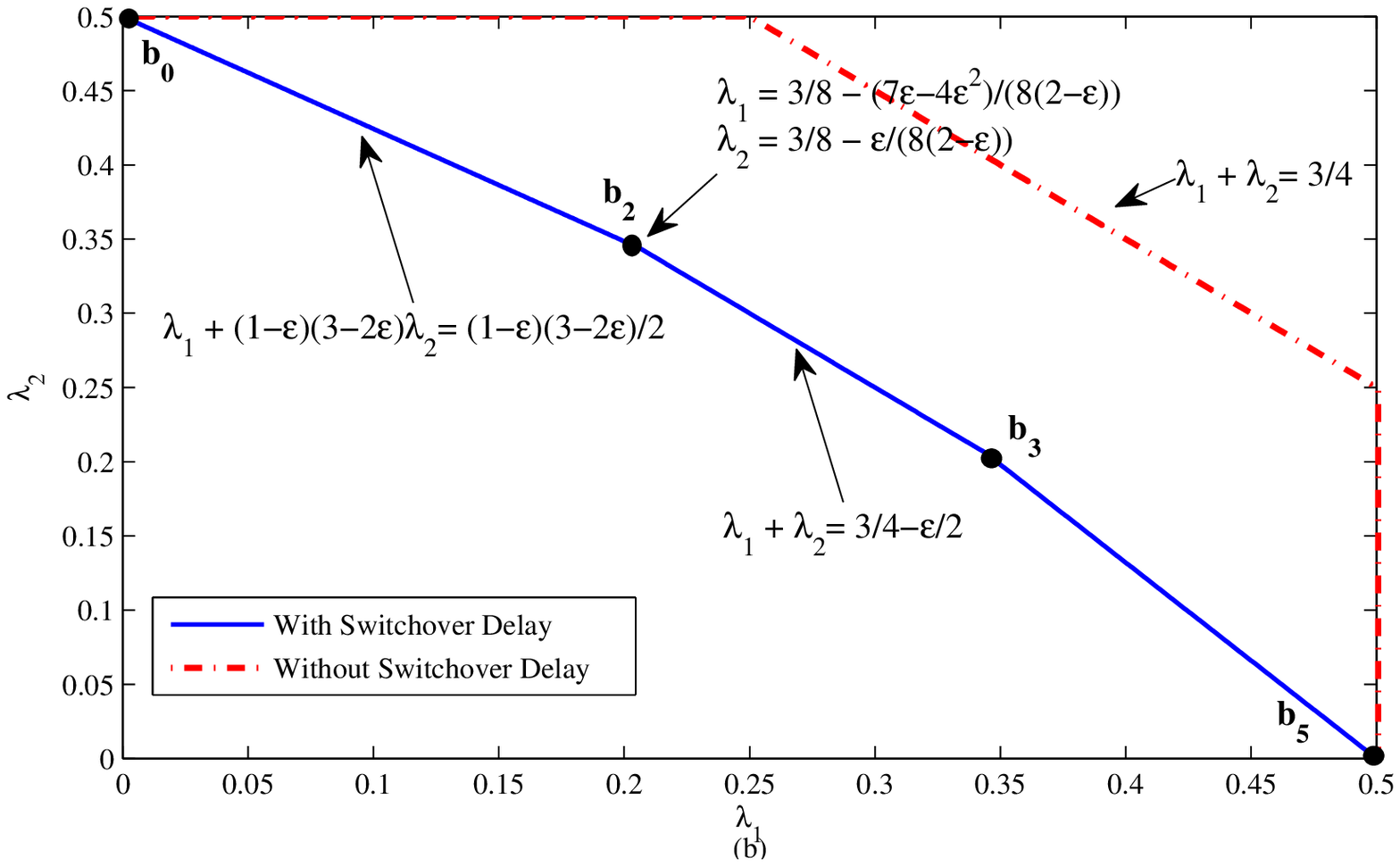}
\end{center}
\captionsetup{font=footnotesize}\caption{Stability
region under channels with memory, with and without switchover delay
for (a) $\epsilon=0.25 < \epsilon_c$ and (b) $\epsilon =0.40 \ge
\epsilon_c$.} \label{Fig:throughput_eps_0_250_write_up}
\end{figure}

The stability regions for the two ranges of $\epsilon$ are
displayed in Fig.~\ref{Fig:throughput_eps_0_250_write_up}~(a) and
(b). As $\epsilon \rightarrow 0.5$, the stability region converges
to that of the i.i.d. channels with ON probability equal to $0.5$.
In this regime, knowledge of the current channel state is of no
value.
%
%
%
%
As $\epsilon\rightarrow 0$ the stability region converges to that
for the system with no-switchover time in \cite{tass93}. In this
regime, the channels are likely to stay the same for many
consecutive time slots, therefore, the effect of switching delay is
negligible.

The rate region $\mathbf{\Lambda}_s$ for the case of non-symmetric Gilbert-Elliot channels is given in Appendix B.

\vspace{1mm}
\begin{remark}
The stability region
characterization in terms of state-action frequencies 
is general. For instance, \emph{this technique can be used to establish the
stability regions of systems with more than two queues, arbitrary
switchover times, and more complicated Markovian channel processes.} Of course, explicit characterization as in Theorem~\ref{thm:stab} may not always be possible.
\end{remark}
%

\subsection{Frame Based Dynamic Control (FBDC) Policy}\label{Sec:Opt_Pol}

We propose a frame-based dynamic control (FBDC) policy inspired by
the characterization of the stability region in terms of
state-action frequencies and prove that it is throughput-optimal
asymptotically in the frame length. The motivation behind the FBDC
policy is that a policy $\boldsymbol{\pi}^*$ that achieves the optimization in
(\ref{eq:LP}) for given weights $\alpha_1$ and $\alpha_2$ for the
saturated system should achieve a \emph{good} performance in
the original system when
the queue sizes $Q_1$ and $Q_2$ are used as weights. 
This is because first, the policy $\boldsymbol{\pi}^*$ will lead to similar
average departure rates in both systems for sufficiently high
queue sizes, and second, the usage of queue sizes as weights
creates self adjusting policies that capture the dynamic changes due
to stochastic arrivals 
similar to Max-Weight scheduling in
\cite{tass92}. Specifically, we divide the time into equal-size
intervals of $T$ slots and let $Q_1(jT)$ and $Q_2(jT)$ be the queue
lengths at the beginning of the $j$th interval. We find the
deterministic policy that optimally solves (\ref{eq:LP}) when
$Q_1(jT)$ and $Q_2(jT)$ are used as weights 
and then apply this policy in each time slot of the
frame. The FBDC policy is described in Algorithm \ref{alg:FBDC}
in details.

\begin{algorithm}[h]
\caption{\textsc{Frame Based Dynamic Control (FBDC) Policy}}
\label{alg:FBDC}
\begin{algorithmic}[1]

\STATE Find the policy $\boldsymbol{\pi}^*$ that optimally solves the following LP
\beqn
\displaystyle \mbox{max.}_{\{r_1, r_2\}} &  Q_1(jT) r_1 + Q_2(jT) r_2 \nonumber \\
\mbox{subject to}  &   \left( r_1, r_2\right) \in \mathbf{\Lambda}_s
\label{objcon}
\eeqn where $\mathbf{\Lambda}_s$ is the rate polytope derived in
Section \ref{Sec:Stab_Reg}.
\STATE  Apply $\boldsymbol{\pi}^*$ in each time slot of the frame.
\end{algorithmic}
\end{algorithm}
%
There exists an optimal solution $(r_1^{*},
r_2^{*})$ of the LP in (\ref{objcon}) that is a corner point of $\mathbf{\Lambda}_s$ \cite{BertTsit97} and the policy $\boldsymbol{\pi}^*$ that corresponds to this point is a stationary-deterministic policy by Corollary~\ref{cor:X_vortex}.
\vspace{0.1cm}
\begin{theorem}\label{thm:FBDC}
\emph{For any $\delta>0$, there exists a large enough frame length $T$ such that the FBDC policy stabilizes the system for all arrival
rates 
within the
$\mathbf{\delta}$-stripped stability region
$\mathbf{\Lambda}_s^{\delta} = \mathbf{\Lambda}_s - \delta\mathbf{1}$.} 
\end{theorem}
An immediate corollary to this theorem is as follows:
\begin{corollary}
\emph{The FBDC policy is throughput-optimal asymptotically in the frame length.}
\end{corollary}
\vspace{0.5mm}
The proof of Theorem~\ref{thm:FBDC} is given in Appendix D. It performs a drift
analysis using the standard quadratic Lyapunov function. However, \emph{it
is novel in utilizing the state-action frequency framework of MDP theory
within the Lyapunov drift arguments}. The basic idea is that, for sufficiently large queue lengths, when the
optimal policy solving (\ref{objcon}), $\boldsymbol{\pi}^*$,
is applied over
a sufficiently long frame of $T$ slots, the average 
output rates of both the actual system and the corresponding saturated
system converge to $\mathbf{r}^*$. 
For the saturated system, the probability of a large difference between empirical and steady state rates
decreases exponentially fast in $T$
\cite{shie05}, similar to the convergence of a positive recurrent Markov chain to its steady state.
Therefore, for sufficiently large
queue lengths, the difference between the empirical rates in the
actual system and $r^*$ 
also decreases with $T$. This
ultimately results in a negative Lyapunov drift when $\boldsymbol \lambda$ is
inside the $\delta(T)$-stripped stability region
since from (\ref{objcon}) we have
$Q_1(jT)r_1^* + Q_2(jT)r_2^* > Q_1(jT)\lambda_1 + Q_2(jT)\lambda_2$.


The FBDC policy is easy to implement since it does not require the arrival rate information for stabilizing the
system for arrival rates in $\boldsymbol{\Lambda}-\delta(T )\boldsymbol{1}$, and it does not require
the solution of the LP (\ref{objcon}) for each frame. Instead, one can solve the LP (\ref{objcon}) for
all possible $(Q_1,Q_2)$ pairs only \emph{once} in advance and
create a mapping from $(Q_1,Q_2)$ pairs to the corners of the stability region.
Then, this mapping can be used to
find the corresponding optimal saturated-system policy 
to be applied during each frame.
Solving the LP in (\ref{objcon}) for all
possible $(Q_1, Q_2)$ pairs is possible because first,
the solution of the LP will be one of the corner points of the stability region in Fig.~\ref{Fig:throughput_eps_0_250_write_up},
and second, the weights $(Q_1,Q_2)$, which are the inputs to the LP, determine which corner point is optimal.
The theory of Linear Programming suggests that the solution to the LP in (\ref{objcon}) depends only on the relative value of the weights $(Q_1,Q_2)$ with respect to each other.
Namely, changing the queue size ratio $Q_2/Q_1$ varies the slope of the objective function of the LP in (\ref{objcon}), and the value of this slope $Q_2/Q_1$ with respect to
the slopes of the
lines in the stability region in
Fig.~\ref{Fig:throughput_eps_0_250_write_up} determine which corner point the FBDC policy operates on.
These mappings from the queue size ratios to the corners of the stability region are shown in Table~\ref{Fig:optimal_table_eps_small} for the case of $\epsilon <
\epsilon_c$ and in Table~\ref{Fig:optimal_table} for the case of $\epsilon \ge
\epsilon_c$. The corresponding mappings for the FBDC policy for the case of non-symmetric Gilbert-Elliot channels are shown in Appendix B.
%
Given these tables, one no longer needs to solve the LP (\ref{objcon}) for each frame, but just has to perform a simple table look-up to determine the optimal policy to use in each frame.
%

\begin{table}
\centering \psfrag{1\r}[l][][.8]{$0$}
\psfrag{2\r}[l][][.8]{$\!\!T^*_1(\epsilon)$}
\psfrag{3\r}[l][][.8]{\!\!$T^*_2(\epsilon)$}
\psfrag{4\r}[l][][.8]{\!$1$}
\psfrag{5\r}[l][][.8]{\!\!\!\!$T^*_3(\epsilon)$}
\psfrag{0\r}[l][][.8]{\!\!\!\!$T^*_4(\epsilon)$}
\psfrag{9\r}[l][][.8]{\!\!\!\!\!\!\!\!\!\!\!\!\!\!\!\!\!\!\textrm{corner}
$b_0$}
\psfrag{z\r}[l][][.8]{\!\!\!\!\!\!\!\!\!\!\!\!\!\!\!\!\!\textrm{corner}
$b_1$}
\psfrag{x\r}[l][][.8]{\!\!\!\!\!\!\!\!\!\!\!\!\textrm{corner}
$b_2$}
\psfrag{6\r}[l][][.8]{\!\!\!\!\!\!\!\!\!\!\!\!\!\!\textrm{corner}
$b_3$}
\psfrag{7\r}[l][][.8]{\!\!\!\!\!\!\!\!\!\!\!\!\!\!\!\!\!\textrm{corner}
$b_4$}
\psfrag{8\r}[l][][.8]{\!\!\!\!\!\!\!\!\!\!\!\!\!\!\!\!\!\textrm{corner}
$b_5$}
\psfrag{a\r}[l][][.75]{(1,1,1): \textrm{switch}}
\psfrag{b\r}[l][][.75]{(1,1,0): \textrm{switch}}
\psfrag{d\r}[l][][.75]{(1,0,1): \textrm{switch}}
\psfrag{e\r}[l][][.75]{(1,0,0): \textrm{switch}}
\psfrag{f\r}[l][][.75]{(2,1,1): \textrm{stay}}
\psfrag{g\r}[l][][.75]{(2,1,0): \textrm{stay}}
\psfrag{h\r}[l][][.75]{(2,0,1): \textrm{stay}}
\psfrag{i\r}[l][][.75]{(2,0,0): \textrm{stay}}
\psfrag{j\r}[l][][.75]{(1,1,1): \textrm{stay}}
\psfrag{k\r}[l][][.75]{(1,1,0): \textrm{stay}}
\psfrag{l\r}[l][][.75]{(2,1,0): \textrm{switch}}
\psfrag{m\r}[l][][.75]{(1,0,0): \textrm{stay}}
\psfrag{o\r}[l][][.75]{(2,0,0): \textrm{switch}}
\psfrag{p\r}[l][][.75]{(1,0,1): \textrm{stay}}
\psfrag{r\r}[l][][.75]{(2,1,1): \textrm{switch}}
\psfrag{s\r}[l][][.75]{(2,0,1): \textrm{switch}}
\psfrag{t\r}[l][][.75]{(2,0,0): \textrm{stay}}
\psfrag{u\r}[l][][.75]{(2,0,1): \textrm{stay}}
\psfrag{v\r}[l][][.75]{(2,0,0): \textrm{switch}}
\psfrag{w\r}[l][][.75]{(1,0,0): \textrm{stay}}
\psfrag{c\r}[l][][1]{$\!\!\!\!\!\!\frac{Q_2}{Q_1}$}
\includegraphics[width=0.50\textwidth]{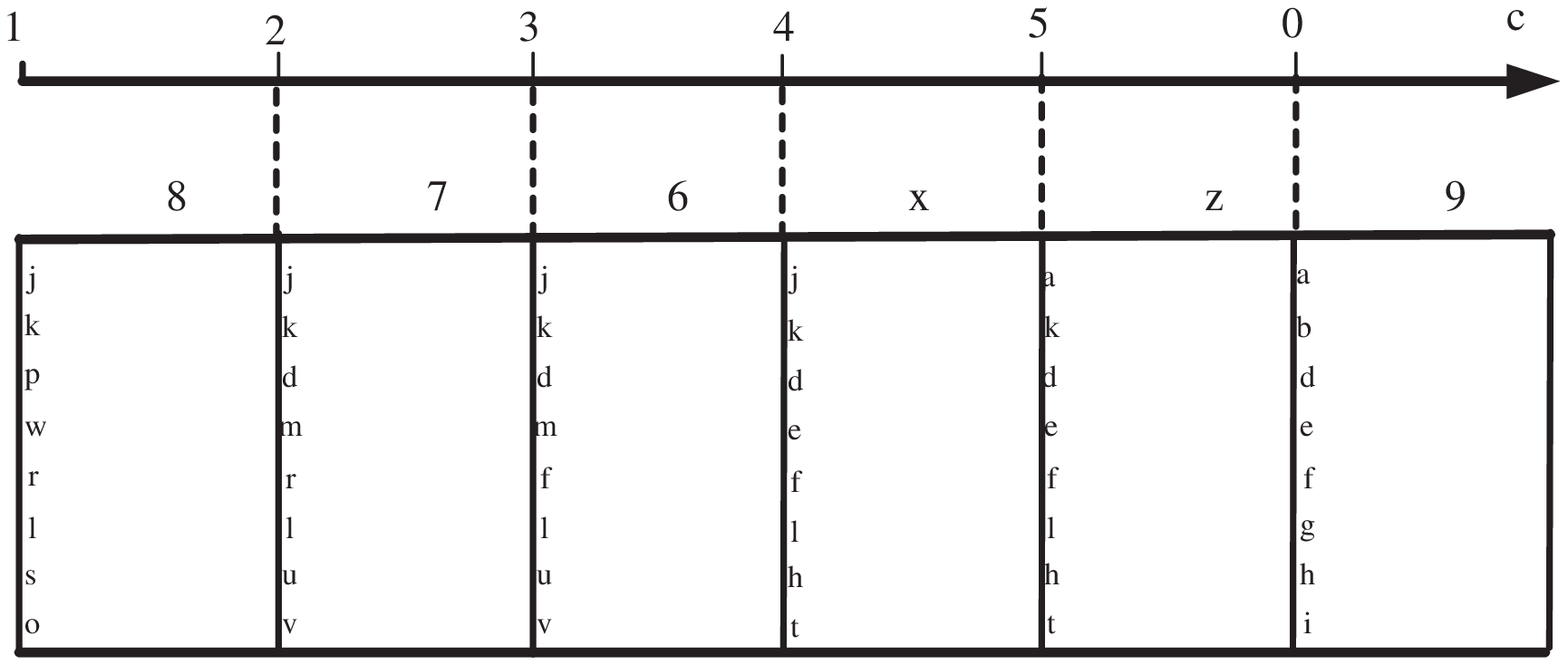}
\captionsetup{font=footnotesize} \caption{FBDC policy mapping from
the queue sizes to the corners of $\mathbf{\Lambda_s}$,
$b_0,b_1,b_2,b_3,b_4,b_5$ shown in Fig.
\ref{Fig:throughput_eps_0_250_write_up} (a), for $\epsilon <
\epsilon_c$.
For each state $\mathbf{s}=(m(t),C_1(t),C_2(t))$ the
optimal action is specified. The thresholds
on $Q_2/Q_1$ 
are $0,T^*_1=\epsilon/(1-\epsilon)^2,T^*_2=(1-\epsilon)/(1+\epsilon-\epsilon^2),1,
T^*_3=(1+\epsilon-\epsilon^2)/(1-\epsilon),T^*_4=(1-\epsilon)^2/\epsilon$.} \label{Fig:optimal_table_eps_small}
\end{table}s

\begin{table}
\centering \psfrag{1\r}[l][][.8]{$0$}
\psfrag{2\r}[l][][.8]{$\!\!T'_1(\epsilon)$}
\psfrag{3\r}[l][][.8]{\!\!$T^*_1(\epsilon)$}
\psfrag{4\r}[l][][.8]{\!$1$}
\psfrag{5\r}[l][][.8]{\!\!\!\!$T^*_2(\epsilon)$}
\psfrag{0\r}[l][][.8]{\!\!\!\!$T'_2(\epsilon)$}
\psfrag{9\r}[l][][.8]{\!\!\!\!\!\!\!\!\!\!\!\!\!\!\!\!\!\!\textrm{corner}
$b_0$}
\psfrag{6\r}[l][][.8]{\!\!\!\!\!\!\!\!\!\!\!\!\!\!\!\!\!\textrm{corner}
$b_1$} \psfrag{7\r}[l][][.8]{\!\!\!\!\!\!\!\!\!\!\!\!\textrm{corner}
$b_2$}
\psfrag{8\r}[l][][.8]{\!\!\!\!\!\!\!\!\!\!\!\!\!\!\textrm{corner}
$b_3$}
\psfrag{a\r}[l][][.8]{\!\!(1,1,1): \textrm{switch}}
\psfrag{b\r}[l][][.8]{\!\!(1,1,0): \textrm{switch}}
\psfrag{d\r}[l][][.8]{\!\!(1,0,1): \textrm{switch}}
\psfrag{e\r}[l][][.8]{\!\!(1,0,0): \textrm{switch}}
\psfrag{f\r}[l][][.8]{\!\!(2,1,1): \textrm{stay}}
\psfrag{g\r}[l][][.8]{\!\!(2,1,0): \textrm{stay}}
\psfrag{h\r}[l][][.8]{\!\!(2,0,1): \textrm{stay}}
\psfrag{i\r}[l][][.8]{\!\!(2,0,0): \textrm{stay}}
\psfrag{j\r}[l][][.8]{\!\!(1,1,1): \textrm{stay}}
\psfrag{k\r}[l][][.8]{\!\!(1,1,0): \textrm{stay}}
\psfrag{l\r}[l][][.8]{\!\!(2,1,0): \textrm{switch}}
\psfrag{m\r}[l][][.8]{\!\!(1,0,0): \textrm{stay}}
\psfrag{o\r}[l][][.8]{\!\!(2,0,0): \textrm{switch}}
\psfrag{p\r}[l][][.8]{\!\!(1,0,1): \textrm{stay}}
\psfrag{r\r}[l][][.8]{\!\!(2,1,1): \textrm{switch}}
\psfrag{s\r}[l][][.8]{\!\!(2,0,1): \textrm{switch}}
\psfrag{t\r}[l][][.8]{\!\!(2,0,0): \textrm{stay}}
\psfrag{u\r}[l][][.8]{\!\!(2,0,1): \textrm{stay}}
\psfrag{v\r}[l][][.8]{\!\!(2,0,0): \textrm{switch}}
\psfrag{w\r}[l][][.8]{\!\!(1,0,0): \textrm{stay}}
\psfrag{c\r}[l][][1]{$\!\!\!\!\!\!\frac{Q_2}{Q_1}$}
\includegraphics[width=0.50\textwidth]{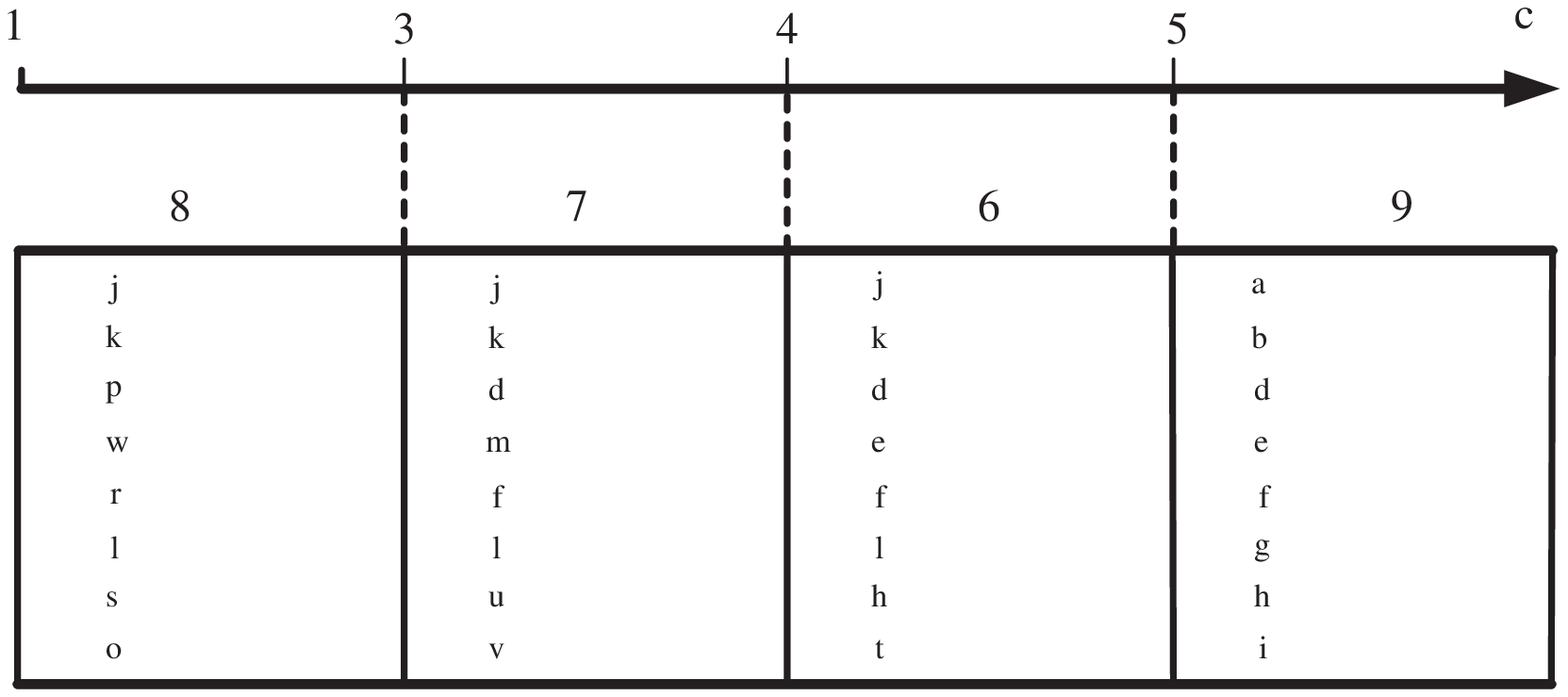}
\captionsetup{font=footnotesize} \caption{FBDC policy mapping from
the queue sizes to the corners of $\mathbf{\Lambda_s}$,
$b_0,b_1,b_2,b_3$ shown in Fig.
\ref{Fig:throughput_eps_0_250_write_up} (b), for $\epsilon \ge
\epsilon_c$. For each state $\mathbf{s}=(m(t),C_1(t),C_2(t))$ the
optimal action is specified. The thresholds
on $Q_2/Q_1$ 
are $0,T^*_1=1/((1-\epsilon)(3-2\epsilon)),1,T^*_2=(1-\epsilon)(3-2\epsilon)$.} \label{Fig:optimal_table}
\end{table}


%
%

In the next subsection we provide an upper bound to the long-run packet-average delay under the FBDC policy, which is linear in $T$. This suggest that the packet delay increases with increasing frame lengths as expected. However, such increases are at most linear in $T$.
Note that the FBDC policy can also be implemented without any frames by setting $T=1$, i.e., by solving the LP in Algorithm~\ref{alg:FBDC} in each time slot. The simulation results in Section~\ref{Sec:Sim} suggest that the FBDC policy implemented without frames
has a similar throughput performance to the original FBDC policy.
This is because for large queue lengths, the optimal solution of the LP in (\ref{objcon}) depends on the queue length ratios, and hence, the policy $\boldsymbol{\pi}^*$ that solves the LP optimally does not change fast when the queue lengths get large. When the policy is implemented without the use of frames, it becomes more adaptive to dynamic changes in the queue lengths, which results in a better delay performance than the frame-based implementations.

\subsubsection{Delay Upper Bound}\label{Sec:FBDC_delay}

The delay upper bound in this section is easily derived once the stability of the FBDC algorithm is established.
The stability proof utilizes the following quadratic Lyapunov function
\begin{equation*}
L(\mathbf{Q}(t)) =\sum_{i=1}^2 Q_i^2(t),
\end{equation*}
which represents a quadratic measure of the total load in the system at time slot $t$.
Let $t_k$ denote the time slots at the frame boundaries, $k=0,1,...,$ and define the $T$-step conditional
drift
\begin{equation*}
\Delta_{T}(t_k)\triangleq \mathbb{E} \left[ L(\mathbf{Q}(t_k+T)) -
L(\mathbf{Q}(t_k))\big|\mathbf{Q}(t_k)\right],
\end{equation*}
The following drift expression follows from the stability analysis in Appendix~C:
\begin{equation*}
 \frac{\Delta_T(t_k)}{2T} \leq BT - \Big(\sum_{i}Q_i(t_k)
\Big)\xi, 
\end{equation*}
where $B=1+A_{\max}^2$, $\boldsymbol{\lambda}$ is strictly inside the
$\mathbf{\delta}$-\emph{stripped} stability region
$\mathbf{\Lambda}- \delta\mathbf{1}$, and
$\xi >0$ represents a measure of the distance of $\boldsymbol{\lambda}$ to the boundary of $\mathbf{\Lambda}- \delta\mathbf{1}$.
Taking expectations with respect to $\mathbf{Q}(t_k)$,
writing a similar expression over the frame boundaries $t_k, k \in
\{0,1,2,...,K \}$, summing them and telescoping these expressions lead to
\begin{equation*}
L(\mathbf{Q}(t_{K})) -
L(\mathbf{Q}(0)) \le 2KBT^2  - 2\xi T \!\!\sum_{k=0}^{K-1} \!\!\mathbb{E}\left[\sum_{i}Q_i(t_k)\right]. 
\end{equation*}
Using $L(\mathbf{Q}(t_K))\ge 0$ and $L(\mathbf{Q}(0))=0$, we have 
\begin{equation*}
\displaystyle \limsup_{K\rightarrow \infty} \frac{1}{K} \sum_{k=0}^{K-1} \sum_i \mathbb{E}[Q_i(t_{k})] \le \frac{BT}{\xi}.
\end{equation*}
%
For $t\in(t_k,t_{k+1})$ we have
$Q_i(t) \le Q_i(t_k) + \sum_{\tau=0}^{T-1}A_i(t_k+\tau)$. Therefore,
$\mathbb{E}[Q_i(t)] \le \mathbb{E}[Q_i(t_k)] + T\lambda_i \le \mathbb{E}[Q_i(t_k)] + T A_{\max}$. Therefore, for $T_K \doteq KT$ we have
\begin{eqnarray*}
\displaystyle
&~& \!\!\!\!\!\!\!\!\!\!\!\!\!\limsup_{T_K\rightarrow \infty}\frac{1}{KT}   \sum_{t=0}^{KT-1}\sum_i \mathbb{E}[Q_i(t)] \\
&\!\!\!\!\!\!\!\!\le& \!\!\!\!\!\!\!\!\limsup_{K \rightarrow \infty} \frac{1}{TK}\!\! \sum_{k=0}^{K-1}\!\!\sum_i \!T\mathbb{E}[Q_i(t_{k})]\!+\!T^2\!A_{\max} \!\le\! \frac{(B \!\!+\!\!A_{\max}\xi) T}{\xi}.
\end{eqnarray*}
Dividing by the total arrival rate into the system $\sum_i\lambda_i$ and applying Little's law, the average delay is upper bounded by an expression that is linear in the frame length $T$.

In the next section we consider Myopic policies that do not require
the solution of an LP and that are able to stabilize the network for arrival rates within over 
$90\%$ of the stability region. 
Simulation results in Section \ref{Sec:Sim} suggest that the Myopic
policies may in fact achieve the full stability region while
providing better delay performance than
the FBDC policy for most arrival rates. 

\subsection{Myopic Control Policies}\label{Sec:Myopic_Pol}

We investigate the performance of simple \emph{Myopic}
policies that make scheduling/switching decisions according to weight
functions that are products of the queue lengths
and the channel
predictions for a small number of slots into the future. We refer to
a Myopic policy considering $k$ future time slots as the
$k$-Lookahead Myopic policy.
We implement these policies over frames of length $T$ time slots where during the $j$th frame,
the queue lengths at the beginning of the frame, $Q_1(jT)$ and $Q_2(jT)$, are used for weight calculations during the frame.
Specifically, in the $1$-Lookahead Myopic policy,
assuming that the server is with queue $1$ at some $t\in\{jT,...,j(T+1)-1\}$, the weight of queue $1$ is the product of
$Q_1(jT)$ and the summation of the current state of the channel process $C_1$ and the probability that $C_1$ will be in the ON state
at $t+1$. The weight of queue $2$ is calculated similarly, however, the current state of the channel process $C_2$ is not included in the weight since queue $2$ is not available to the server in the current time slot.
%
%
%
%
The detailed description of the $1$-Lookahead Myopic policy is given in Algorithm~\ref{alg:Myop1}
below. 
%
\begin{algorithm}[h]
\caption{\textsc{$1$-Lookahead Myopic Policy}} \label{alg:Myop1}
\begin{algorithmic}[1]

\STATE Assuming that the server is currently with queue $1$ and the
system is at the $j$th frame, calculate the following weights in
each time slot of the current frame; \beqn
W_1(t) \!\!\!\!&=\!\!\!\!&
Q_1(jT)\Big(C_1(t)+\mathbb{E}\big[C_1(t+1)|C_1(t)\big]\Big)\nonumber\\
W_2(t) \!\!\!\!&=\!\!\!\!&
Q_2(jT)\mathbb{E}\big[C_2(t+1)|C_2(t)\big].\label{eq:myopic}
\eeqn \STATE
If $W_1(t)\ge W_2(t)$ stay with queue 1, 
otherwise, switch to the other queue. A similar rule applies for queue
2. 
%
\end{algorithmic}
\end{algorithm}

%
%


Next, we establish a lower bound on the stability region of the
1-Lookahead Myopic Policy 
by comparing its drift over a frame to the drift of the FBDC policy.
\begin{theorem}\label{thm:Myopic_1LH}
\emph{The 1-Lookahead Myopic policy achieves at least
$\gamma$-fraction of the stability region $\mathbf{\Lambda_s}$
asymptotically in $T$ where $\gamma \ge 90\%$.}
\end{theorem}
%
The proof is constructive and will be establish in various steps in the following.
The basic idea behind the proof is
that the 1-Lookahead Myopic (OLM) policy produces a mapping from the
set of queue sizes to the stationary deterministic policies
corresponding to the corners of the stability region. This mapping
is similar to that of the FBDC policy, however, the thresholds on
the queue size ratios $Q_2/Q_1$ are determined according to
(\ref{eq:myopic}): \\\\
%
%
%
\textbf{Mapping from queue sizes to actions. Case-1: $\epsilon < \epsilon_c$ }\\
For $\epsilon < \epsilon_c$, there are 6 corners in the stability
region 
denoted by $b_0,
b_1,...,b_5$ where $b_0$ is $(0,0.5)$ and $b_5$ is $(0.5,0)$ 
as shown in Fig.~\ref{Fig:throughput_eps_0_250_write_up} (a). We
derive conditions on $Q_2/Q_1$ such that the OLM policy chooses the
stationary deterministic decisions that correspond to a given corner
point.\\
\textbf{Corner $b_0$}:\\
Optimal actions are to stay at queue-2 for every channel condition.
Therefore, the server chooses queue-2 even when the channel state is
$C_1(t),C_2(t)=(1,0)$. Therefore, using (\ref{eq:myopic}), for the
Myopic policy to take the deterministic actions corresponding to
$b_0$ we need
\begin{equation*}
Q_1.(1-\epsilon) < Q_2.(\epsilon)\;\;\Rightarrow \frac{Q_2}{Q_1} >
\frac{1-\epsilon}{\epsilon}.
\end{equation*}
This means that if we apply the Myopic policy with coefficients
$Q_1, Q_2$ such that $Q_2/Q_1 > (1-\epsilon)/\epsilon$, then the
system output rate will be driven towards the corner point $b_0$
(both in the saturated system or in the actual system with large enough arrival rates). \\
\textbf{Corner $b_1$}:\\
The optimal actions for the corner point $b_1$ are as follows: At
queue-1, for the channel state $10$:stay, for the channel states
$11$, $01$ and $00$: switch. At queue-2, for the channel state $10$:
switch, for the channel states $11$, $01$ and $00$: stay. The most
limiting conditions are $11$ at queue-1 and $10$ at queue-2.
Therefore we need, $Q_1(2-\epsilon) < Q_2 (1-\epsilon)$ and
$Q_1(1-\epsilon) > Q_2\epsilon$. Combining these we have
\begin{equation*}
 \frac{2-\epsilon}{1-\epsilon}<\frac{Q_2}{Q_1} <
\frac{1-\epsilon}{\epsilon}.
\end{equation*}
Note that the condition $\epsilon < \epsilon_c = 1-\sqrt{2}/2$
implies that $\frac{1-\epsilon}{\epsilon}
> \frac{2-\epsilon}{1-\epsilon}$.\\
\textbf{Corner $b_2$}:\\
The optimal actions for the corner point $b_1$ are as follows: At
queue-1, for the channel state $10$ and $11$:stay, for the channel
states $01$ and $00$: switch. At queue-2, for the channel states
$10$: switch, for the channel states $11$, $01$ and $00$: stay. The
most limiting conditions are $11$ at queue-1 and $00$. Therefore we
need, $Q_1(2-\epsilon) > Q_2 (1-\epsilon)$ and $Q_1 < Q_2$.
Combining these we have
\begin{equation*}
1 <\frac{Q_2}{Q_1} <  \frac{2-\epsilon}{1-\epsilon}.
\end{equation*}
The conditions for the rest of the corners are symmetric and can be
found similarly to obtain the mapping in Fig.~\ref{Fig:myopic_table_eps_small}.\\
\\
\textbf{Mapping from queue sizes to actions. Case-2: $\epsilon \ge
\epsilon_c$ }\\
In this case there are 4 corner points in the
throughput region. We enumerate these corners as $b_0, b_2,b_3,b_5$
where $b_0$ is $(0,0.5)$ and $b_5$ is $(0.5,0)$.\\
\textbf{Corner $b_0$}:\\
The analysis is the same as the $b_0$ analysis in the previous case
and we obtain that for the Myopic policy to take the deterministic
actions corresponding to $b_0$ we need
\begin{equation*}
 \frac{Q_2}{Q_1} >
\frac{1-\epsilon}{\epsilon}.
\end{equation*}\\
\textbf{Corner $b_2$}:\\
This is the same corner point as in the previous case corresponding
to the same deterministic policy: At queue-1, for the channel state
$10$ and $11$:stay, for the channel states $01$ and $00$: switch. At
queue-2, for the channel states $10$: switch, for the channel states
$11$, $01$ and $00$: stay. The most limiting conditions are $10$ at
queue-2 (since $\epsilon \ge \epsilon_c$ we have
$\frac{1-\epsilon}{\epsilon} < \frac{2-\epsilon}{1-\epsilon}$) and
$00$. Therefore we need, $Q_1(1-\epsilon)
> Q_2\epsilon$ and $Q_1 < Q_2$. Combining these we have
\begin{equation*}
1 <\frac{Q_2}{Q_1} <  \frac{1-\epsilon}{\epsilon}.
\end{equation*}
The conditions for the rest of the corners are symmetric and can be
found similarly to obtain the mapping in Fig. \ref{Fig:myopic_optimal_table}
for $\epsilon \ge\epsilon_c$.\\

The conditions for the corners $b_2$ and $b_3$ are symmetric, completing the mapping from the queue
sizes to the corners of $\mathbf{\Lambda_s}$ for $\epsilon \ge \epsilon_c$ 
shown in
Table~\ref{Fig:myopic_optimal_table}. This mapping is in general different from the corresponding mapping of the FBDC policy in Table~\ref{Fig:optimal_table}. Therefore, for a given ratio of the
queue sizes $Q_2/Q_1$,
the FBDC and
the OLM policies \emph{may} apply different stationary deterministic policies corresponding to different corner
points of $\mathbf{\Lambda_s}$, denoted by $\mathbf{r}^*$ and
$\hat{\mathbf{r}}$ respectively.
The shaded intervals of $Q_2/Q_1$ in Table~\ref{Fig:myopic_optimal_table} are
the intervals in which the OLM and the FBDC policies apply different policies.
A similar mapping can be obtained for the
OLM policy
for $\epsilon < \epsilon_c$.
The corresponding mapping for the OLM policy for the case of non-symmetric Gilbert-Elliot channels is given in Appendix~B.

The following lemma is proved in
Appendix E and completes the proof by establishing the $90\%$
bound on the weighted average departure rate of the OLM policy
w.r.t. to that of the FBDC policy.
\begin{lemma}\label{lem:My_FBDC_comp}
\emph{We have that
\begin{equation}
\Psi \doteq \frac{\sum_{i} Q_i(t) \hat{r}_i}{\sum_{i} Q_i(t)
r_i^{*}} \ge 90\%. \label{eq:Psi_bound}
\end{equation}
Furthermore, $\Psi \ge 90\%$ is a sufficient condition for the OLM
policy to achieve at least $90\%$ of $\mathbf{\Lambda_s}$
asymptotically in $T$.}
\end{lemma}

\begin{table}
\centering \psfrag{1\r}[l][][.8]{$0$}
\psfrag{2\r}[l][][.8]{\,$T^*_1(\epsilon)$}
\psfrag{n\r}[l][][.8]{$\!\!\!\!\!\!\!T_1(\epsilon)$}
\psfrag{3\r}[l][][.8]{$\,T^*_2(\epsilon)$}
\psfrag{q\r}[l][][.8]{\!\!\!\!\!$T_2(\epsilon)$}
\psfrag{4\r}[l][][.8]{\!$1$}
\psfrag{5\r}[l][][.8]{\!\!\!\!\!\!\!$T^*_3(\epsilon)$}
\psfrag{y\r}[l][][.8]{\!\!$T_3(\epsilon)$}
\psfrag{0\r}[l][][.8]{\!\!\!\!\!\!\!$T^*_4(\epsilon)$}
\psfrag{w\r}[l][][.8]{\!\!\!\!$T_4(\epsilon)$}
\psfrag{9\r}[l][][.8]{\!\!\!\!\!\!\!\!\!\!\!\!\!\!\!\!\!\!\textrm{corner}
$b_0$}
\psfrag{z\r}[l][][.8]{\!\!\!\!\!\!\!\!\!\!\!\!\!\!\!\!\!\textrm{corner}
$b_1$}
\psfrag{x\r}[l][][.8]{\!\!\!\!\!\!\!\!\!\!\!\!\textrm{corner}
$b_2$}
\psfrag{6\r}[l][][.8]{\!\!\!\!\!\!\!\!\!\!\!\!\!\!\textrm{corner}
$b_3$}
\psfrag{7\r}[l][][.8]{\!\!\!\!\!\!\!\!\!\!\!\!\!\!\!\!\!\textrm{corner}
$b_4$}
\psfrag{8\r}[l][][.8]{\!\!\!\!\!\!\!\!\!\!\!\!\!\!\!\!\!\textrm{corner}
$b_5$}
\psfrag{a\r}[l][][.75]{(1,1,1): \textrm{switch}}
\psfrag{b\r}[l][][.75]{(1,1,0): \textrm{switch}}
\psfrag{d\r}[l][][.75]{(1,0,1): \textrm{switch}}
\psfrag{e\r}[l][][.75]{(1,0,0): \textrm{switch}}
\psfrag{f\r}[l][][.75]{(2,1,1): \textrm{stay}}
\psfrag{g\r}[l][][.75]{(2,1,0): \textrm{stay}}
\psfrag{h\r}[l][][.75]{(2,0,1): \textrm{stay}}
\psfrag{i\r}[l][][.75]{(2,0,0): \textrm{stay}}
\psfrag{j\r}[l][][.75]{(1,1,1): \textrm{stay}}
\psfrag{k\r}[l][][.75]{(1,1,0): \textrm{stay}}
\psfrag{l\r}[l][][.75]{(2,1,0): \textrm{switch}}
\psfrag{m\r}[l][][.75]{(1,0,0): \textrm{stay}}
\psfrag{o\r}[l][][.75]{(2,0,0): \textrm{switch}}
\psfrag{p\r}[l][][.75]{(1,0,1): \textrm{stay}}
\psfrag{r\r}[l][][.75]{(2,1,1): \textrm{switch}}
\psfrag{s\r}[l][][.75]{(2,0,1): \textrm{switch}}
\psfrag{t\r}[l][][.75]{(2,0,0): \textrm{stay}}
\psfrag{u\r}[l][][.75]{(2,0,1): \textrm{stay}}
\psfrag{v\r}[l][][.75]{(2,0,0): \textrm{switch}}
\psfrag{c\r}[l][][1]{$\!\!\!\!\!\!\frac{Q_2}{Q_1}$}
\includegraphics[width=0.50\textwidth]{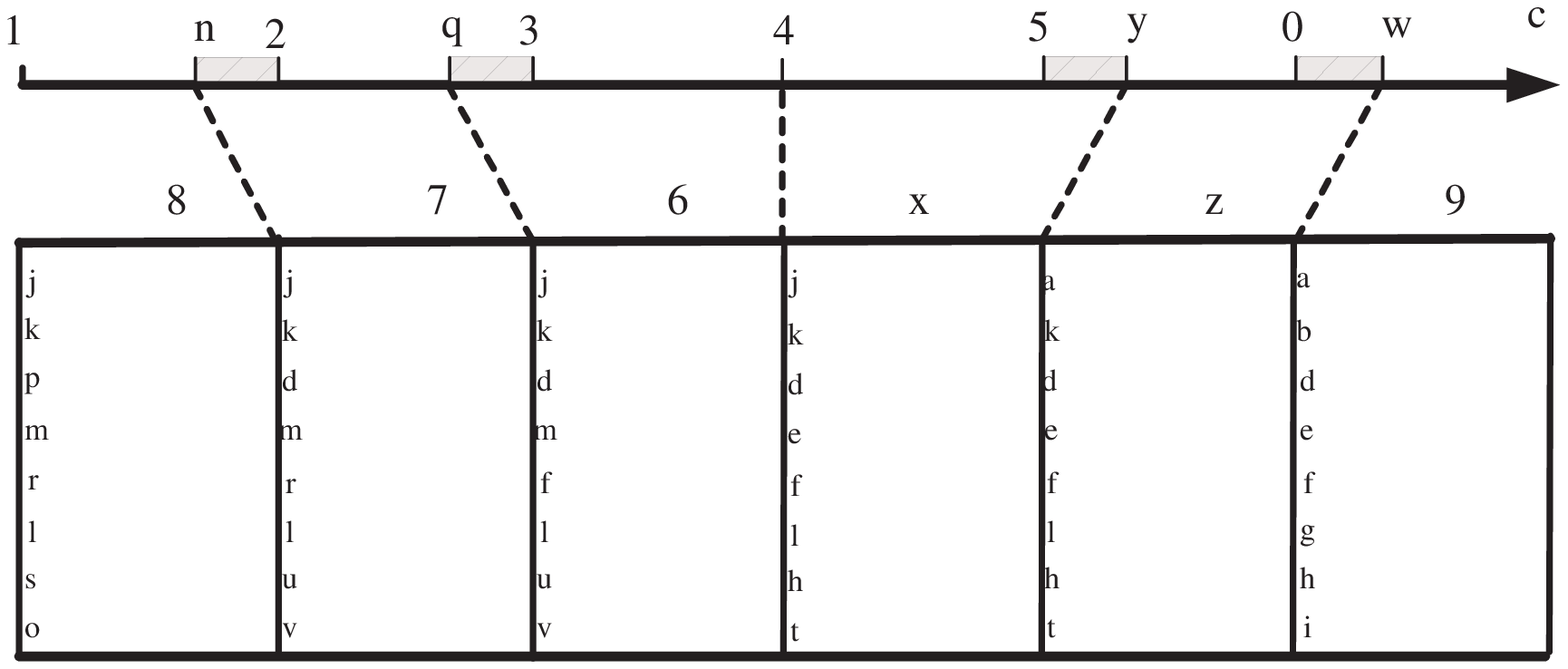}
\captionsetup{font=footnotesize} \caption{1 Lookahead Myopic policy mapping from
the queue sizes to the corners of $\mathbf{\Lambda_s}$,
$b_0,b_1,b_2,b_3,b_4,b_5$ shown in Fig.
\ref{Fig:throughput_eps_0_250_write_up} (a), for $\epsilon <
\epsilon_c$. For each state $\mathbf{s}=(m(t),C_1(t),C_2(t))$ the
optimal action is specified. The thresholds
on $Q_2/Q_1$ 
are 
$0,T_1=\epsilon/(1-\epsilon),T_2=(1-\epsilon)/(2-\epsilon),1,
T_3=(2-\epsilon)/(1-\epsilon),T_4=(1-\epsilon)/\epsilon$.
The corresponding thresholds for the FBDC policy are $0,T^*_1,T^*_2,1,T^*_3,T^*_4$.
For example, corner
$b_2$ is chosen in the FBDC policy if $1\le Q_2/Q_1 < T^*_3$, whereas in
the OLM policy if $1 \le Q_2/Q_1 < T_3$.}
\label{Fig:myopic_table_eps_small}
\end{table}

\begin{table}
\centering \psfrag{1\r}[l][][.8]{$0$}
\psfrag{2\r}[l][][.8]{$\!\!T_1(\epsilon)$}
\psfrag{3\r}[l][][.8]{\!\!$T^*_1(\epsilon)$}
\psfrag{4\r}[l][][.8]{\!$1$}
\psfrag{5\r}[l][][.8]{\!\!\!\!$T^*_2(\epsilon)$}
\psfrag{0\r}[l][][.8]{\!\!\!\!$T_2(\epsilon)$}
\psfrag{9\r}[l][][.8]{\!\!\!\!\!\!\!\!\!\!\!\!\!\!\!\!\!\!\textrm{corner}
$b_0$}
\psfrag{6\r}[l][][.8]{\!\!\!\!\!\!\!\!\!\!\!\!\!\!\!\!\!\textrm{corner}
$b_2$} \psfrag{7\r}[l][][.8]{\!\!\!\!\!\!\!\!\!\!\!\!\textrm{corner}
$b_3$}
\psfrag{8\r}[l][][.8]{\!\!\!\!\!\!\!\!\!\!\!\!\!\!\textrm{corner}
$b_5$}
\psfrag{a\r}[l][][.8]{\!\!(1,1,1): \textrm{switch}}
\psfrag{b\r}[l][][.8]{\!\!(1,1,0): \textrm{switch}}
\psfrag{d\r}[l][][.8]{\!\!(1,0,1): \textrm{switch}}
\psfrag{e\r}[l][][.8]{\!\!(1,0,0): \textrm{switch}}
\psfrag{f\r}[l][][.8]{\!\!(2,1,1): \textrm{stay}}
\psfrag{g\r}[l][][.8]{\!\!(2,1,0): \textrm{stay}}
\psfrag{h\r}[l][][.8]{\!\!(2,0,1): \textrm{stay}}
\psfrag{i\r}[l][][.8]{\!\!(2,0,0): \textrm{stay}}
\psfrag{j\r}[l][][.8]{\!\!(1,1,1): \textrm{stay}}
\psfrag{k\r}[l][][.8]{\!\!(1,1,0): \textrm{stay}}
\psfrag{l\r}[l][][.8]{\!\!(2,1,0): \textrm{switch}}
\psfrag{m\r}[l][][.8]{\!\!(1,0,0): \textrm{stay}}
\psfrag{o\r}[l][][.8]{\!\!(2,0,0): \textrm{switch}}
\psfrag{p\r}[l][][.8]{\!\!(1,0,1): \textrm{stay}}
\psfrag{r\r}[l][][.8]{\!\!(2,1,1): \textrm{switch}}
\psfrag{s\r}[l][][.8]{\!\!(2,0,1): \textrm{switch}}
\psfrag{t\r}[l][][.8]{\!\!(2,0,0): \textrm{stay}}
\psfrag{u\r}[l][][.8]{\!\!(2,0,1): \textrm{stay}}
\psfrag{v\r}[l][][.8]{\!\!(2,0,0): \textrm{switch}}
\psfrag{w\r}[l][][.8]{\!\!(1,0,0): \textrm{stay}}
\psfrag{c\r}[l][][1]{$\!\!\!\!\!\!\frac{Q_2}{Q_1}$}
\includegraphics[width=0.50\textwidth]{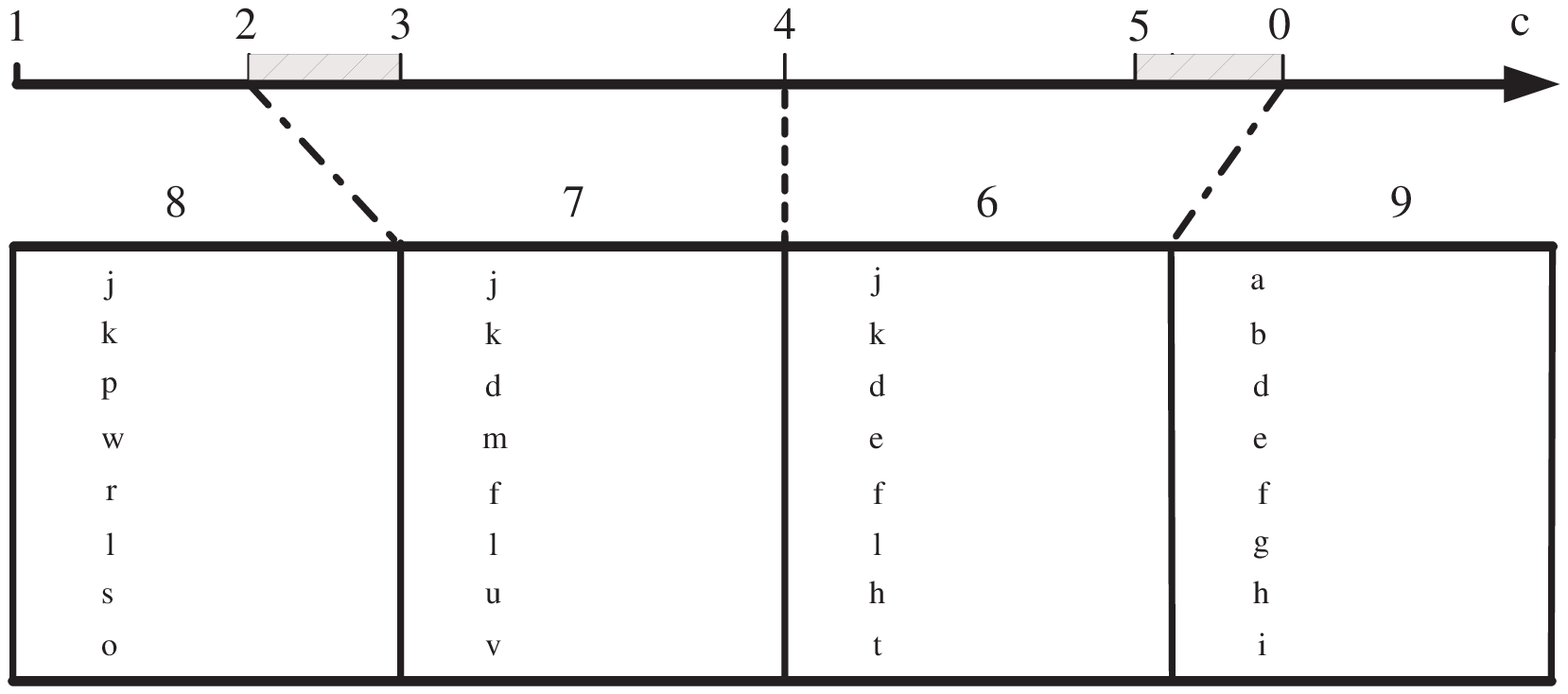}
\captionsetup{font=footnotesize} \caption{1 Lookahead Myopic policy mapping from
the queue sizes to the corners of $\mathbf{\Lambda_s}$,
$b_0,b_1,b_2,b_3$ shown in Fig.
\ref{Fig:throughput_eps_0_250_write_up} (b), for $\epsilon \ge
\epsilon_c$. For each state $\mathbf{s}=(m(t),C_1(t),C_2(t))$ the
optimal action is specified. The thresholds
on $Q_2/Q_1$ 
are $0,T_1=\epsilon/(1-\epsilon),1,T_2=(1-\epsilon)/\epsilon$.
The corresponding thresholds for the FBDC policy are $0,T^*_1,1,T^*_2$.
For example, corner
$b_1$ is chosen in the FBDC policy if $1\le Q_2/Q_1 < T^*_2$, whereas in
the OLM policy if $1 \le Q_2/Q_1 < T_2$.}
\label{Fig:myopic_optimal_table}
\end{table}

A similar analysis shows that the \emph{2-Lookahead Myopic Policy}
achieves at least $94\%$ of $\mathbf{\Lambda}_s$, while the
\emph{3-Lookahead Myopic Policy} achieves at least $96\%$ of
$\mathbf{\Lambda}_s$. The $k$-Lookahead Myopic Policy is the same as
before except that the following weight functions are used for
scheduling decisions: Assuming the server is with queue 1 at time
slot $t$,\\
$W_1(t)
=Q_1(jT)\big(C_1(t)+\sum_{\tau=1}^k\mathbb{E}[C_1(t+\tau)|C_1(t)]\big)$
and $W_2(t) =Q_2(jT)\sum_{\tau=1}^k\mathbb{E}[C_2(t+\tau)|C_2(t)]$.

\section{General System}\label{Sec:General_sys}
In this section we extend the results developed in the previous section to the general case of an arbitrary number of queues in the system.

\subsection{Model}\label{Sec:model_gen}

Consider the same model as in Section \ref{Sec:Model} with $N>1$ queues for some $N\in \mathbb{N}$ as shown in Fig. \ref{Fig:Parallel_queues_servers_switchover_Rconnect}. Let the i.i.d. process $A_i(t)$ with arrival rate $\lambda_i$ denote the number of arrivals to queue $i$ at time slot $t$, where
$\mathbb{E}[A_i^2(t)]\le
A_{\max}^2$, $i \in \{1,2,...,N\}$. 
Let $C_i(t)$
be the channel (connectivity) process of queue $i$, $i \in \{1,2,...,N\}$,
that forms the two-state
Markov chain with transition probabilities $p_{01}$ and $p_{10}$  as
shown in Fig.~\ref{Fig:channel_MC}.
%
%
We assume that the
processes $A_i(t),i \in \{1,...,N\}$ and $C_i(t),i\in \{1,...,N\}$ are independent. 
%
It takes one slot for the server to switch from one queue to the
other, and $m(t) \in \{1,...,N\}$ denotes
the queue at which the server is present at 
slot $t$. 
Let $\mathbf{s}_t=(m(t),C_1(t),...,C_N(t)) \in \mathcal{S}$ denote the state of the corresponding saturated system at time $t$ where $\mathcal{S}$ is the set of all states. The action $a(t)$ in each time slot is to choose the queue at which the server will be present in the next time slot, i.e., $a_t \in \{1,...,N\} \doteq \mathcal{A}$ where $\mathcal{A}$ is the
set of all actions at each state.

\subsection{Stability Region}\label{Sec:stab_reg_gen}


In this section we characterize the stability region of the general system under non-symmetric channel models\footnote{For Markovian (Gilbert-Elliot) channels, we preserve the symmetry of the channel processes across the queues.}.
For the case of i.i.d. channel processes we explicitly characterize the stability region and the throughput-optimal policy. For Markovian channel models, we extend the stability region characterization in terms of state-action frequencies to the general system.
Furthermore, we develop a tight outer bound on the stability region using an upper bound on the sum-throughput and show that a simple myopic policy
achieves this upper bound for the corresponding saturated system.

A dynamic server allocation problem over parallel channels with randomly varying connectivity and limited channel sensing has been investigated in \cite{AhmadTaraKrish09,AhmadLiu,ZhaoKrishLiu08} under the Gilbert-Elliot channel model. The goal in \cite{AhmadTaraKrish09,AhmadLiu,ZhaoKrishLiu08} is to maximize the sum-rate for the saturated system, where it is proved that a myopic policy is optimal. In this section we prove that a myopic policy is sum-rate optimal under the Gilbert-Elliot channel model and 1-slot server switchover delay. Furthermore, our goal is to characterize the set of all achievable rates, i.e., the stability region, together with a throughput-optimal scheduling algorithm for the dynamic queuing system.

\vspace{2mm}
\subsubsection{Memoryless Channels}
The results established in Section \ref{Sec:Motiv} for the case of i.i.d. connectivity processes can easily be extended to the system of $N$ queues with non-symmetric i.i.d. channels as the same intuition applies for the general case. We state this result in the following theorem whose proof can be found in Appendix~A.
\begin{theorem}
\emph{For a system of $N$ queues with 
\emph{arbitrary} switchover times 
and i.i.d. channels with probabilities
$p_i, i \in \{1,...,N \}$, the stability region $\mathbf{\Lambda}$ is given by
\begin{equation*}
\mathbf{\Lambda}=\left\{ \boldsymbol{\lambda} \ge \mathbf{0} \bigg | \sum_{i=1}^N \frac{\lambda_i}{p_i} \le 1\right\}.
\end{equation*}
In addition, the simple Exhaustive (Gated) policy is
throughput-optimal.}
\end{theorem}
As for the case of two queues, the simultaneous presence of randomly varying connectivity and the switchover delay significantly reduces the stability region as compared to the corresponding system without switchover delay analyzed in \cite{tass93}. Furthermore, when the channel processes are memoryless, no policy can take advantage of the channel diversity
as the simple queue-blind Exhaustive-type policies are throughput-optimal.

In the next section, we show that, similar to the case of two queues, the memory in the channel improves the stability region of the general system.

\vspace{2mm}
\subsubsection{Channels With Memory}

Similar to Section \ref{Sec:Stab_Reg}, we start by establishing the rate region $\mathbf{\Lambda}_s$ by formulating an MDP for rate maximization in the corresponding saturated system. 
The reward
functions in this case are given as follows:
\beqn
\!\!\!\!\!\!\overline{r}_i(\mathbf{s},a) \!\!\!\!\!&\doteq& \!\!\!\!\!1 \textrm{ if } m=i, C_i=1  \textrm{, and } a\!=\!i,\; i=1,...,N, \label{eq:r_i-N}
\eeqn
and $\overline{r}_i(\mathbf{s},a)\doteq0$ otherwise, where $m$ denotes the queue at which the server is present. That is, one reward is
obtained when the server stays at a queue with an ON channel. Given some
$\alpha_i \ge 0, \, i\in\{1,...,N\}$,  
$\sum_i\alpha_i=1$,
we define the system reward at time $t$ as
\beqn
\overline{r}(\mathbf{s},a) \doteq \sum_{i=1}^N \alpha_i \overline{r}_i(\mathbf{s},a).\nonumber
\eeqn
The average reward of policy $\boldsymbol{\pi}$ is defined as
\begin{equation*}
r^{\boldsymbol{\pi}} \doteq \displaystyle \lim_{K\rightarrow \infty} \frac{1}{K} E
\Big \{ \sum_{t=1}^K \overline{r}(\mathbf{s}_t,a_t^{\boldsymbol{\pi}}) \Big \}.
\end{equation*}
Therefore, the problem of maximizing
the time average expected reward over all policies,
$r^{*} \doteq \max_{\boldsymbol{\pi}} r^{\boldsymbol{\pi}}$, 
is a discrete time 
MDP characterized by the state transition probabilities
$\mathbf{P}(\mathbf{s}'|\mathbf{s},a)$ with $N2^N$ states and $N$ possible actions per state.
Furthermore, similar to the two-queue system, there exists a positive probability path between any given pair of states under some stationary-deterministic policy. Therefore, this MDP belongs to the class of \emph{Weakly Communicating}
MDPs \cite{puterman05} for which there exists a stationary-deterministic optimal policy independent of the initial state \cite{puterman05}.
The \emph{state-action polytope}, $\mathbf{X}$ is the set of
state-action frequency vectors $\mathbf{x}$ that satisfy the balance
equations
\beqn
\sum_{a\in \mathcal{A}}\mathbf{x}(\mathbf{s},a) = \sum_{\mathbf{s'} \in \mathcal{S}}
\sum_{a\in \mathcal{A}} \mathbf{P}\big(\mathbf{s}|\mathbf{s}',a\big)\mathbf{x}(\mathbf{s}',a) ,
\;\forall\; \mathbf{s} \in \mathcal{S}, \label{eq:saf_balance_N}
\eeqn
the normalization
condition
\beqn
\sum_{\mathbf{s} \in \mathcal{S}}\sum_{a\in \mathcal{A}}\mathbf{x}(\mathbf{s},a)=1,  \nonumber
\eeqn
and the nonnegativity constraints
\beqn
\mathbf{x}(\mathbf{s},a)\ge 0,  \:
\mbox{for} \: \mathbf{s} \in \mathcal{S}, a\in \mathcal{A}, \nonumber
\eeqn
where the transition probabilities $\mathbf{P}\big(\mathbf{s}|\mathbf{s}',a\big)$ are functions of the channel parameters $p_{10}$ and $p_{01}$.
The following linear transformation of the state-action polytope
$\mathbf{X}$ defines the \emph{rate polytope} $\mathbf{\Lambda}_s$, namely, the set of all time average
expected rate pairs that can be obtained in the saturated
system.
\begin{equation*}
\mathbf{\Lambda}_s \!=\!\!
\Big\{\mathbf{r}\big| r_i \!=\!\!
\sum_{\mathbf{s}\in \mathcal{S}}\sum_{a\in\mathcal{A}}\!\mathbf{x}(\mathbf{s},a)\overline{r}_i(\mathbf{s},a),
\mathbf{x} \in \mathbf{X}, i \in \{1,2,...,N\} \Big\},
\end{equation*}
where the reward functions $\overline{r}_i(\mathbf{s},a), i \in \{1,...,N \}$, are defined in (\ref{eq:r_i-N}). Algorithm \ref{alg:LP_for_Lambda_s_N}
gives an alternative characterization of the rate region $\mathbf{\Lambda}_s$.

\begin{figure}
\begin{center}
\includegraphics[width=0.45\textwidth]{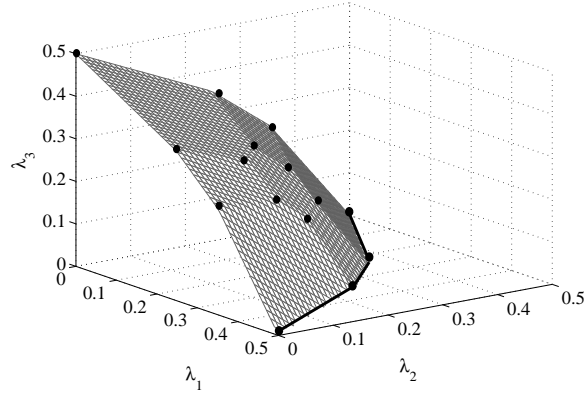}
\end{center}
\captionsetup{font=footnotesize}\caption{Stability region for
3 parallel queues for $p_{10}=p_{01}=0.3$.} \label{Fig:stab_n3_m1eps03}
\end{figure}

\begin{algorithm}[h]
\caption{\emph{Stability Region Characterization}}
\label{alg:LP_for_Lambda_s_N}
\begin{algorithmic}[1]

\STATE Given $\alpha_1,...,\alpha_N\!\! \ge\! 0$, $\sum_i\alpha_i=1$, solve the following LP
\begin{eqnarray}
\max_{\mathbf{x}}.  \quad  \sum_{i=1}^N \alpha_i r_i(\mathbf{x})\nonumber  \\
\mbox{subject to} \quad  \mathbf{x} \in \mathbf{X}.\label{eq:LP_N}
\end{eqnarray}

\STATE There exists an optimal solution
$(r_1^{*},...,r_N^{*})$ of this LP that lies at a corner point of 
$\mathbf{\Lambda_s}$. Find all possible corner points and take their
convex combination.
\end{algorithmic}
\end{algorithm}
Similar to the two-queue case, the fundamental theorem of Linear Programming guarantees existence of an optimal solution to (\ref{eq:LP_N}) at a corner point of the polytope $\mathbf{X}$ \cite{BertTsit97}.
We will establish in the next section that the rate region, $\boldsymbol{\Lambda}_s$ is in fact achievable in the dynamic queueing system, which will imply that $\mathbf{\Lambda}=\mathbf{\Lambda}_s$. For the case of 3 queues, Fig.~\ref{Fig:stab_n3_m1eps03} shows the stability region $\Lambda$. As expected, the stability region is significantly reduced as compared to the corresponding system with zero switchover delays analyzed in \cite{tass93}.

\vspace{3mm}
\subsubsection{Analytical Outer Bound For The Stability Region}

In this section we first derive an upper bound to the sum-throughput in the saturated system and then use it to characterize an outer bound to the rate region $\mathbf{\Lambda}_s$. Let $C_0^{(N)} \doteq \frac{p_{10}^N}{(p_{10}+p_{01})^N}$ denote the probability that all channels are in OFF state in steady state.

\begin{lemma}\label{lemma:sum_tp_Nqueues}
\emph{An upper bound on the sum-rate in the saturated system is given by}
\begin{equation}\label{eq:sum-rate_ub}
\sum_{i=1}^N r_i \le 1-C_0^{(N)} - \Big(p_{10}(1-C_0^{(N)}) - p_{01}C_0^{(N)}\Big).
\end{equation}
\end{lemma}
The proof is given in Appendix F. In the next section we propose a simple myopic policy for the saturated system that achieves this upper bound. Similar to the case of two-queues, the surface $\sum_{i=1}^N r_i \le 1-C_0^{(N)}$ is one of the boundaries of the stability region for the system without switchover delay analyzed in \cite{tass93}, where the probability that at least 1 channel is in ON state in steady state is $1-C_0^{(N)}$. Therefore, $p_{10}(1-C_0^{(N)}) - p_{01}C_0^{(N)}$ is \emph{the throughput loss due to 1 slot switchover delay} in our system. The analysis of the myopic policy in the next section shows that this throughput loss due to switchover delay
corresponds to the probability that the server is at a queue with OFF state when at least one other queue is in ON state. For the case of $N=3$ queues, the sum-throughput upper bound in Lemma \ref{lemma:sum_tp_Nqueues} is the hexagonal region at the center of the plot in Fig.~\ref{Fig:stab_n3_m1eps03}.

\begin{figure}
\begin{center}
\includegraphics[width=0.45\textwidth]{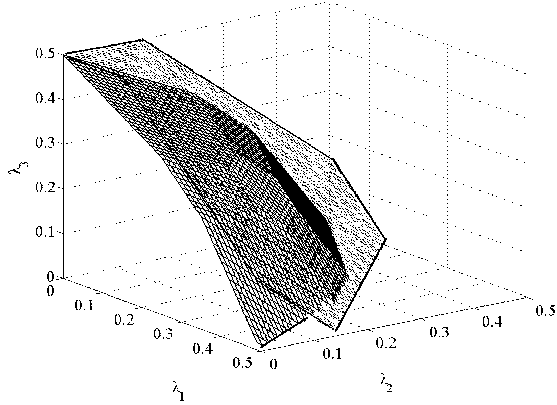}
\end{center}
\captionsetup{font=footnotesize}\caption{Stability region outer bound for
3 parallel queues for $p_{10}=p_{01}=0.3$.} \label{Fig:stab_n3_m1eps03_ub}
\end{figure}

Because any convex combination of $r_i, i \in \{1,...,N\}$, must lie under the sum-rate surface, (\ref{eq:sum-rate_ub}) is in fact an outer bound on the whole rate region $\mathbf{\Lambda}_s$. Furthermore, no queue can achieve a time average expected rate that is greater than the steady state probability that the corresponding queue is in ON state, i.e., $p_{01}/(p_{10}+p_{01})$. Therefore, the intersection of these $N+1$ surfaces in the $N$ dimensional space constitutes an outer bound for the rate-region $\mathbf{\Lambda}_s$. Note that this outer bound is 
tight in that the sum-rate surface of the maximum rate region $\mathbf{\Lambda}_s$, as well as the corner points $p_{01}/(p_{10}+p_{01})$ coincide with the outer bound. 
This outer bound with respect to the rate region are displayed in Fig.~\ref{Fig:stab_n3_m1eps03_ub} for the case of $N=3$ nodes.

\subsection{Myopic Policy for the Saturated System}\label{Sec:Myopic_gen}

We show in this section that a simple and intuitive policy, termed the Greedy Myopic (GM) policy, achieves the sum rate maximization for the saturated system. This policy is a greedy policy in that under the policy, if the current queue is available to serve, the server serves it.
%
%
Otherwise, the server switches to a queue with ON channel state, if such a queue exists. The policy is described in Algorithm~\ref{alg:Sat_Myopic_Pol}. Recall that $m(t)$ denotes the queue the server is present at time slot $t$.

\begin{algorithm}[h]
\caption{\emph{Greedy Myopic Policy}}
\label{alg:Sat_Myopic_Pol}
\begin{algorithmic}[1]

\STATE For all time slots $t$, if $C_{m(t)}(t)=1$, serve queue $m(t)$.

\STATE Otherwise, if $\exists j \in \{1,...,N\}, j\neq m(t)$, such that $C_j(t)=1$, among the queues that have ON channel state, switch to the queue with the smallest index in a cyclic order starting from queue $m(t)$.
\end{algorithmic}
\end{algorithm}
%
%
The cyclic switching order under the GM policy is as follows:
If the server is at queue $i$ and the decision is to switch, then the server switches to queue $j$, where
for $i=N$, $j=\arg \min_{j \in \{1,...,N-1\}}(C_j(t)=1)$ and for
$i\neq N$ if $\exists j \in \{i+1,...,N\}$, such that $C_j(t)=1$, we have $j=\arg \min_{j \in \{i+1,...,N\}}(C_j(t)=1)$, if not, then $j=\arg \min_{j \in \{1,...,i-1\}}(C_j(t)=1)$.

\vspace{3mm}
\begin{theorem}
\emph{The GM policy achieves the sum-rate upper bound.}
\end{theorem}
\begin{proof}
Given a fixed decision rule at each state, the system state forms a finite state space, irreducible and positive recurrent Markov chain. Therefore, under the GM policy, the system state converges to a steady state distribution. We partition the total probability space into three disjoint events:\\
$E_1$: the event that all the channels are in OFF state,\\
$E_2$: the event that at least 1 channel is in ON state and the server is at a queue with ON state\\
$E_3$: the event that at least 1 channel is in ON state and the server is at a queue with OFF state\\
Since these events are disjoint we have,
\begin{equation*}
1 = \mathbb{P}(E_1) + \mathbb{P}(E_2) + \mathbb{P}(E_3).
\end{equation*}
We have $\mathbb{P}(E_1) = C_0^{(N)}$ by definition. 
Since the GM policy decides to serve the current queue if it is in ON state, $\mathbb{P}(E_2)$ gives the sum throughput $\sum_ir_i$ for the GM policy. Therefore, we have that under the GM policy
\begin{equation*}
\sum_i r_i = 1 - C_0^{(N)} -  \mathbb{P}(E_3).
\end{equation*}
We show that $\mathbb{P}(E_3) = p_{10}(1-C_0^{(N)}) - p_{01}C_0^{(N)}$. Consider a time slot $t$ in steady state and let $\kappa(t)$ be the number of channels with ON states at time slot $t$ and let $E_0(t)$ be the event that the server is at a queue with OFF state at time slot $t$. We have
\begin{eqnarray*}
\mathbb{P}(E_3) \!\!\!\!\!&=& \!\!\!\!\!\mathbb{P}(E_0(t) \textrm{ and } \!1\!\le\! \kappa(t)\! \le\! N\!\!-\!\!1\!) = \!\mathbb{P}(E_0(t) \!\textrm{ and }  \!\kappa(t) \!\ge\! 1 ) \\
&=&\!\!\!\!\!\mathbb{P}(E_0(t) \textrm{ and }  \kappa(t) \ge 1 | \kappa(t-\!1)\ge1)\mathbb{P}(\kappa(t-\!1)\ge1)\\
&+& \!\!\!\!\!\mathbb{P}(E_0(t) \textrm{ and }  \kappa(t) \ge 1 | \kappa(t-\!1) = 0)\mathbb{P}(\kappa(t-\!1)=0).
\end{eqnarray*}
Since $t$ is a time slot in steady state, we have that\\$\mathbb{P}(\kappa(t-1)=0) = C_0^{(N)}$. Therefore, $\mathbb{P}(E_3)$ is given by
\begin{eqnarray*}
\!\!\!\!&\!\!\!\!\!\!\! \mathbb{P}&\!\!\!\!(\kappa(t) \!\ge\! 1| E_0(t),\kappa(t\!-\!1)\!\ge\!1\!)\mathbb{P}(E_0(t)|\kappa(t\!-\!1)\!\ge\!1\!)(1\!-\!C_0^{(N)}\!)\\
\!\!+&\!\!\!\!\!\!\mathbb{P}&\!\!\!\!(\kappa(t)\! \ge\! 1| E_0(t),\kappa(t\!-\!1)\!=\!0)\mathbb{P}(E_0(t)|\kappa(t\!-\!1)\!=\!0) \tiny{C_0^{(N)}}.
\end{eqnarray*}
We have $\mathbb{P}(E_0(t)|\kappa(t\!-\!1)\!\ge\!1\!) = p_{10}$ since the GM policy chooses a queue with ON state if there is such a queue and $\mathbb{P}(E_0(t)|\kappa(t\!-\!1)\!=0) $ is the probability that the queue chosen by the GM policy keeps its OFF channel state, given by $1-p_{01}$. 
\begin{eqnarray*}
\mathbb{P}(E_3) \!\!\!\!\!&=& \!\!\!\!\! (1-\mathbb{P}(\kappa(t) =0| E_0(t),\kappa(t\!-\!1)\!\ge\!1))p_{10}(1-C_0^{(N)}\!)\\
\!\!&+&\!\!\!\!\!(1-\mathbb{P}(\kappa(t) =0 | E_0(t),\kappa(t\!-\!1)\!=\!0))(1-p_{01}) \tiny{C_0^{(N)}}\\\\
\!\!\!\!\!&=& \!\!\!\!\! p_{10}(1-C_0^{(N)}\!)\bigg(1 \!-\! \frac{\mathbb{P}(\kappa(t) =0, E_0(t)|\kappa(t\!-\!1)\!\ge\!1)}{\mathbb{P}(E_0(t)|\kappa(t\!-\!1)\!\ge\!1\!)}\!\bigg)\\
\!\!&+&\!\!\!\!\!(1-p_{01}) \tiny{C_0^{(N)}}\bigg(1\!-\! \frac{\mathbb{P}(\kappa(t) =0, E_0(t)|\kappa(t\!-\!1)\!=\!0)}{\mathbb{P}(E_0(t)|\kappa(t\!-\!1)\!=\!0\!)}\!\bigg)\\
\!\!\!\!\!&=& \!\!\!\!\! p_{10}(1-C_0^{(N)}\!)\bigg(1\!-\! \frac{\mathbb{P}(\kappa(t) =0|\kappa(t\!-\!1)\!\ge\!1)}{p_{10}}\bigg)\\
\!\!&+&\!\!\!\!\!(1-p_{01}) \tiny{C_0^{(N)}}\bigg(1\!-\! \frac{\mathbb{P}(\kappa(t) =0|\kappa(t\!-\!1)\!=\!0)}{1-p_{01}}\bigg).
\end{eqnarray*}
We have that $\mathbb{P}(\kappa(t) =0|\kappa(t\!-\!1)\!\ge\!1)$ is given by \\
\begin{equation*}
\frac{\mathbb{P}(\kappa(t) =0) - \mathbb{P}(\kappa(t) =0|\kappa(t-1)=0)\mathbb{P}(\kappa(t-1) =0)}{\mathbb{P}(\kappa(t-1) \ge 1)},
\end{equation*}
which is equivalent to
$(C_0^{(N)} - (1-p_{01})^N C_0^{(N)})/(1-C_0^{(N)})$.
Therefore, $\mathbb{P}(E_3)$ is given by
\begin{eqnarray*}
\mathbb{P}(E_3)
\!\!\!\!\!&=& \!\!\!\!\! p_{10}(1-C_0^{(N)}\!)\bigg(1- \frac{C_0^{(N)} - (1-p_{01})^N C_0^{(N)}}{p_{10}(1-C_0^{(N)})}\bigg)\\
\!\!&+&\!\!\!\!\!(1-p_{01}) \tiny{C_0^{(N)}}\bigg(1- \frac{(1-p_{01})^N}{1-p_{01}}\bigg)\\
\!\!\!\!\!&=& \!\!\!\!\! p_{10}(1-C_0^{(N)}\!) - p_{01} C_0^{(N)}.
\end{eqnarray*}
As mentioned in the previous section, $\mathbb{P}(E_3)$ is the throughput loss due to switching as it represents the fraction of time the server is at a queue with OFF state when there are queues with ON state in the system.
\end{proof}

\subsection{Frame-Based Dynamic Control Policy}\label{Sec:FBDC_gen}

In this section we generalize the FBDC policy to the general system and show that it is throughput-optimal asymptotically in the frame length for the general case.
The FBDC algorithm for the general system is very similar to the FBDC algorithm described for two queues in Section~\ref{Sec:Opt_Pol}. Specifically, the time is divided into equal-size intervals of $T$ slots.
We find the
stationary-deterministic policy that optimally solves (\ref{eq:LP_N}) for the saturated system when
$Q_1(jT),...,Q_N(jT) $ are used as weights 
and then apply this policy in each time slot of the
frame in the actual system.
The FBDC policy is described in Algorithm \ref{alg:FBDC_Nqueues}
in details.

\begin{algorithm}[h]
\caption{\textsc{Frame Based Dynamic Control (FBDC) Policy}}
\label{alg:FBDC_Nqueues}
\begin{algorithmic}[1]

\STATE Find the optimal solution to the following LP
\beqn
\displaystyle \mbox{max.}_{\left\{\mathbf{r}\right\}} & \!\!\!\!\!\!\!\!\!\!\!\!\!\!\!\!\!\!\!\sum_{i=1}^N Q_i(jT)r_i \nonumber \\
\mbox{subject to}  &    \mathbf{r}=(r_1,...,r_N) \in \mathbf{\Lambda}_s
\label{objcon_N}
\eeqn where $\mathbf{\Lambda}_s$ is the rate region for the saturated system.
\STATE The optimal solution $(r_1^*,...,r_N^*)$ in step 1 is a corner point of $\mathbf{\Lambda}_s$ that
corresponds to a stationary-deterministic policy denoted by $\boldsymbol{\pi}^*$.
Apply $\boldsymbol{\pi}^*$ in each time slot of the frame.
\end{algorithmic}
\end{algorithm}
\vspace{0.1cm}
\begin{theorem}\label{thm:FBDC_N}
 \emph{For any $\delta>0$, there exists a large enough frame length $T$ such that the FBDC policy stabilizes the system for all arrival
rates 
within the
$\mathbf{\delta}$-stripped stability region
$\mathbf{\Lambda}_s^{\delta} = \mathbf{\Lambda}_s - \delta\mathbf{1}$.}
\end{theorem}
The proof is very similar to the proof of Theorem \ref{thm:FBDC} and is omitted. The theorem establishes the asymptotic throughput-optimality of the FBDC policy for the general system. 

\vspace{1mm}
\begin{remark}\label{Rem:FBDC_gen}
The FBDC policy provides \emph{a new framework for developing
throughput-optimal policies for network control}. Namely, given any
queuing system whose corresponding saturated system is Markovian
with finite state and action spaces, throughput-optimality is
achieved by solving an LP in order to find the stationary MDP
solution of the corresponding saturated system and applying this
solution over a frame in the actual system. In particular, \emph{the
FBDC policy can stabilize systems with arbitrary switchover times and more complicated Markov modulated
channel structures.} 
The FBDC policy can also be used to achieve throughput-optimality for
classical network control problems such as the parallel queueing systems in \cite{neely03},\cite{tass93}, 
scheduling in switches in \cite{ShahWis06} or 
scheduling under delayed channel state information \cite{shakk08}.
\end{remark}

Similar to the delay analysis in Section~\ref{Sec:FBDC_delay} for the two-queue system, a delay upper bound that is linear in the frame length $T$ can be obtained for
the FBDC policy for the the general system.
Moreover, the FBDC policy for the general system can also be implemented without any frames by setting $T=1$, i.e., by solving the LP in Algorithm~\ref{alg:FBDC_Nqueues} in each time slot. The simulation results 
regarding such implementations suggest that the FBDC policy implemented without frames
has a similar throughput performance and an improved delay performance as compared to the original FBDC policy.

\subsubsection{Discussion}

For systems with switchover delay, it is well-known that the celebrated Max-Weight scheduling policy is not throughput-optimal \cite{Berz}. In the absence of randomly varying connectivity, variable frame based generalizations of the Max-Weight policy are throughput-optimal \cite{Celik_ISIT11}. However, when the switchover delay and randomly varying connectivity are simultaneously present in the system, the FBDC policy is the only policy to achieve throughput-optimality and it has a significantly different structure from the Max-Weight policy.

The FBDC policy for a fixed frame length $T$ does not require the arrival rate information for stabilizing the system for arrival rates in $\mathbf{\Lambda} - \delta(T)\mathbf{1}$, however, it requires the knowledge of the channel connectivity parameters $p_{10},p_{01}$. 
To deal with this problem one can estimate the channel parameters periodically and use these estimates to solve the LP in (\ref{objcon_N}). This approach, of course, incurs a throughput loss depending on how large the estimation error is.

As mentioned in Remark~\ref{Rem:FBDC_gen}, the FBDC policy can stabilize a large class of network control problems whose corresponding saturated system is weakly communicating Markovian with a finite state and action spaces. However, one caveat of the FBDC policy is that the state space of the LP that needs to be solved increases exponentially with the number of links in the system.
The celebrated Max-Weight policy (which is not stabilizing for the system considered here) has linear complexity for the single server system considered in this paper. However, for general multi-server systems with $N$ servers or for a single hop network with $N$ interfering links, the Max-Weight policy has to solve a maximum-independent set problem over all links at each time slot, which is a hard problem whose state space is also exponential in the number of links $N$. The FBDC policy on the other hand, only has to solve an LP, for which there are standard solvers available such as CPLEX.
Furthermore, the FBDC policy has to solve the LP once per frame, whereas the the Max-Weight policy performs maximum-independent set computation each time slot. If the frame length for the FBDC policy is chosen to be bigger than the computational complexity of the LP in (\ref{objcon_N}), then the per-slot computational complexity of the algorithm is reduced to $O(1)$. Such a frame-based implementation is also possible for the Max-Weight policy to reduce its complexity to $O(1)$ per time slot. On the other hand, the shortcoming of such an approach for both policies is the increase in delay as a result of the larger frame length. This outlines a tradeoff between complexity and delay, whereby tuning a reduction in complexity by adjusting the frame length comes at the expense of delay.

The celebrated Max-Weight policy was first introduced in \cite{tass92} for multi-hop networks and, despite its exponential complexity in number of links, it provided a useful structure for designing queue length based scheduling algorithms. Later, this structure suggested by the Max-Weight policy lead to suboptimal but low-complexity algorithms, as well as distributed implementations of the Max-Weight policy for certain systems (see e.g., Greedy-Maximal network scheduling in \cite{WuSri06}). Our aim in proposing the FBDC policy and the state-action frequency framework for network scheduling is to give a structure for throughput-optimal algorithms for systems with time-varying channels and switchover delays, and hopefully to provide
insight into designing scalable algorithms that can stabilize such systems. The Myopic control policies we discuss in the next section constitute a first approach towards characterizing the structure of some more scalable algorithms. 

\subsection{Myopic Control Policies}\label{Sec:Myopic_Pol_N}

In this section, we generalize Myopic policies that we introduced for the two-queue system in Section~\ref{Sec:2queues} to the general system. Myopic policies make scheduling decisions based on queue lengths
and simple channel
predictions into the future.
%
%
We present an implementation of these policies over frames of length $T$ time slots where during the $j$th frame,
the queue lengths at the beginning of the frame, $Q_1(jT),...,Q_N(jT)$, are used for weight calculations during the frame. 
We describe the $1$-Lookahead Myopic (OLM) policy for the general system in 
Algorithm~\ref{alg:Myop1_N}.
%
\begin{algorithm}[h]
\caption{\textsc{$1$-Lookahead Myopic Policy}} \label{alg:Myop1_N}
\begin{algorithmic}[1]

\STATE Assuming that the server is currently with queue $1$ and the
system is at the $j$th frame, calculate the following weights in
each time slot of the current frame; \beqn
W_1(t) \!\!\!\!&=\!\!\!\!&
Q_1(jT)\Big(C_1(t)+\mathbb{E}\big[C_1(t+1)|C_1(t)\big]\Big)\nonumber\\
W_i(t) \!\!\!\!&=\!\!\!\!&
Q_i(jT)\mathbb{E}\big[C_i(t+1)|C_i(t)\big].\label{eq:myopicN}
\eeqn \STATE
If $W_1(t)\ge W_i(t),\; \forall i \in \{2,...,N\}$, then stay with queue 1. 
Otherwise, switch to a queue that achieves
\begin{equation*}
\max_i Q_i(jT)\mathbb{E}\big[C_i(t+1)|C_i(t)\big].
\end{equation*}
A similar rule applies when the server is at other queues.
\end{algorithmic}
\end{algorithm}

Similar to the FBDC policy, the Myopic policies can be implemented without the use of frames by setting $T=1$, i.e., by utilizing the current queue lengths for updating the decision rules every time slot. This could potentially lead to more delay-efficient policies that are more adaptive to dynamic changes in queue lengths. We elaborate on this via the numerical results in the next section.

These policies have very low complexity and they are simpler to implement as compared to the FBDC policy.
As suggested by the simulation results in Section~\ref{Sec:Sim}, 
the Myopic
policies may achieve the full stability region while
providing better delay performance than
the FBDC policy in most cases. 

\subsection{Myopic Control Policies}\label{Sec:Myopic_Pol_N}

In this section, we generalize Myopic policies that we introduced for the two-queue system in Section~\ref{Sec:2queues} to the general system. Myopic policies make scheduling decisions based on queue lengths
and simple channel
predictions into the future.
%
%
We present an implementation of these policies over frames of length $T$ time slots where during the $j$th frame,
the queue lengths at the beginning of the frame, $Q_1(jT),...,Q_N(jT)$, are used for weight calculations during the frame. 
We describe the $1$-Lookahead Myopic (OLM) policy for the general system in 
Algorithm~\ref{alg:Myop1_N}.
%
\begin{algorithm}[h]
\caption{\textsc{$1$-Lookahead Myopic Policy}} \label{alg:Myop1_N}
\begin{algorithmic}[1]

\STATE Assuming that the server is currently with queue $1$ and the
system is at the $j$th frame, calculate the following weights in
each time slot of the current frame; \beqn
W_1(t) \!\!\!\!&=\!\!\!\!&
Q_1(jT)\Big(C_1(t)+\mathbb{E}\big[C_1(t+1)|C_1(t)\big]\Big)\nonumber\\
W_i(t) \!\!\!\!&=\!\!\!\!&
Q_i(jT)\mathbb{E}\big[C_i(t+1)|C_i(t)\big].\label{eq:myopicN}
\eeqn \STATE
If $W_1(t)\ge W_i(t),\; \forall i \in \{2,...,N\}$, then stay with queue 1. 
Otherwise, switch to a queue that achieves
\begin{equation*}
\max_i Q_i(jT)\mathbb{E}\big[C_i(t+1)|C_i(t)\big].
\end{equation*}
A similar rule applies when the server is at other queues.
\end{algorithmic}
\end{algorithm}

The technique used for the case of two queues for analyzing the stability region achieved by the OLM policy is extremely cumbersome to generalize to the general system with $N$ queues. Therefore, for the general system, we have investigated the performance of the OLM policy in simulations. The simulation results in Section~\ref{Sec:Sim} suggest that 
the OLM policy may achieve the full stability region while
providing a better delay performance as compared to
the FBDC policy. 

Similar to the FBDC policy, the Myopic policies can be implemented without the use of frames by setting $T=1$, i.e., by utilizing the current queue lengths for updating the decision rules every time slot. This could potentially lead to more delay-efficient policies that are more adaptive to dynamic changes in queue lengths. We elaborate on this via the numerical results in the next section.

%

Similar to the system with two queues, the $k$-Lookahead Myopic Policy is the same as
before except that the following weight functions are used for
scheduling decisions: Assuming the server is with queue 1 at time
slot $t$,\\
$W_1(t)
=Q_1(jT)\big(C_1(t)+\sum_{\tau=1}^k\mathbb{E}[C_1(t+\tau)|C_1(t)]\big)$
and $W_i(t) =Q_i(jT)\sum_{\tau=1}^k\mathbb{E}[C_i(t+\tau)|C_i(t)],\; i \in \{2,...,N\}$.\\
These policies have very low complexity and they are simpler to implement as compared to the FBDC policy.

\section{Numerical Results}\label{Sec:Sim}

We performed simulation experiments that present average queue
occupancy results for the FBDC, the One-Lookahead Myopic (OLM) and the Max-Weight (MW) policies for systems with $N=2$ or $N=3$ queues. We first
verified that in the simulation results for the FBDC policy, queue sizes grow unbounded only for arrival rates outside the stability region, 
and then performed experiments for the 1-Lookahead Myopic (OLM) policy. 
In all the reported results, we have $\boldsymbol{\lambda} \in \mathbf{\Lambda}$ with $0.01$ increments.
For each point at the boundary of $\mathbf{\Lambda}$, we simulated one
point outside the stability region. Furthermore, for each
data point, the arrival processes were i.i.d., the channel processes
were Markovian as in Fig. \ref{Fig:channel_MC} and the simulation
length was
$T_s=100,000$ slots. 

\begin{figure}
\centering
\includegraphics[width=0.4\textwidth]{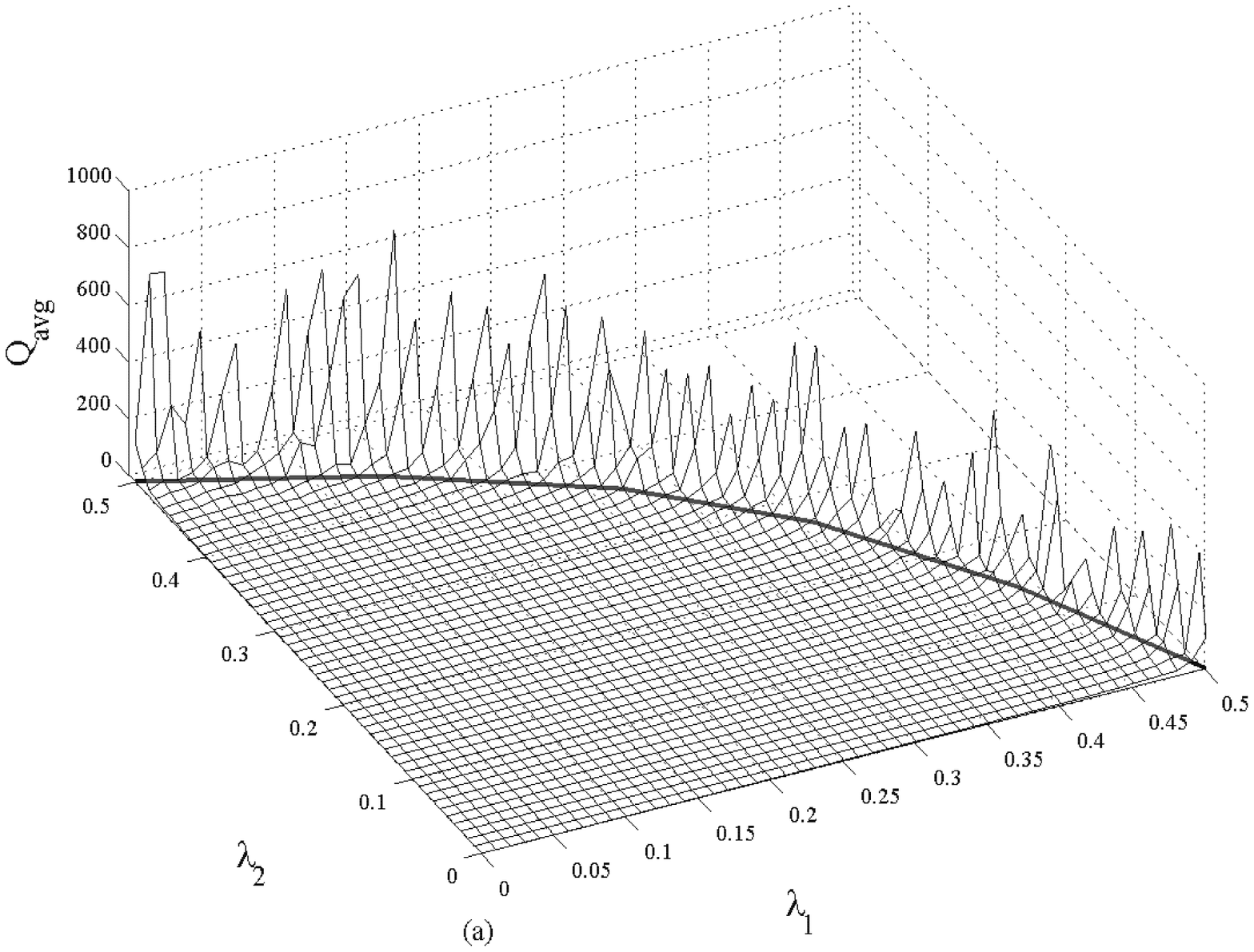}\\
\includegraphics[width=0.4\textwidth]{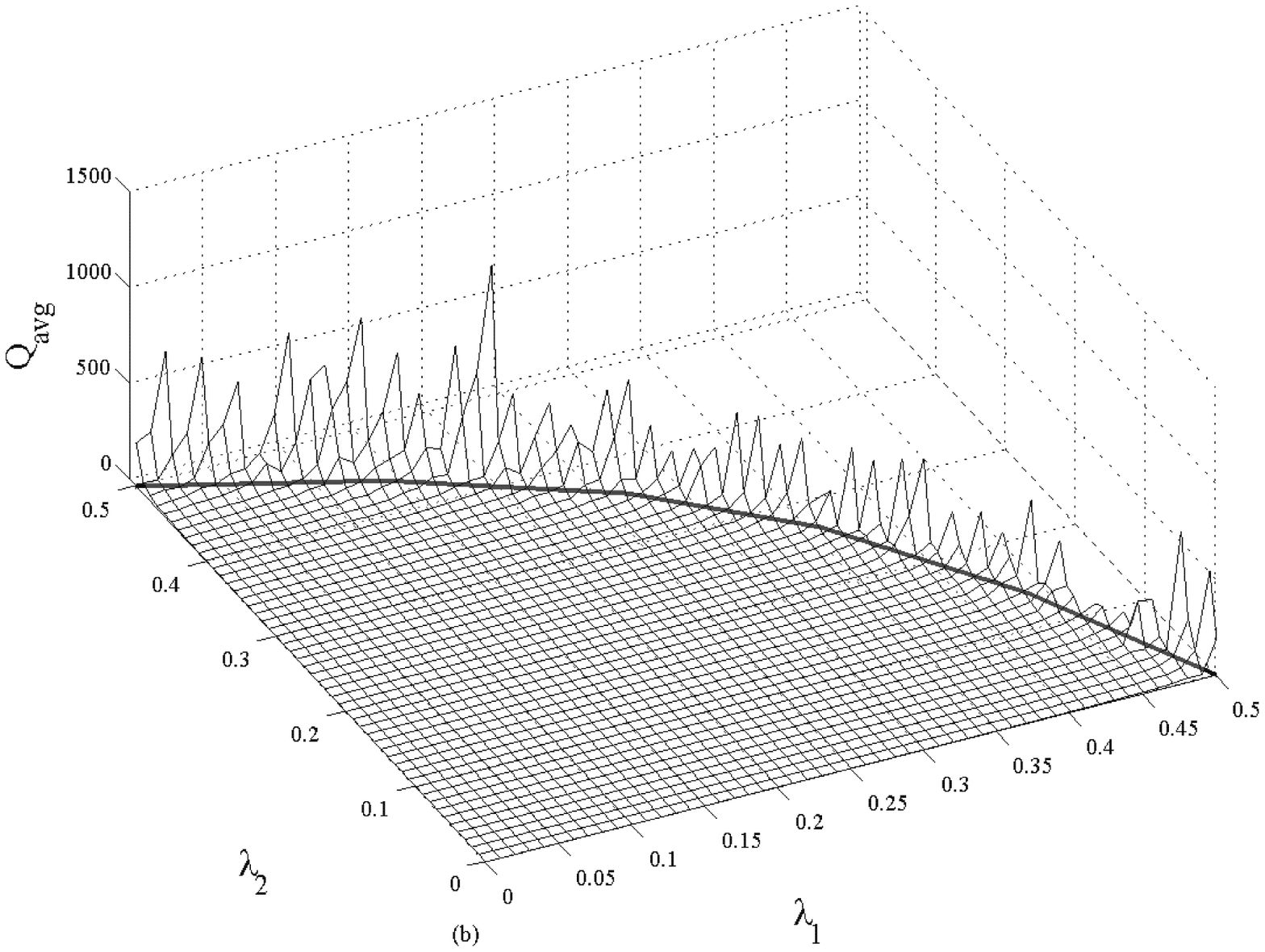}
\captionsetup{font=footnotesize}\caption{The total
average queue size for (a) the FBDC policy and (b) the Myopic policy
for $T=25$ and $\epsilon = 0.25$.}\label{Fig:FBDCvsMyopicT25eps025}
\end{figure}

Fig.\ \ref{Fig:FBDCvsMyopicT25eps025}~(a) presents the total average
queue size, $\mathbf{Q}_{avg}\doteq\sum_{t=1}^{T_s}(Q_1(t)+Q_2(t))/T_s$,
under the FBDC policy for $N=2$ queues, $\epsilon=0.25<\epsilon_c$, and a frame size of $T=25$ slots.
The boundary
of the stability region is shown by (red) lines on the two dimensional
$\lambda_1-\lambda_2$ plane. We observe that the average queue sizes
are small for all $(\lambda_1,\lambda_2) \in \mathbf{\Lambda}_s$ and
the big jumps in queue sizes occur for points outside
$\mathbf{\Lambda}$.
Fig.\ \ref{Fig:FBDCvsMyopicT25eps025}~(b) presents the performance of the OLM policy with $T=25$ slot frames for the
same system.
The simulation results suggest that 
there is no appreciable difference between the stability regions of the
FBDC and the OLM policies. Note that the total average queue size is proportional to long-run packet-average delay in the system through Little's law. For these two
figures, the average delay under the OLM policy is less than that under the FBDC policy 
for $81\%$ of all arrival rates considered.

\begin{figure}
\centering
\includegraphics[width=0.4\textwidth]{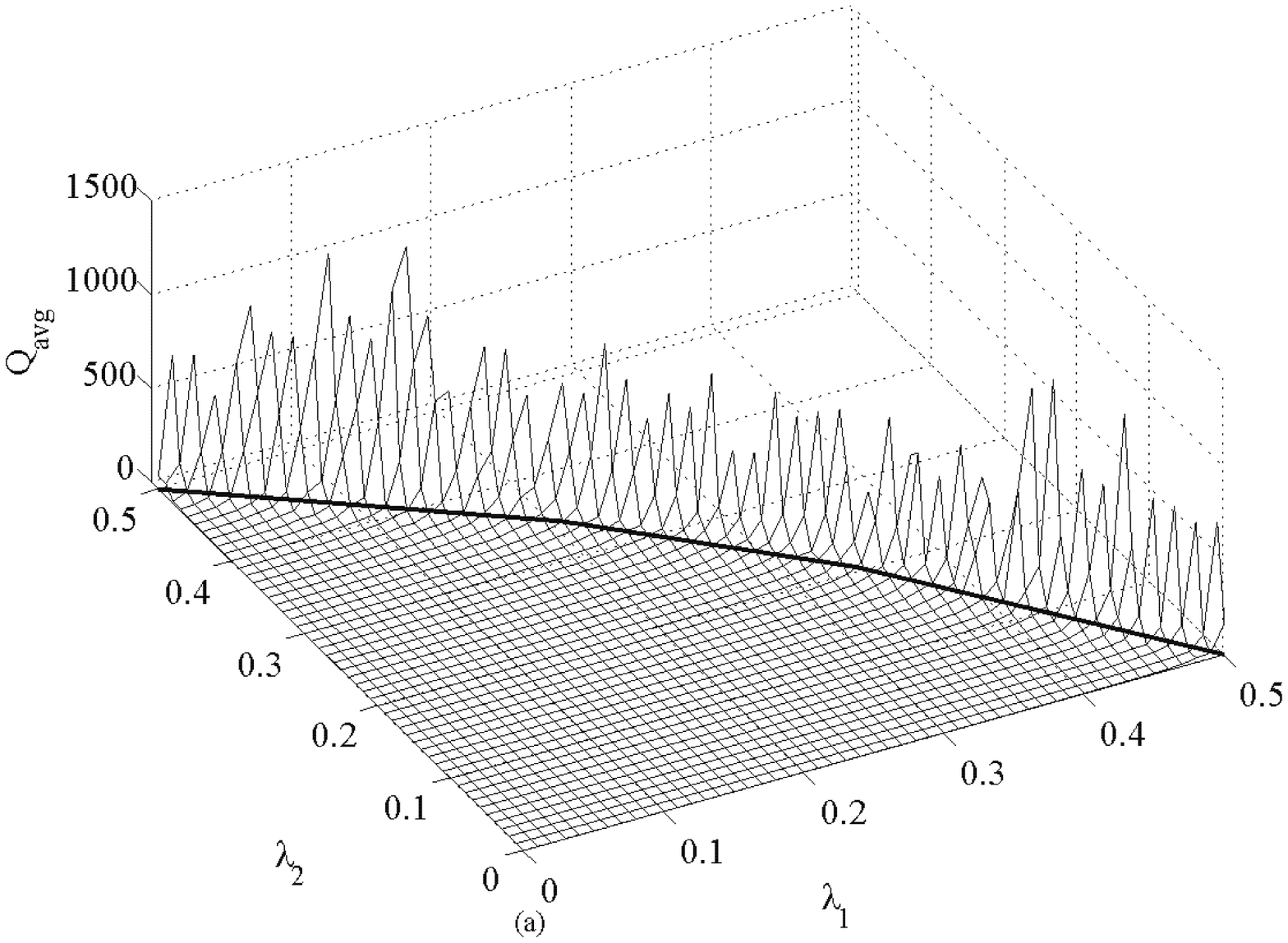}\\
\includegraphics[width=0.4\textwidth]{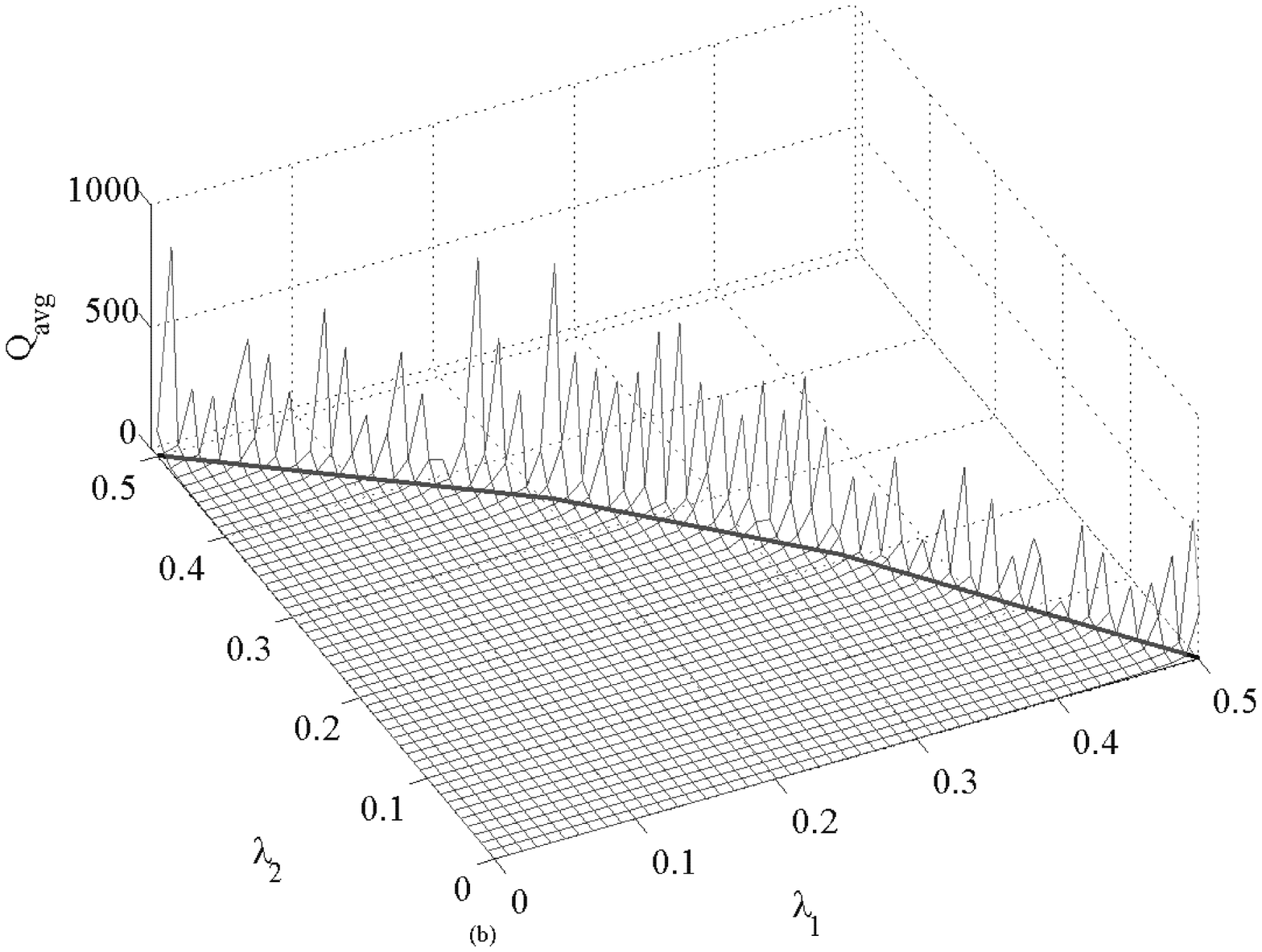}
\captionsetup{font=footnotesize}\caption{The total
average queue size for (a) the FBDC policy and (b) the Myopic Policy implemented without the use of frames (i.e., for $T=1$) and $\epsilon = 0.40$.}\label{Fig:FBDCvsMyopicT1eps040}
\end{figure}

Next, we implemented the FBDC and the OLM policies without the use of any frames (i.e., for $T=1$). When there are no frames, the FBDC policy solves the LP in Algorithm~\ref{alg:FBDC} in each time slot, and the OLM policy utilizes the queue length information in the
\emph{current} time slot for the weight calculations in
(\ref{eq:myopic}).
Fig.\ \ref{Fig:FBDCvsMyopicT1eps040}~(a) and (b) present the total average
queue size under the FBDC and the OLM policies for $N=2$ queues, $T=1$, and $\epsilon=0.40>\epsilon_c$. Similar to the frame based implementations, we observe that the average queue sizes
are small for all $(\lambda_1,\lambda_2) \in \mathbf{\Lambda}$ for both policies and
the big jumps in queue sizes occur for points outside
$\mathbf{\Lambda}$, which suggests that the the non-frame based implementation of the FBDC and the OLM
policies may achieve the full stability region. The reason why the FBDC and the OLM policies provide stability without the use of frames is because for large queue lengths, the corner point that these policies choose to apply depend completely on the queue length ratios, and hence, the choice of corner points and the associated saturated-system policies utilized in the FBDC and the OLM policies do not change fast when the queue lengths get large. Furthermore, the no-frame implementations of these policies are more adaptive to dynamic changes in the queue sizes as compared to implementations with large frames.

\begin{figure}
\centering
\includegraphics[width=0.4\textwidth]{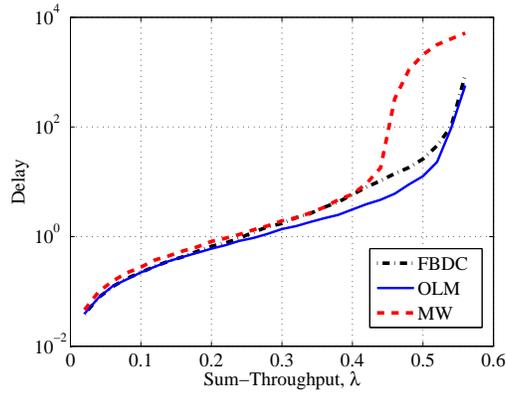}
\captionsetup{font=footnotesize}\caption{Delay vs Sum-throughput for the FBDC, the OLM, and the Max-Weight policies implemented without the use of frames (i.e., for $T=1$) for $N=2$ queues and $\epsilon = 0.40$.}\label{Fig:Delay_FBDC_OLM_MW_T1_eps040}
\end{figure}

For the same system (i.e., $N=2$ queues and $\epsilon=0.40>\epsilon_c$), Fig.~\ref{Fig:Delay_FBDC_OLM_MW_T1_eps040} presents the long-run packet average delay as a function of the sum-throughput $\lambda_1 + \lambda_2$ along the main diagonal line (i.e., $\lambda_1=\lambda_2$). We compare the delay performance of the FBDC and the OLM policies with $T=1$, and the Max-Weight policy which, in each time slot $t$, chooses the queue that achieves $\max_i Q_i(t)C_i(t)$.  The maximum sum-throughput is $0.75-\epsilon/2 = 0.55$ as suggested by Theorem~\ref{thm:stab}. Fig.~\ref{Fig:Delay_FBDC_OLM_MW_T1_eps040} shows that while FBDC and the OLM policies stabilize the system for all $ \lambda_1+\lambda_2< 0.55$, the system becomes unstable under the Max-Weight policy around $\lambda_1+\lambda_2 =0.45$. This result also confirms that the OLM policy has a much better delay performance than the FBDC and the Max-Weight policies.

\begin{figure}
\centering
\includegraphics[width=0.4\textwidth]{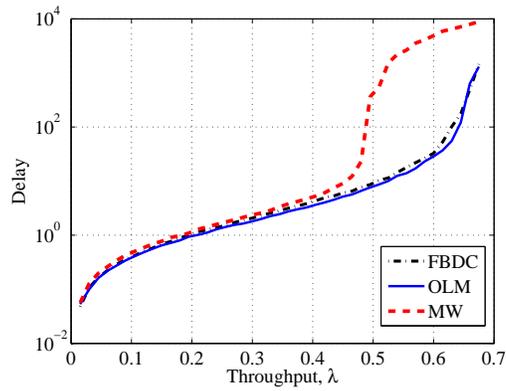}
\captionsetup{font=footnotesize}\caption{Delay vs Sum-throughput for the FBDC, the OLM, and the Max-Weight policies implemented without the use of frames (i.e., for $T=1$) for $N=3$ queues and $\epsilon = 0.30$.}\label{Fig:Delay_N3_FBDC_OLM_MW_T1_eps030}
\end{figure}

For $N=3$ queues and $\epsilon=0.30$, Fig.~\ref{Fig:Delay_N3_FBDC_OLM_MW_T1_eps030} presents the long-run packet-average delay as a function of the sum-throughput $\sum_i\lambda_i$ along the main diagonal line (i.e., $\lambda_1=\lambda_2=\lambda_3$).
The maximum sum-throughput is $1-0.5^N - \epsilon*(1-0.5^N) + \epsilon*0.5^N = 0.65$
$1-\epsilon - 0.5^N(1 - 2\epsilon) = 0.65$
as suggested by Lemma~\ref{lemma:sum_tp_Nqueues}. Similar to the previous case, Fig.~\ref{Fig:Delay_N3_FBDC_OLM_MW_T1_eps030} shows that FBDC and the OLM policies stabilize the system for all $ \sum_i \lambda_i< 0.65$, the system becomes unstable under the Max-Weight policy around $\sum_i\lambda_i=0.48$. This result also confirm that the OLM policy has a delay performance than the other two policies.

The delay results in this section show that the OLM policy is not only simpler to
implement as compared to the FBDC policy, but it can also be more delay
efficient.

\section{Conclusions}\label{Sec:Conc}

We investigated the dynamic server allocation problem with
\emph{randomly varying connectivity} and \emph{server switchover
time}. For the case of two queues, we analytically characterized the stability region of the system using state-action frequencies that are stationary solutions to an MDP formulation for the corresponding saturated system.
We developed the throughput-optimal FBDC policy. We also developed simple \emph{Myopic Policies} that achieve
a large fraction of the stability region. We extended the stability region characterization in terms of state-action frequencies and the throughput-optimality of the FBDC policy to the general system with arbitrary number of queues. We characterized tight analytical outer bounds on the stability region using an upper bound on the sum-rate and showed that a simple greedy-myopic policy achieves this sum-rate bound. 
The stability region
characterization in terms of the state-action frequencies of the
saturated system and the throughput-optimality of the FBDC policy
hold for systems with arbitrary
switchover times and general Markovian channels. Furthermore, the
FBDC policy provides a new framework for developing
throughput-optimal policies for network control as this policy can be used to stabilize a large class of other network control problems. 

In the future, we intend to 
explicitly characterize the
stability region 
of systems with multiple-slot
switchover times with general Markov modulated channels. 
We intend to develop throughput-optimal Myopic policies 
for general system models. 
Finally, joint scheduling and routing in multihop networks with
dynamically changing channels and server switchover times is a challenging future direction.



\section*{Appendix A-Proof of Theorem \ref{thm:stab_iid} }\label{Sec:App-A}
%
%

We prove Theorem \ref{thm:stab_iid} for a more general system with
$N$-queues and travel time between queue-$i$ and queue-$j$ given by
$D_{ij}$ slots. We call the term $\sum_{i=1}^N \lambda_i/p_i$ the
system load and denote it by $\rho$ since it is the rate with which
the work is entering the system in the form of service slots.\\\\
\textit{Necessity}\\
We prove that a necessary condition for the stability of any policy is $\rho
=\sum_{i=1}^N \lambda_i/p_i<1$.
\begin{proof}
Since queues have memoryless channels, for any received packet, as
soon as the server switches to some queue $i$, the expected time to ON
state is $1/p_i$. Namely, the time to ON state is a geometric random
variable with parameter $p_i$, and hence, $1/p_i$ is the
``service time per packet'' for queue $i$. 
%
%
In
a multiuser single-server system \emph{with or without switchover
delays}, with i.i.d. arrivals whose average arrival rates are
$\lambda_i, i \in \{1,...,N\}$, and i.i.d. service times independent of
arrivals with mean $1/p_i, i \in \{1,...,N\}$, a
necessary condition for stability is given by the system load,
$\rho$, less than 1 \cite{Boxma90}.
\end{proof}
To futher ellaborate on this, consider the polling system with
zero switchover times, i.i.d. arrivals of mean $\lambda_i$ and
i.i.d. service times of mean $1/p_i$. The throughput region of
this system is an upperbound on the throughput region of the
corresponding system with nonzero switchover times. This is because for the same
sample path of arrival and service processes, the system with zero
switchover time can achieve exactly the same departure process as
the system with nonzero switchover times by making the server idle
when necessary. A necessary condition for the stability of the
former system is $\rho = \lambda_1/p_1 +
\lambda_N/p_N+...+\lambda_1/p_N <1$, (e.g., \cite{Walrand}, \cite{Boxma90}). Finally, note that this necessary stability condition can also be
derived by utilizing the state-action frequency approach of Section \ref{Sec:Stab_Reg} for the system with i.i.d. connectivity processes.\\\\
\textit{Sufficiency}
\begin{proof}
Under the Gated cyclic policy, we have a Polling system with i.i.d. arrivals with mean
$\lambda_i, i \in \{1,...,N\}$, i.i.d. service times independent of
arrivals with mean $1/p_i, i \in \{1,...,N\}$, and finite and constant switchover delays. It is shown in \cite{Boxma90}
that the Gated cyclic policy results in an ergodic system if $\rho
=\sum_{i=1}^N \frac{\lambda_i}{p_i} <1$, the expected per-message waiting times in steady-state are finite, and they satisfy a pseudo-conversation law. Through Little's law (\cite[pp. 139]{Gallager} or \cite[pp. 1109]{AltmanLevy94}), this implies that the expected number of packets in the system in steady state is finite, which in turn implies that the system is stable.
\end{proof}

\section*{Appendix B-Proof of Theorem \ref{thm:stab}}\label{Sec:App-C}

We enumerate the states as follows:
\begin{equation}\label{eq:state_enum}
\begin{array}{llll}
s=(1,1,1)\equiv1,\;\;\; &s=(1,1,0)\equiv2,\;\;\;
&s=(1,0,1)\equiv3,\;\;\;  & s=(1,0,0)\equiv4,\;\;\;\\
s=(2,1,1) \equiv5 & s=(2,1,0)\equiv6 & s=(2,0,1) \equiv7& s=(2,0,0) \equiv8.\\
\end{array}
\end{equation}
We rewrite the balance equations in (\ref{eq:LP}) in more details.
\beqn x(1;1)+x(1;0) = (1-\epsilon)^2 \big(x(1;1)+x(5;0)\big)
\!\!\!\!\!&+&
\!\!\!\!\!\epsilon(1-\epsilon)\big(x(2;1)+ x(6;0)\big)\nonumber\\
+ \epsilon(1-\epsilon)\big(x(3;1)+x(7;0)\big) \!\!\!\!\!&+&\!\!\!\!\! \epsilon^2\big(x(4;1)+x(8;0)\big)\label{eq:saf1}  \\
x(2;1)+x(2;0) = \epsilon(1-\epsilon) \big(x(1;1)+x(5;0)\big)
\!\!\!\!\!&+&\!\!\!\!\!
(1-\epsilon)^2\big(x(2;1) + x(6;0)\big) \nonumber\\
\!\!\!\!\!+ \epsilon^2\big(x(3;1)+x(7;0)\big)
\!\!\!\!\!&+&\!\!\!\!\!
\epsilon(1-\epsilon)\big(x(4;1)+x(8;0)\big)\label{eq:saf2}\\
\!\!\!\!\!&\ldots& \nonumber\\
x(5;1)+x(5;0) = (1-\epsilon)^2 \big(x(5;1)+x(1;0)\big)
\!\!\!\!\!\!\!\!&+\!\!\!\!\!&
\epsilon(1-\epsilon)\big(x(6;1)+ x(2;0)\big)\nonumber\\
+ \epsilon(1-\epsilon)\big(x(7;1)+x(3;0)\big) \!\!\!\!\!&+&\!\!\!\!\! \epsilon^2\big(x(8;1)+x(4;0)\big) \label{eq:saf3} \\
x(7;1)+x(7;0)= \epsilon(1-\epsilon) \big(x(5;1)+x(1;0)\big)
\!\!\!\!\!&+&\!\!\!\!\!
\epsilon^2\big(x(6;1) + x(2;0)\big) \nonumber\\
+ (1-\epsilon)^2 \big(x(7;1)+x(3;0)\big) \!\!\!\!\!&+&\!\!\!\!\!
\epsilon(1-\epsilon)\big(x(8;1)+x(4;0)\big)\label{eq:saf4}\\
&\ldots\nonumber  \eeqn\\
The following equations hold the channel state pairs
$(C_1,C_2)$. \beqn
x(1;1)+x(1;0)+ x(5;1)+x(5;0) = 1/4 \label{eq:state_1}\\
x(2;1)+x(2;0)+ x(6;1)+x(6;0) = 1/4 \label{eq:state_2}\\
x(3;1)+x(3;0)+ x(7;1)+x(7;0) = 1/4 \label{eq:state_3}\\
x(4;1)+x(4;0)+ x(8;1)+x(8;0) = 1/4 \label{eq:state_4}.
 \eeqn\\
 \\
Let $u_1=(x(1;1)+x(2;1)$ and $u_2=(x(5;1)+x(7;1))$. Summing up
(\ref{eq:saf1}) with (\ref{eq:saf2}) and (\ref{eq:saf3}) with
(\ref{eq:saf4}) we have
\begin{eqnarray*}
\epsilon u_1 \!\!\!\!&=&\!\!\!\! - \big(x(1;0)+x(2;0)\big) +
\epsilon\big( x(3;1)+x(4;1) \big) + \epsilon\big(x(7;0)+ x(8;0)
\big) + (1-\epsilon) \big( x(5;0)+x(6;0) \big)\\
\epsilon u_2 \!\!\!\!&=&\!\!\!\! - \big(x(5;0)+x(7;0)\big) +
\epsilon\big( x(6;1)+x(8;1) \big) + \epsilon\big(x(2;0)+ x(4;0)
\big) + (1-\epsilon) \big( x(1;0)+x(3;0) \big)\\
\end{eqnarray*}
\\
Rearranging and using (\ref{eq:state_1})-(\ref{eq:state_4}) we have
\begin{eqnarray}
u_1 \!\!\!\!&=&\!\!\!\! \frac{1-\epsilon}{2}+\epsilon\big(
x(3;1)+x(4;1)+x(7;0)+ x(8;0)\big) - (2-\epsilon) \big( x(1;0)+x(2;0) \big) - (1-\epsilon) \big( x(5;1)+x(6;1)\big)\nonumber\\
\label{eq:saf_work1}\\
u_2 \!\!\!\!&=&\!\!\!\! \frac{2-\epsilon}{4}+\epsilon\big(
x(2;0)-x(4;1) + x(6;1)-x(8;0) \big) - (2-\epsilon) \big(
x(5;0)+x(7;0) \big) - (1-\epsilon) \big(
x(1;1)+x(3;1)\big)\nonumber\\
\label{eq:saf_work2}
\end{eqnarray}
\\
\\
Using (\ref{eq:saf3}) in (\ref{eq:saf_work1})  and (\ref{eq:saf1})
in (\ref{eq:saf_work2}) we have
\begin{eqnarray}
u_1 \!\!\!\!&=&\!\!\!\! \frac{1-\epsilon}{2}+\epsilon\big(
x(3;1)+x(4;1)+x(7;0)+ x(8;0)\big) -
\frac{\epsilon(1-\epsilon)}{2-\epsilon} \big( x(4;0)+x(8;1) \big)
-\frac{(1-\epsilon)(3-2\epsilon)}{2-\epsilon} x(6;1)\nonumber \\
&+& \frac{1-\epsilon}{\epsilon(2-\epsilon)} x(5;0)
-\frac{1+\epsilon-\epsilon^2}{\epsilon(2-\epsilon)} x(1;0) -
\frac{(1-\epsilon)^2}{2-\epsilon} \big( x(3;0)+x(7;1) \big)
-\Big( 2-\epsilon + \frac{(1-\epsilon)^2}{2-\epsilon}\Big) x(2;0)\label{eq:saf_work3} \\
u_2 \!\!\!\!&=&\!\!\!\! \frac{2-\epsilon}{4}+\epsilon\big( x(2;0)+
x(6;1)\big) - \Big(
\epsilon+\frac{\epsilon(1-\epsilon)}{2-\epsilon}\Big) \big(
x(4;1)+x(8;0) \big)
-\frac{(1-\epsilon)(3-2\epsilon)}{2-\epsilon} x(3;1)\nonumber \\
&+& \frac{1-\epsilon}{\epsilon(2-\epsilon)} x(1;0)
-\frac{1+\epsilon-\epsilon^2}{\epsilon(2-\epsilon)} x(5;0) -
\frac{(1-\epsilon)^2}{2-\epsilon} \big( x(2;1)+x(6;0) \big) -\Big(
2-\epsilon + \frac{(1-\epsilon)^2}{2-\epsilon}\Big)
x(7;0).\label{eq:saf_work4}
\end{eqnarray}
\\
\\
Using (\ref{eq:state_3}) and (\ref{eq:state_4}) in
(\ref{eq:saf_work3}) and (\ref{eq:state_2}) in (\ref{eq:saf_work4})
we have
\begin{eqnarray}
u_1 \!\!\!\!&=&\!\!\!\!
\frac{(1-\epsilon)(3-2\epsilon)}{4(2-\epsilon)}+ \Big(\epsilon+
\frac{\epsilon(1-\epsilon)}{2-\epsilon} \Big) \big(
x(4;1)+x(8;0)\big) +  \frac{1}{2-\epsilon} \big(x(3;1)+x(7;0)\big)\nonumber\\
&-&\frac{(1-\epsilon)(3-2\epsilon)}{2-\epsilon} x(6;1)  +
\frac{1-\epsilon}{\epsilon(2-\epsilon)} x(5;0)
-\frac{1+\epsilon-\epsilon^2}{\epsilon(2-\epsilon)} x(1;0)
-\Big( 2-\epsilon + \frac{(1-\epsilon)^2}{2-\epsilon}\Big) x(2;0)\label{eq:saf_work5} \\
u_2 \!\!\!\!&=&\!\!\!\! \frac{3-2\epsilon}{4(2-\epsilon)}
-\Big(\epsilon+ \frac{\epsilon(1-\epsilon)}{2-\epsilon} \Big) \big(
x(4;1)+x(8;0)\big) + \frac{1}{2-\epsilon} \big(x(2;0)+x(6;1)\big)\nonumber\\
&-&\frac{(1-\epsilon)(3-2\epsilon)}{2-\epsilon} x(3;1) +
\frac{1-\epsilon}{\epsilon(2-\epsilon)} x(1;0) -
\frac{1+\epsilon-\epsilon^2}{\epsilon(2-\epsilon)} x(5;0) - \Big(
2-\epsilon + \frac{(1-\epsilon)^2}{2-\epsilon}\Big)
x(7;0).\label{eq:saf_work6}
\end{eqnarray}
\\

Consider the LP objective function $\alpha_1 \big(x(1;1)+x(2;1)\big)
+ \alpha_2\big(x(5;1)+x(7;1)\big)$, and note that the solution to
this LP is a stationary deterministic policy for any given
$\alpha_1$ and $\alpha_2$. This means that, for any state $s$ either
$x(s;1)$ or $x(s;0)$ has to be zero. In order to maximize $\alpha_1
\big(x(1;1)+x(2;1)\big) + \alpha_2\big(x(5;1)+x(7;1)\big)$ we need

\begin{eqnarray*}
x(7;0) \!\!\!\!&=&\!\!\!\! 0 \;\;\;\;\textrm{if  }\;\; \frac{\alpha_2}{\alpha_1} \ge \frac{1}{(2-\epsilon)^2+(1-\epsilon)^2},\\
x(3;1) \!\!\!\!&=&\!\!\!\! 0 \;\;\;\;\textrm{if  }\;\;\frac{\alpha_2}{\alpha_1} \ge \frac{1}{(1-\epsilon)(3-2\epsilon)},\\
x(5;0) \!\!\!\!&=&\!\!\!\! 0 \;\;\;\;\textrm{if  }\;\;\frac{\alpha_2}{\alpha_1} \ge \frac{1-\epsilon}{1+\epsilon-\epsilon^2},\\
x(8;0)=x(4;1)\!\!\!\!&=&\!\!\!\! 0 \;\;\;\;\textrm{if  }\;\;\frac{\alpha_2}{\alpha_1} \ge 1,\\
x(6;0) \!\!\!\!&=&\!\!\!\! 0 \;\;\;\;\textrm{if  }\;\;\frac{\alpha_2}{\alpha_1} \ge (1-\epsilon)(3-2\epsilon),\\
x(1;1) \!\!\!\!&=&\!\!\!\! 0
\;\;\;\;\textrm{if}\;\;\frac{\alpha_2}{\alpha_1} \ge
\frac{1+\epsilon-\epsilon^2}{1-\epsilon}.
\end{eqnarray*}
\\
Note that we have
\begin{eqnarray*}
(2-\epsilon)^2+(1-\epsilon)^2 &\ge& (1-\epsilon)(3-2\epsilon) \ge 1\\
(2-\epsilon)^2+(1-\epsilon)^2 &\ge& \frac{1+\epsilon-\epsilon^2}{1-\epsilon} \ge 1\\
\end{eqnarray*}
holding for all $\epsilon \in [0,0.5]$. Consider the following two
cases:
\\
\\
\textbf{Case-1: $\epsilon > \epsilon_c = 1-\sqrt{2}/2$}
\\
In this case we have $(1-\epsilon)(3-2\epsilon) <
(1+\epsilon-\epsilon^2)/(1-\epsilon)$. This means that we have the
following optimal
policies depending on the value of $\alpha_2/\alpha_1$.\\
\\
$1 \le \frac{\alpha_2}{\alpha_1} \le (1-\epsilon)(3-2\epsilon)$:
\begin{equation*}
\begin{array}{llll}
\textrm{@queue 1}: (1,1,1):\textrm{stay}, &(1,1,0):\textrm{stay}, &(1,0,1):\textrm{switch}, &(1,0,0):\textrm{switch}. \\
\textrm{@queue 2}: (2,1,1):\textrm{stay}, &(2,1,0):\textrm{switch},
&(2,0,1):\textrm{stay}, &(2,0,0):\textrm{stay}.
\end{array}
\end{equation*}
Substituting zero values for the state action pairs that are not chosen into (\ref{eq:saf_work5}) and
(\ref{eq:saf_work6}), it can be seen that this policy achieves the
rate pair
\begin{equation*}
r_1= \frac{(1-\epsilon)(3-2\epsilon)}{4(2-\epsilon)},\;\;\;\;\; r_2
=\frac{3-2\epsilon}{4(2-\epsilon)}.
\end{equation*}
\\
\\
$\frac{\alpha_2}{\alpha_1} > (1-\epsilon)(3-2\epsilon)$:
\begin{equation*}
\begin{array}{llll}
\textrm{@queue 2}: (2,1,1):\textrm{stay}, &(2,1,0):\textrm{stay},
&(2,0,1):\textrm{stay}, &(2,0,0):\textrm{stay}.
\end{array}
\end{equation*}
In this case it is optimal to stay at queue 2 for all channel
conditions. The decisions at queue 1 are to switch to queue 2.
Namely, it is sufficient that at least one state corresponding to
server being at queue 1 to take a switch decision, which is the case
for $\alpha_2/\alpha_1\ge ((1-\epsilon)(3-2\epsilon))$, since $x(3;1)
=0$ if $\alpha_2/\alpha_1\ge 1/((1-\epsilon)(3-2\epsilon))$. Since
the policy decides to always stay at queue 2, it achieves the rate pair
\begin{equation*}
r_1= 0,\;\;\;\;\; r_2 =0.5.
\end{equation*}
Note that the case for $\alpha_2/\alpha_1 < 1$ is symmetric and can
be obtained similarly.\\
\\
\\
\textbf{Case-2: $\epsilon < \epsilon_c = 1-\sqrt{2}/2$}\\
In this case we have $(1-\epsilon)(3-2\epsilon) >
(1+\epsilon-\epsilon^2)/(1-\epsilon)$. This means that before the
state $x(6;0)$ becomes zero, namely for
$(1+\epsilon-\epsilon^2)/(1-\epsilon) < \alpha_2/\alpha_1 <
(1-\epsilon)(3-2\epsilon)$, having $x(1;1)=0$ is optimal. This means
that there is one more corner point of the rate region for $\epsilon
< \epsilon_c$. We have the following optimal
policies.
\\
\\
$1 \le \frac{\alpha_2}{\alpha_1} \le
\frac{1+\epsilon-\epsilon^2}{1-\epsilon}$:\\
\begin{equation*}
\begin{array}{llll}
\textrm{@queue 1}: (1,1,1):\textrm{stay}, &(1,1,0):\textrm{stay}, &(1,0,1):\textrm{switch}, &(1,0,0):\textrm{switch}. \\
\textrm{@queue 2}: (2,1,1):\textrm{stay}, &(2,1,0):\textrm{switch},
&(2,0,1):\textrm{stay}, &(2,0,0):\textrm{stay}.
\end{array}
\end{equation*}
This policy is the same policy as in the previous case and it
achieves the rate pair
\begin{equation*}
r_1= \frac{(1-\epsilon)(3-2\epsilon)}{4(2-\epsilon)},\;\;\;\;\; r_2
=\frac{3-2\epsilon}{4(2-\epsilon)}.
\end{equation*}
\\
\\
$ \frac{\alpha_2}{\alpha_1}
> \frac{1+\epsilon-\epsilon^2}{1-\epsilon}$:\\
\\
We have the following deterministic actions.
\begin{equation*}
\begin{array}{llll}
\textrm{@queue 1}: (1,1,1):\textrm{switch}, &(1,1,0):\textrm{?}, &(1,0,1):\textrm{switch}, &(1,0,0):\textrm{switch}. \\
\textrm{@queue 2}: (2,1,1):\textrm{stay}, &(2,1,0):\textrm{?},
&(2,0,1):\textrm{stay}, &(2,0,0):\textrm{stay}.
\end{array}
\end{equation*}
In order to find the final threshold on $\alpha_2/\alpha_1$, we
substitute the above deterministic decisions in (\ref{eq:saf2}),
(\ref{eq:saf3}) and (\ref{eq:saf4}). Utilizing also
(\ref{eq:state_1}), (\ref{eq:state_2}), (\ref{eq:state_3}) and
(\ref{eq:state_4}) we obtain
\begin{eqnarray}
x(2;1) \!\!\!\!&=&\!\!\!\! \frac{(1-\epsilon)^2}{4} - (1-\epsilon)^2 x(6;1)\label{eq:saf_work7}\\
x(5;1)+x(7;1) \!\!\!\!&=&\!\!\!\! \frac{2-\epsilon}{4} + \epsilon x(6;1)\label{eq:saf_work8}
\end{eqnarray}
The previous threshold on $\alpha_2/\alpha_1$ for $x(6;0)$ to be
zero, i.e., $(1-\epsilon)(3-2\epsilon)$, is valid for the case where
$x(1;0)=0$. Other decisions staying the same, when $x(1;0)$ is
positive and $x(1;1)=0$, $r_2$ increases and $r_1$ decreases.
Therefore the threshold on $\alpha_2/\alpha_1$ for $x(6;0)$ to be
zero changes, in particular, a simple derivation shows that it becomes $\alpha_2/\alpha_1
> (1-\epsilon)^2/\epsilon$. This gives the following two regions:
\\
\\
$ \frac{1+\epsilon-\epsilon^2}{1-\epsilon} \le
\frac{\alpha_2}{\alpha_1}
\le \frac{(1-\epsilon)^2}{\epsilon}$:\\
\\
The optimal policy is
\begin{equation*}
\begin{array}{llll}
\textrm{@queue 1}: (1,1,1):\textrm{switch}, &(1,1,0):\textrm{stay}, &(1,0,1):\textrm{switch}, &(1,0,0):\textrm{switch}. \\
\textrm{@queue 2}: (2,1,1):\textrm{stay}, &(2,1,0):\textrm{switch},
&(2,0,1):\textrm{stay}, &(2,0,0):\textrm{stay}.
\end{array}
\end{equation*}
From (\ref{eq:saf_work7}) and (\ref{eq:saf_work8}) it is easy to see
that this policy achieves
\begin{equation*}
r_1= \frac{(1-\epsilon)^2}{4},\;\;\;\;\; r_2 =\frac{2-\epsilon}{4}.
\end{equation*}
\\
\\
$\frac{\alpha_2}{\alpha_1}
> \frac{(1-\epsilon)^2}{\epsilon}$:\\
\\
The optimal decisions at queue 1 are to switch to queue 2. This policy achieves
\begin{equation*}
\begin{array}{llll}
\textrm{@queue 2}: (2,1,1):\textrm{stay}, &(2,1,0):\textrm{stay},
&(2,0,1):\textrm{stay}, &(2,0,0):\textrm{stay}.
\end{array}
\end{equation*}
This policy achives
\begin{equation*}
r_1= 0,\;\;\;\;\; r_2 =0.5.
\end{equation*}
Similar to Case-1, the case $\alpha_2/\alpha_1 < 1$ is symmetric and
can be solved similarly.

Thus we have characterized the corner point of the stability region for the two regions of $\epsilon$.
Using these corner points, it is easy to derive the expressions for the lines connecting these corner points,
which are given in Theorem \ref{thm:stab}.

\section*{Appendix C - Generalization to Non-symmetric Gilbert-Elliot Channels}
In the following, we state results analogous to the results established in Section \ref{Sec:2queues} for symmetric Gilbert-Elliot channels to the case of non-symmetric Gilbert-Elliot channel model as given in Fig. \ref{Fig:channel_MC}.
\begin{theorem}\label{thm:stab_nonsymGE}
\emph{When the connectivity processes follow the non-symmetric Gilbert-Elliot channel model, the rate region $\mathbf{\Lambda}_s$ is the set of all rates $r_1\ge 0$, $r_2 \ge 0$ that for $p_{01} < 
\frac{(1-p_{10})^2}{2-p_{10}}$ satisfy}
\beqn
p_{01}r_1 + h_1r_2 &\le& h_1\frac{p_{01}}{p_{01}+p_{10}}\nonumber\\
(1-p_{10})r_1 + h_2r_2  &\le& 1- \frac{p_{10}^2}{(p_{10}+p_{01})^2} - \frac{p_{10}p_{01}}{p_{01}+p_{10}}\nonumber\\
r_1 + r_2  &\le& 1- \frac{p_{10}^2}{(p_{10}+p_{01})^2} - \frac{p_{10}p_{01}}{p_{01}+p_{10}}\nonumber\\
h_2r_1 + (1-p_{10})r_2  &\le& 1- \frac{p_{10}^2}{(p_{10}+p_{01})^2} - \frac{p_{10}p_{01}}{p_{01}+p_{10}}\nonumber\\
h_1r_1 + p_{01}r_2 &\le&
h_1\frac{p_{01}}{p_{01}+p_{10}},\nonumber 
\eeqn
\emph{where $h_1 = (1-p_{01})(1-p_{10})$, $h_2=1+p_{10}-p_{10}^2$, and for $p_{01} \ge \frac{(1-p_{10})^2}{2-p_{10}} $ satisfy}
\beqn
p_{01}r_1 + h_3r_2 &\le& h_3\frac{p_{01}}{p_{10}+p_{01}}\nonumber\\
r_1 +r_2  &\le& 1- \frac{p_{10}^2}{(p_{10}+p_{01})^2} - \frac{p_{10}p_{01}}{p_{01}+p_{10}}\nonumber\\
h_3r_1 + p_{01}r_2  &\le&
h_3\frac{p_{01}}{p_{10}+p_{01}},\nonumber
 \eeqn
 \emph{where $h_3=(1-p_{10})(p_{10}+ (p_{10}+p_{01})(1-p_{10}))$}.
\end{theorem}
\begin{proof}
We enumerate the states as follows:
\begin{equation}\label{eq:state_enum}
\begin{array}{llll}
s=(1,1,1)\equiv1,\;\;\; &s=(1,1,0)\equiv2,\;\;\;
&s=(1,0,1)\equiv3,\;\;\;  & s=(1,0,0)\equiv4,\;\;\;\\
s=(2,1,1) \equiv5 & s=(2,1,0)\equiv6 & s=(2,0,1) \equiv7& s=(2,0,0) \equiv8.\\
\end{array}
\end{equation}
We rewrite the balance equations in (\ref{eq:LP}) in more details.
\beqn x(1;1)+x(1;0) = (1-p_{10})^2 \big(x(1;1)+x(5;0)\big)
\!\!\!\!\!&+&
\!\!\!\!\!p_{01}(1-p_{10})\big(x(2;1)+ x(6;0)\big)\nonumber\\
+ p_{01}(1-p_{10})\big(x(3;1)+x(7;0)\big) \!\!\!\!\!&+&\!\!\!\!\! p_{01}^2\big(x(4;1)+x(8;0)\big)\label{eq:safns1}  \\
x(2;1)+x(2;0) = p_{10}(1-p_{10}) \big(x(1;1)+x(5;0)\big)
\!\!\!\!\!&+&\!\!\!\!\!
(1-p_{10})(1-p_{01})\big(x(2;1) + x(6;0)\big) \nonumber\\
\!\!\!\!\!+ p_{01}p_{10}\big(x(3;1)+x(7;0)\big)
\!\!\!\!\!&+&\!\!\!\!\!
p_{01}(1-p_{01})\big(x(4;1)+x(8;0)\big)\label{eq:safns2}\\
\!\!\!\!\!&\ldots& \nonumber\\
x(5;1)+x(5;0) = (1-p_{10})^2 \big(x(5;1)+x(1;0)\big)
\!\!\!\!\!\!\!\!&+\!\!\!\!\!&
p_{01}(1-p_{10})\big(x(6;1)+ x(2;0)\big)\nonumber\\
+ p_{01}(1-p_{10})\big(x(7;1)+x(3;0)\big) \!\!\!\!\!&+&\!\!\!\!\! p_{01}^2\big(x(8;1)+x(4;0)\big) \label{eq:safns3} \\
x(7;1)+x(7;0)= p_{10}(1-p_{10}) \big(x(5;1)+x(1;0)\big)
\!\!\!\!\!&+&\!\!\!\!\!
p_{10}p_{01}\big(x(6;1) + x(2;0)\big) \nonumber\\
+ (1-p_{01})(1-p_{10}) \big(x(7;1)+x(3;0)\big) \!\!\!\!\!&+&\!\!\!\!\!
p_{01}(1-p_{01})\big(x(8;1)+x(4;0)\big)\label{eq:safns4}\\
&\ldots\nonumber  \eeqn\\
The following equations hold the channel state pairs
$(C_1,C_2)$. \beqn
x(1;1)+x(1;0)+ x(5;1)+x(5;0) = \frac{p_{01}^2}{(p_{01}+p_{10})^2} \label{eq:state_1ns}\\
x(2;1)+x(2;0)+ x(6;1)+x(6;0) = \frac{p_{01}p_{10}}{(p_{01}+p_{10})^2} \label{eq:state_2ns}\\
x(3;1)+x(3;0)+ x(7;1)+x(7;0) = \frac{p_{01}p_{10}}{(p_{01}+p_{10})^2} \label{eq:state_3ns}\\
x(4;1)+x(4;0)+ x(8;1)+x(8;0) = \frac{p_{10}^2}{(p_{01}+p_{10})^2}. \label{eq:state_4ns}
 \eeqn\\
 \\
Let $u_1=(x(1;1)+x(2;1)$ and $u_2=(x(5;1)+x(7;1))$. Summing up
(\ref{eq:safns1}) with (\ref{eq:safns2}) and (\ref{eq:safns3}) with
(\ref{eq:safns4}) we have
\begin{eqnarray*}
p_{10} u_1 \!\!\!\!&=&\!\!\!\! - \big(x(1;0)+x(2;0)\big) +
p_{01}\big( x(3;1)+x(4;1) \big) + p_{01}\big(x(7;0)+ x(8;0)
\big) + (1-p_{10}) \big( x(5;0)+x(6;0) \big)\\
p_{10} u_2 \!\!\!\!&=&\!\!\!\! - \big(x(5;0)+x(7;0)\big) +
p_{01}\big( x(6;1)+x(8;1) \big) + p_{01}\big(x(2;0)+ x(4;0)
\big) + (1-p_{10}) \big( x(1;0)+x(3;0) \big)\\
\end{eqnarray*}
\\
Rearranging and using (\ref{eq:state_1ns})-(\ref{eq:state_4ns}) we have
\begin{eqnarray}
u_1 \!\!\!\!\!\!&=&\!\!\!\!\! \frac{p_{01}(1\!-\!p_{10})}{p_{01}\!+\!p_{10}}\!+\!p_{01}\big(
x(3;1)+x(4;1)+x(7;0)+ x(8;0)\big) - (2-p_{10}) \big( x(1;0)+x(2;0) \big) - (1-p_{10}) \big( x(5;1)+x(6;1)\big)\nonumber\\
\!\!\!\!\!\!\label{eq:safns_work1}\\
u_2 \!\!\!\!\!\!&=&\!\!\!\!\! \frac{p_{01}(p_{01}\!+\!p_{10} \!-\! p_{01}p_{10})}{(p_{01}\!+\!p_{10})^2}\!+\!p_{01}\big(
x(2;0)\!-\!x(4;1) \!+ \!x(6;1)\!-\!x(8;0) \big) \!- \!(2\!-\!p_{10}) \big(x(5;0)\!+\!x(7;0) \big) \!- \!(1\!-\!p_{10}) \big(
x(1;1)\!+\!x(3;1)\big).\nonumber\\
\label{eq:safns_work2}
\end{eqnarray}
Using (\ref{eq:safns3}) in (\ref{eq:safns_work1})  and (\ref{eq:safns1})
in (\ref{eq:safns_work2}) we have
\begin{eqnarray}
u_1 \!\!\!\!\!&=&\!\!\!\!\! \frac{p_{01}(1\!-\!p_{10})}{p_{01}\!+\!p_{10}}\!+\!p_{01}\big(
x(3;1)\!+\!x(4;1)\!+\!x(7;0)\!+ \!x(8;0)\big) \!-\!
\frac{p_{01}^2(1\!-\!p_{10})}{p_{10}(2\!-\!p_{10})} \big( x(4;0)\!+\!x(8;1) \big)
\!-\!(1\!-\!p_{10})(1\!+\!\frac{p_{01}(1-p_{10})}{p_{10}(2-p_{10})} ) x(6;1)\nonumber \\
&+& \frac{1-p_{10}}{p_{10}(2-p_{10})} x(5;0)
-\frac{1+p_{10}-p_{10}^2}{p_{10}(2-p_{10})} x(1;0) -
\frac{p_{01}(1-p_{10})^2}{p_{10}(2-p_{10})} \big( x(3;0)+x(7;1) \big)
\!-\!\Big( 2\!-\!p_{10} \!+\! \frac{p_{01}(1-p_{10})^2}{p_{10}(2-p_{10})}\Big) x(2;0)\label{eq:safns_work3} \\
~\nonumber\\
u_2 \!\!\!\!\!\!&=&\!\!\!\!\! \frac{p_{01}(p_{01}\!+\!p_{10} \!-\! p_{01}p_{10})}{(p_{01}\!+\!p_{10})^2}
\!+\!p_{10}\big( x(2;0)\!+\!
x(6;1)\big)\! - \!\Big(
p_{01}\!+\!\frac{p_{01}^2(1\!-\!p_{10})}{p_{10}(2\!-\!p_{10})}\Big) \big(
x(4;1)\!+\!x(8;0) \big)
\!-\!(1\!-\!p_{10})(1\!+\!\frac{p_{01}(1\!-\!p_{10})}{p_{10}(2\!-\!p_{10})}) x(3;1)\nonumber \\
&+& \frac{1-p_{10}}{p_{10}(2-p_{10})} x(1;0)
-\frac{1+p_{10}-p_{10}^2}{p_{10}(2-p_{10})} x(5;0) -
\frac{p_{01}(1-p_{10})^2}{p_{10}(2-p_{10})} \big( x(2;1)+x(6;0) \big) \!-\!\Big( 2\!-\!p_{10} \!+\! \frac{p_{01}(1-p_{10})^2}{p_{10}(2-p_{10})}\Big)
x(7;0).\label{eq:safns_work4}
\end{eqnarray}
\\
\\
Using (\ref{eq:state_3ns}) and (\ref{eq:state_4ns}) in
(\ref{eq:safns_work3}) and (\ref{eq:state_2ns}) in (\ref{eq:safns_work4})
we have
\begin{eqnarray}
u_1 \!\!\!\!\!&=&\!\!\!\!\! \frac{p_{01}(1\!-\!p_{10})}{p_{01}\!+\!p_{10}}\!
-
\frac{p_{01}^2(1\!-\!p_{10})}{(2-p_{10})(p_{01}\!+\!p_{10})^2}\!
+
\frac{p_{01}^2(1\!-\!p_{10})}{p_{10}(2\!-\!p_{10})} \big(
x(4;1)+x(8;0)\big) +  \frac{p_{01}}{p_{10}(2-p_{10})} \big(x(3;1)+x(7;0)\big)\nonumber\\
\!\!\!&-&\!\!\!\!\!(1\!-\!p_{10})(1\!+\!\frac{p_{01}(1-p_{10})}{p_{10}(2-p_{10})} ) x(6;1)  \!+\!
\frac{1\!-\!p_{10}}{p_{10}(2\!-\!p_{10})} x(5;0)
\!-\!\frac{1\!+\!p_{10}\!-\!p_{10}^2}{p_{10}(2\!-\!p_{10})} x(1;0)
\!-\!\Big( 2\!\!-\!p_{10} \!+\! \frac{p_{01}(1\!-\!p_{10})^2}{p_{10}(2\!-\!p_{10})}\Big)  x(2;0)\label{eq:safns_work5} \\
u_2 \!\!\!\!\!\!&=&\!\!\!\!\! \frac{p_{01}(p_{01}\!+\!p_{10}
\!-\! p_{01}p_{10})}{(p_{01}\!+\!p_{10})^2}
-
\frac{p_{01}^2(1\!-\!p_{10})^2}{(2-p_{10})(p_{01}\!+\!p_{10})^2}\!
- \!\Big(
p_{01}\!+\!\frac{p_{01}^2(1\!-\!p_{10})}{p_{10}(2\!-\!p_{10})}\Big) \big(
x(4;1)\!+\!x(8;0) \big)
+
\frac{p_{01}}{p_{10}(2-p_{10})} \big(x(2;0)+x(6;1)\big)\nonumber\\
&-&
(1\!-\!p_{10})(1\!+\!\frac{p_{01}(1\!-\!p_{10})}{p_{10}(2\!-\!p_{10})}) x(3;1)
+
\frac{1-p_{10}}{p_{10}(2-p_{10})} x(1;0)
-
\frac{1+p_{10}-p_{10}^2}{p_{10}(2-p_{10})} x(5;0)
\!-\!
\Big( 2\!-\!p_{10} \!+\! \frac{p_{01}(1-p_{10})^2}{p_{10}(2-p_{10})}\Big)
x(7;0).\label{eq:safns_work6}
\end{eqnarray}
\\
Consider the LP objective function $\alpha_1 \big(x(1;1)+x(2;1)\big)
+ \alpha_2\big(x(5;1)+x(7;1)\big)$, and note that the solution to
this LP is a stationary deterministic policy for any given
$\alpha_1$ and $\alpha_2$. This means that, for any state $s$ either
$x(s;1)$ or $x(s;0)$ has to be zero. In order to maximize $\alpha_1
\big(x(1;1)+x(2;1)\big) + \alpha_2\big(x(5;1)+x(7;1)\big)$ we need

\begin{eqnarray*}
x(7;0) \!\!\!\!&=&\!\!\!\! 0 \;\;\;\;\textrm{if  }\;\; \frac{\alpha_2}{\alpha_1} \ge \frac{p_{01}}{p_{10}(2-p_{10})^2+p_{01}(1-p_{10})^2},\\
x(3;1) \!\!\!\!&=&\!\!\!\! 0 \;\;\;\;\textrm{if  }\;\;\frac{\alpha_2}{\alpha_1} \ge \frac{p_{01}}{(1-p_{10})(p_{10}(2-p_{10})+p_{01}(1-p_{10}))},\\
x(5;0) \!\!\!\!&=&\!\!\!\! 0 \;\;\;\;\textrm{if  }\;\;\frac{\alpha_2}{\alpha_1} \ge \frac{1-p_{10}}{1+p_{10}-p_{10}^2},\\
x(8;0)=x(4;1)\!\!\!\!&=&\!\!\!\! 0 \;\;\;\;\textrm{if  }\;\;\frac{\alpha_2}{\alpha_1} \ge 1,\\
x(6;0) \!\!\!\!&=&\!\!\!\! 0 \;\;\;\;\textrm{if  }\;\;\frac{\alpha_2}{\alpha_1} \ge \frac{(1-p_{10})(p_{10}(2-p_{10})+p_{01}(1-p_{10}))}{p_{01}},\\
x(1;1) \!\!\!\!&=&\!\!\!\! 0
\;\;\;\;\textrm{if}\;\;\frac{\alpha_2}{\alpha_1} \ge
\frac{1+p_{10}-p_{10}^2}{1-p_{10}}.
\end{eqnarray*}
\\
Note that we have that
\begin{eqnarray*}
\frac{p_{10}(2-p_{10})^2+p_{01}(1-p_{10})^2}{p_{01}} &\ge& \frac{1+p_{10}-p_{10}^2}{1-p_{10}} \ge 1\\
\frac{p_{10}(2-p_{10})^2+p_{01}(1-p_{10})^2}{p_{01}} &\ge&  \frac{(1-p_{10})(p_{10}(2-p_{10})+p_{01}(1-p_{10}))}{p_{01}} \ge 1,\\
\end{eqnarray*}
whenever $p_{10}+p_{10} <1$ (the condition for positive correlation). Consider the following two
cases:
\\
\\
\textbf{Case-1: $p_{01} \ge \frac{(1-p_{10})^2}{2-p_{10}}$}
\\
In this case we have $\frac{(1-p_{10})(p_{10}(2-p_{10})+p_{01}(1-p_{10}))}{p_{01}} \le \frac{1+p_{10}-p_{10}^2}{1-p_{10}}$. This means that we have the
following optimal
policies depending on the value of $\alpha_2/\alpha_1$.\\
\\
$1 \le \frac{\alpha_2}{\alpha_1} \le \frac{(1-p_{10})(p_{10}(2-p_{10})+p_{01}(1-p_{10}))}{p_{01}}$:
\begin{equation*}
\begin{array}{llll}
\textrm{@queue 1}: (1,1,1):\textrm{stay}, &(1,1,0):\textrm{stay}, &(1,0,1):\textrm{switch}, &(1,0,0):\textrm{switch}. \\
\textrm{@queue 2}: (2,1,1):\textrm{stay}, &(2,1,0):\textrm{switch},
&(2,0,1):\textrm{stay}, &(2,0,0):\textrm{stay}.
\end{array}
\end{equation*}
Substituting the above zero variables into (\ref{eq:safns_work5}) and
(\ref{eq:safns_work6}), it can be seen that this policy achieves the
rate pair
\begin{equation*}
r_1= \frac{1-p_{10}}{2-p_{10}}\frac{p_{01}p_{10}+p_{01}(1-p_{10})(p_{01}+p_{10})}
{(p_{01}+p_{10})^2},\;\;\;\;\; r_2
=\frac{1}{2-p_{10}}\frac{p_{01}p_{10}+p_{01}(1-p_{10})(p_{01}+p_{10})}
{(p_{01}+p_{10})^2}.
\end{equation*}
\\
\\
$\frac{\alpha_2}{\alpha_1} > \frac{(1-p_{10})(p_{10}(2-p_{10})+p_{01}(1-p_{10}))}{p_{01}}$:
\begin{equation*}
\begin{array}{llll}
\textrm{@queue 2}: (2,1,1):\textrm{stay}, &(2,1,0):\textrm{stay},
&(2,0,1):\textrm{stay}, &(2,0,0):\textrm{stay}.
\end{array}
\end{equation*}
In this case it is optimal to stay at queue 2 for all channel
conditions. The decisions at queue 1 are to switch to queue 2.
Namely, it is sufficient that at least one state corresponding to
server being at queue 1 to take a switch decision, which is the case
for $\alpha_2/\alpha_1\ge \frac{(1-p_{10})(p_{10}(2-p_{10})+p_{01}(1-p_{10}))}{p_{01}}$, since $x(3;1)
=0$ if $\alpha_2/\alpha_1\ge (\frac{(1-p_{10})(p_{10}(2-p_{10})+p_{01}(1-p_{10}))}{p_{01}})^{-1}$. Since
the policy decides to always stay at queue 2, it achieves the rate pair
\begin{equation*}
r_1= 0,\;\;\;\;\; r_2 =\frac{p_{01}}{p_{01}+p_{10}}.
\end{equation*}
Note that the case for $\alpha_2/\alpha_1 < 1$ is symmetric and can
be obtained similarly.\\
\\
\\
\textbf{Case-2: $p_{01} < \frac{(1-p_{10})^2}{2-p_{10}}$}\\
In this case we have $\frac{(1-p_{10})(p_{10}(2-p_{10})+p_{01}(1-p_{10}))}{p_{01}} > \frac{1+p_{10}-p_{10}^2}{1-p_{10}}$. This means that before the
state $x(6;0)$ becomes zero, namely for
$\frac{1+p_{10}-p_{10}^2}{1-p_{10}} < \alpha_2/\alpha_1 < \frac{(1-p_{10})(p_{10}(2-p_{10})+p_{01}(1-p_{10}))}{p_{01}}$, having $x(1;1)=0$ is optimal. This means
that there is one more corner point to the rate region for $p_{01} < \frac{(1-p_{10})^2}{2-p_{10}}$. We have the following optimal
policies.
\\
\\
$1 \le \frac{\alpha_2}{\alpha_1} \le
\frac{1+p_{10}-p_{10}^2}{1-p_{10}}$:\\
\begin{equation*}
\begin{array}{llll}
\textrm{@queue 1}: (1,1,1):\textrm{stay}, &(1,1,0):\textrm{stay}, &(1,0,1):\textrm{switch}, &(1,0,0):\textrm{switch}. \\
\textrm{@queue 2}: (2,1,1):\textrm{stay}, &(2,1,0):\textrm{switch},
&(2,0,1):\textrm{stay}, &(2,0,0):\textrm{stay}.
\end{array}
\end{equation*}
This policy is the same policy as in the previous case and it
achieves the rate pair
\begin{equation*}
r_1= \frac{1-p_{10}}{2-p_{10}}\frac{p_{01}p_{10}+p_{01}(1-p_{10})(p_{01}+p_{10})}
{(p_{01}+p_{10})^2},\;\;\;\;\; r_2
=\frac{1}{2-p_{10}}\frac{p_{01}p_{10}+p_{01}(1-p_{10})(p_{01}+p_{10})}
{(p_{01}+p_{10})^2}.
\end{equation*}
\\
\\
$ \frac{\alpha_2}{\alpha_1}
> \frac{1+p_{10}-p_{10}^2}{1-p_{10}}$:\\
\\
We have the following deterministic actions.
\begin{equation*}
\begin{array}{llll}
\textrm{@queue 1}: (1,1,1):\textrm{switch}, &(1,1,0):\textrm{?}, &(1,0,1):\textrm{switch}, &(1,0,0):\textrm{switch}. \\
\textrm{@queue 2}: (2,1,1):\textrm{stay}, &(2,1,0):\textrm{?},
&(2,0,1):\textrm{stay}, &(2,0,0):\textrm{stay}.
\end{array}
\end{equation*}
In order to find the final threshold on $\alpha_2/\alpha_1$, we
substitute the above deterministic decisions in (\ref{eq:safns2}),
(\ref{eq:safns3}) and (\ref{eq:safns4}). Utilizing also
(\ref{eq:state_1ns}), (\ref{eq:state_2ns}), (\ref{eq:state_3ns}) and
(\ref{eq:state_4ns}) we obtain
\begin{eqnarray}
x(2;1) \!\!\!\!&=&\!\!\!\! \frac{(1-\epsilon)^2}{4} - (1-\epsilon)^2 x(6;1)\label{eq:safns_work7}\\
x(5;1)+x(7;1) \!\!\!\!&=&\!\!\!\! \frac{p_{01}(p_{10}-p_{10}p_{01}+p_{01})}{(p_{01}+p_{10})^2} + p_{01} x(6;1)\label{eq:safns_work8}
\end{eqnarray}
The previous threshold on $\alpha_2/\alpha_1$ for $x(6;0)$ to be
zero, i.e., $\frac{(1-p_{10})(p_{10}(2-p_{10})+p_{01}(1-p_{10}))}{p_{01}}$, is valid for the case where
$x(1;0)=0$. Other decisions staying the same, when $x(1;0)$ is
positive and $x(1;1)=0$, $r_2$ increases and $r_1$ decreases.
Therefore the threshold on $\alpha_2/\alpha_1$ for $x(6;0)$ to be
zero changes, in particular, a simple derivation shows that it becomes $\alpha_2/\alpha_1
> (1-p_{10})(1-p_{01})/p_{01}$. This gives the following two regions:
\\
\\
$ \frac{1+p_{10}-p_{10}^2}{1-p_{10}} \le
\frac{\alpha_2}{\alpha_1}
\le (1-p_{10})(1-p_{01})/p_{01}$:\\
\\
The optimal policy is
\begin{equation*}
\begin{array}{llll}
\textrm{@queue 1}: (1,1,1):\textrm{switch}, &(1,1,0):\textrm{stay}, &(1,0,1):\textrm{switch}, &(1,0,0):\textrm{switch}. \\
\textrm{@queue 2}: (2,1,1):\textrm{stay}, &(2,1,0):\textrm{switch},
&(2,0,1):\textrm{stay}, &(2,0,0):\textrm{stay}.
\end{array}
\end{equation*}
From (\ref{eq:safns_work7}) and (\ref{eq:safns_work8}) it is easy to see
that this policy achieves
\begin{equation*}
r_1= (1-p_{10})(1-p_{01})\frac{p_{01}p_{10}}{(p_{01}+p_{10})^2},\;\;\;\;\; r_2 =\frac{p_{01}(p_{10}-p_{10}p_{01}+p_{01})}{(p_{01}+p_{10})^2}.
\end{equation*}
\\
\\
$\frac{\alpha_2}{\alpha_1}
> (1-p_{10})(1-p_{01})/p_{01}$:\\
\\
The optimal policy is
\begin{equation*}
\begin{array}{llll}
\textrm{@queue 2}: (2,1,1):\textrm{stay}, &(2,1,0):\textrm{stay},
&(2,0,1):\textrm{stay}, &(2,0,0):\textrm{stay}.
\end{array}
\end{equation*}
The optimal decisions at queue 1 are to switch to queue 2. This policy achieves
\begin{equation*}
r_1= 0,\;\;\;\;\; r_2 = \frac{p_{01}}{p_{01}+p_{10}}.
\end{equation*}
Similar to Case-1, the case $\alpha_2/\alpha_1 < 1$ is symmetric and
can be solved similarly.

Thus we have characterized the corner point of the stability region for the two regions of $p_{01}$ and $p_{10}$.
Using these corner points, it is easy to derive the expressions for the lines connecting these corner points,
which are given in Theorem~\ref{thm:stab_nonsymGE}.
\end{proof}

Closely examining the upper bound on sum-rate $r_1 + r_2$, the term $1-p_{10}^2/(p_{10}+p_{01})^2$ is the steady state probability that at least one channel is in ON state. This is the maximum achievable sum-rate value for the system with zero switchover delay studied in \cite{tass93}. Therefore, the term $\frac{p_{10}p_{01}}{p_{01}+p_{10}}$ is exactly the loss due to switchover delay. It can be shown that, under a sum-rate-optimal policy, this term is equal to the steady state probability that server is at a queue with an OFF channel state when the other queue is at an ON channel state.

The FBDC policy is asymptotically throughput-optimal under the non-symmetric Gilbert-Elliot channel model. This is straightforward as the FBDC policy only needs to solve the LP in Algorithm~\ref{alg:FBDC} for a given Markovian state transition structure, and the non-symmetric Gilbert-Elliot channel model leads to a Markovian state transition structure.
%
For the non-symmetric Gilbert-Elliot channels case, the mappings from the queue sizes to the corner points of the rate region used by the FBDC policy, analogues to the mappings in tables \ref{Fig:optimal_table}, and \ref{Fig:optimal_table_eps_small} can be obtained from the slopes of the lines forming the boundary of the stability region. Furthermore, an analysis very similar to the one in Section~\ref{Sec:Myopic_Pol} gives the corresponding mapping for the One Lookahead Myopic (OLM) policy.
These mappings 
are shown in Table~\ref{Fig:optimal_table_nonsym_GE_small_eps} for the case of $p_{01} < \frac{(1-p_{10})^2}{2-p_{10}}$, and in Table~\ref{Fig:optimal_table_nonsym_GE} for the case of $p_{01} \ge \frac{(1-p_{10})^2}{2-p_{10}}$.

\begin{table}
\centering \psfrag{1\r}[l][][.8]{$0$}
\psfrag{2\r}[l][][.8]{\,$T^*_1(\epsilon)$}
\psfrag{n\r}[l][][.8]{$\!\!\!\!\!\!\!T_1(\epsilon)$}
\psfrag{3\r}[l][][.8]{$\,T^*_2(\epsilon)$}
\psfrag{q\r}[l][][.8]{\!\!\!\!\!$T_2(\epsilon)$}
\psfrag{4\r}[l][][.8]{\!$1$}
\psfrag{5\r}[l][][.8]{\!\!\!\!\!\!\!$T^*_3(\epsilon)$}
\psfrag{y\r}[l][][.8]{\!\!$T_3(\epsilon)$}
\psfrag{0\r}[l][][.8]{\!\!\!\!\!\!\!$T^*_4(\epsilon)$}
\psfrag{w\r}[l][][.8]{\!\!\!\!$T_4(\epsilon)$}
\psfrag{9\r}[l][][.8]{\!\!\!\!\!\!\!\!\!\!\!\!\!\!\!\!\!\!\textrm{corner}
$b_0$}
\psfrag{z\r}[l][][.8]{\!\!\!\!\!\!\!\!\!\!\!\!\!\!\!\!\!\textrm{corner}
$b_1$}
\psfrag{x\r}[l][][.8]{\!\!\!\!\!\!\!\!\!\!\!\!\textrm{corner}
$b_2$}
\psfrag{6\r}[l][][.8]{\!\!\!\!\!\!\!\!\!\!\!\!\!\!\textrm{corner}
$b_3$}
\psfrag{7\r}[l][][.8]{\!\!\!\!\!\!\!\!\!\!\!\!\!\!\!\!\!\textrm{corner}
$b_4$}
\psfrag{8\r}[l][][.8]{\!\!\!\!\!\!\!\!\!\!\!\!\!\!\!\!\!\textrm{corner}
$b_5$}
\psfrag{a\r}[l][][.75]{(1,1,1): \textrm{switch}}
\psfrag{b\r}[l][][.75]{(1,1,0): \textrm{switch}}
\psfrag{d\r}[l][][.75]{(1,0,1): \textrm{switch}}
\psfrag{e\r}[l][][.75]{(1,0,0): \textrm{switch}}
\psfrag{f\r}[l][][.75]{(2,1,1): \textrm{stay}}
\psfrag{g\r}[l][][.75]{(2,1,0): \textrm{stay}}
\psfrag{h\r}[l][][.75]{(2,0,1): \textrm{stay}}
\psfrag{i\r}[l][][.75]{(2,0,0): \textrm{stay}}
\psfrag{j\r}[l][][.75]{(1,1,1): \textrm{stay}}
\psfrag{k\r}[l][][.75]{(1,1,0): \textrm{stay}}
\psfrag{l\r}[l][][.75]{(2,1,0): \textrm{switch}}
\psfrag{m\r}[l][][.75]{(1,0,0): \textrm{stay}}
\psfrag{o\r}[l][][.75]{(2,0,0): \textrm{switch}}
\psfrag{p\r}[l][][.75]{(1,0,1): \textrm{stay}}
\psfrag{r\r}[l][][.75]{(2,1,1): \textrm{switch}}
\psfrag{s\r}[l][][.75]{(2,0,1): \textrm{switch}}
\psfrag{t\r}[l][][.75]{(2,0,0): \textrm{stay}}
\psfrag{u\r}[l][][.75]{(2,0,1): \textrm{stay}}
\psfrag{v\r}[l][][.75]{(2,0,0): \textrm{switch}}
\psfrag{c\r}[l][][1]{$\!\!\!\!\!\!\frac{Q_2}{Q_1}$}
\includegraphics[width=0.50\textwidth]{myopic_table_1_small_eps.eps}
\captionsetup{font=footnotesize} \caption{Mapping from
the queue sizes to the corners of $\mathbf{\Lambda_s}$,
$b_0,b_1,b_2,b_3,b_4,b_5$, for $p_{01} <
\frac{(1-p_{10})^2}{2-p_{10}}$. For each state $\mathbf{s}=(m(t),C_1(t),C_2(t))$ the
optimal action is specified. The thresholds
on $Q_2/Q_1$ for the FBDC policy are
$0,T^*_1=p_{01}/((1-p_{01})(1-p_{10})),T^*_2= (1-p_{10})/(1+p_{10}-p_{10}^2),1,T^*_3 =1/T^*_2,T^*_4=1/T^*_1$ and for the OLM policy
are $0,T_1=p_{01}/(1-p_{10}),T_2=(1-p_{10})/(2-p_{01}),1,T_3=1/T_2,T_4=1/T_1$.
For example corner
$b_1$ is chosen in the FBDC policy if $1\le Q_2/Q_1 < T^*_3$, whereas in
the OLM policy if $1 \le Q_2/Q_1 < T_3$.} \label{Fig:optimal_table_nonsym_GE_small_eps}
\end{table}

\begin{table}
\centering \psfrag{1\r}[l][][.8]{$0$}
\psfrag{2\r}[l][][.8]{$\!\!T_1$}
\psfrag{3\r}[l][][.8]{\!\!$T^*_1$}
\psfrag{4\r}[l][][.8]{\!$1$}
\psfrag{5\r}[l][][.8]{\!\!\!\!$T^*_2$}
\psfrag{0\r}[l][][.8]{\!\!\!\!$T_2$}
\psfrag{9\r}[l][][.8]{\!\!\!\!\!\!\!\!\!\!\!\!\!\!\!\!\!\!\textrm{corner}
$b_0$}
\psfrag{6\r}[l][][.8]{\!\!\!\!\!\!\!\!\!\!\!\!\!\!\!\!\!\textrm{corner}
$b_1$} \psfrag{7\r}[l][][.8]{\!\!\!\!\!\!\!\!\!\!\!\!\textrm{corner}
$b_2$}
\psfrag{8\r}[l][][.8]{\!\!\!\!\!\!\!\!\!\!\!\!\!\!\textrm{corner}
$b_3$}
\psfrag{a\r}[l][][.8]{\!\!(1,1,1): \textrm{switch}}
\psfrag{b\r}[l][][.8]{\!\!(1,1,0): \textrm{switch}}
\psfrag{d\r}[l][][.8]{\!\!(1,0,1): \textrm{switch}}
\psfrag{e\r}[l][][.8]{\!\!(1,0,0): \textrm{switch}}
\psfrag{f\r}[l][][.8]{\!\!(2,1,1): \textrm{stay}}
\psfrag{g\r}[l][][.8]{\!\!(2,1,0): \textrm{stay}}
\psfrag{h\r}[l][][.8]{\!\!(2,0,1): \textrm{stay}}
\psfrag{i\r}[l][][.8]{\!\!(2,0,0): \textrm{stay}}
\psfrag{j\r}[l][][.8]{\!\!(1,1,1): \textrm{stay}}
\psfrag{k\r}[l][][.8]{\!\!(1,1,0): \textrm{stay}}
\psfrag{l\r}[l][][.8]{\!\!(2,1,0): \textrm{switch}}
\psfrag{m\r}[l][][.8]{\!\!(1,0,0): \textrm{stay}}
\psfrag{o\r}[l][][.8]{\!\!(2,0,0): \textrm{switch}}
\psfrag{p\r}[l][][.8]{\!\!(1,0,1): \textrm{stay}}
\psfrag{r\r}[l][][.8]{\!\!(2,1,1): \textrm{switch}}
\psfrag{s\r}[l][][.8]{\!\!(2,0,1): \textrm{switch}}
\psfrag{t\r}[l][][.8]{\!\!(2,0,0): \textrm{stay}}
\psfrag{u\r}[l][][.8]{\!\!(2,0,1): \textrm{stay}}
\psfrag{v\r}[l][][.8]{\!\!(2,0,0): \textrm{switch}}
\psfrag{w\r}[l][][.8]{\!\!(1,0,0): \textrm{stay}}
\psfrag{c\r}[l][][1]{$\!\!\!\!\!\!\frac{Q_2}{Q_1}$}
\includegraphics[width=0.50\textwidth]{myopic_optimal_table_1.eps}
\captionsetup{font=footnotesize} \caption{Mapping from
the queue sizes to the corners of $\mathbf{\Lambda_s}$,
$b_0,b_1,b_2,b_3$, for $p_{01} \ge
\frac{(1-p_{10})^2}{2-p_{10}}$. For each state $\mathbf{s}=(m(t),C_1(t),C_2(t))$ the
optimal action is specified. The thresholds
on $Q_2/Q_1$ for the FBDC policy are
$0,T^*_1=p_{01}/((1-p_{10})(p_{10}+ (p_{10}+p_{01})(1-p_{10}))),1,T^*_2=(1-p_{10})(p_{10}+ (p_{10}+p_{01})(1-p_{10}))/p_{01}$ and for the OLM policy
are $0,T_1=p_{01}/(1-p_{10}),1,T_2=(1-p_{10})/p_{01}$. For example corner
$b_1$ is chosen in the FBDC policy if $1\le Q_2/Q_1 < T^*_2$, whereas in
the OLM policy if $1 \le Q_2/Q_1 < T_2$.} \label{Fig:optimal_table_nonsym_GE}
\end{table}

\section*{Appendix D - Proof of Theorem \ref{thm:FBDC}}
Let $t_k$ be the first slot of the $k$th frame where $t_{k+1}=t_k+T$.
Let $D_i(t)$ 
be the \emph{service opportunity} given to queue $i$ at time slot $t$, where 
$D_i(t)$ is equal to $1$ if queue $i$ is scheduled at time
slot $t$ (regardless of whether queue $i$ is empty or not) and zero otherwise. We have the following queue evolution
relation:
\beq
Q_i(t+1) = \max(Q_i(t) - D_i(t),0) + A_i(t).
\eeq
Similarly, the following $T$-step queue evolution relation holds:
\begin{eqnarray}
Q_i(t_k+T) &\!\!\!\!\!\!\!\!\!\!\!\!\!\!\!\!\!\!\!\!\leq&\!\!\!\!\!\!\!\!\!\!\!\!\!\!\!\!\!\!\!\! \max\!\!  \left\{  Q_i(t_k)  \!-  \!\!\!   \sum_{\tau=0}^{T-1} \!\!  D_i(t_k+\tau), 0   \!\!\right\} \nonumber\\
&\;\;\;\;\;\;+&    \sum_{\tau=0}^{T-1} A_i(t_k+\tau),\label{eq:FBDC_drift0}
\end{eqnarray}
where $\sum_{\tau=0}^{T-1} D_i(t_k+\tau)$ is the total \emph{service opportunity} given to queue $i$ during the $k^{\textrm{th}}$ frame.
To see this, note that if $\sum_{\tau=0}^{T-1} D_i(t_k+\tau)$, the total \emph{service opportunity} given to queue $i$ during the 
$k^{\textrm{th}}$ frame, is smaller than $Q_i(t_k)$, then 
we have an equality. Otherwise, the first term 
is 0 and
we have an inequality. This is because some of the arrivals during the frame might depart before the end of the frame.
Note that $\sum_{\tau=0}^{T-1} D_i(t_k+\tau)$ denotes the link $i$ departures that would happen in the corresponding \emph{saturated system} if we were to apply the \emph{same} switching decisions over $T$ time slots in the corresponding saturated system.
We first prove stability at the frame boundaries.
Define the quadratic Lyapunov function
\begin{equation*}
L(\mathbf{Q}(t)) =\sum_{i=1}^2 Q_i^2(t),
\end{equation*}
which represents a quadratic measure of the total load in the system at time slot $t$.
Define the $T$-step conditional
drift
\begin{equation*}
\Delta_{T}(t_k)\triangleq \mathbb{E} \left[ L(\mathbf{Q}(t_k+T)) -
L(\mathbf{Q}(t_k))\big|\mathbf{Q}(t_k)\right],
\end{equation*}
where the conditional expectation is over the randomness in arrivals and possibly the scheduling decisions. 
%
%
%
Squaring both sides of (\ref{eq:FBDC_drift0}), using $\max(0,x)^2 \leq x^2, \forall x \in \mathbb{N}\cup \{0\}$, and $D_i(t) \le 1,\forall t$ we have
\begin{eqnarray} 
\displaystyle
Q_i(t_k+T)^2 \!\!\!\!\!\!&-& \!\!\!\!\!\!Q_i(t_k)^2 \le T^2  +
\Big(\sum_{\tau=0}^{T-1}A_i(t_k+\tau) \Big)^2\nonumber\\
&\!\!\!\!\!\!\!\!\!\!\!\!\!\!\!\!\!\!\!\!\!\!\!\!\!\!\!\!\!\!\!\!\!\!\!\!\!\!\!\!\!\!\!\!\!\!\!\!\!\!\!\!\!\!\!\!\!\!\!\!\!\!\!\!\!\!-& \!\!\!\!\!\!\!\!\!\!\!\!\!\!\!\!\!\!\!\!\!\!\!\!\!\!\!\!\!\!\!\!\!\!\!\!\!\!2Q_i(t_k)
\Big(\!\!\sum_{\tau=0}^{T-1}\!\!\!D_i(t_k+\tau)
-\!\!\!\sum_{\tau=0}^{T-1}\!\!\!A_i(t_k+\tau)\!\Big).\label{eq:drift_n2}
\end{eqnarray}
Summing (\ref{eq:drift_n2}) over the queues, using 
$\mathbb{E}[A_i(t)^2]
\le A_{\max}^2$ and $\mathbb{E}[A_i(t_1)A_i(t_2)] \le
\sqrt{\mathbb{E}[A_i(t_1)]^2\mathbb{E}[A_i(t_2)]^2}\le
A_{\max}^2$ for all time slots $t_1$ and $t_2$, we can easily derive the
following $T$-step conditional Lyapunov drift 
\begin{eqnarray}
\Delta_{T}(t_k) &\!\!\!\!\!\le& \!\!\!\!\!2BT^2 + 2T \sum_{i} 
Q_i(t_k)\lambda_i \nonumber\\
&\!\!\!\!\!\!\!\!\!\!\!\!\!\!\!\!\!\!\!\!\!\!\!\!\!\!\!\!\!\!\!\!\!\!\!\!\!\!\!\!\!\!-&\!\!\!\!\!\!\!\!\!\!\!\!\!\!\!\!\!\!\!\!\!\!\!\!2\sum_{i}
Q_i(t_k) \mathbb{E}\left[ \sum_{\tau=0}^{T-1} D_i(t_k+\tau)
\big|\mathbf{Q}(t_k)\! \right]\!\!, \label{eq:drift_intermediate}
\end{eqnarray}
where
$B\doteq 1+A_{\max}^2$ 
and we used the fact that the arrival processes are i.i.d. over time, independent of the queue lengths.
%
%
%
Recall the definition of the reward functions $\overline{r}_i(\mathbf{s}_t,\mathbf{a}_t)$ 
in in (\ref{eq:dept_rate1}) and
(\ref{eq:dept_rate})
and let
$\overline{r}_i(\mathbf{s}_t,\mathbf{a}_t)$ be the reward function associated with applying
policy $\boldsymbol{\pi}^*$ given in the definition of the FBDC policy in Algorithm~\ref{alg:FBDC} to the saturated system.
Let $\overline{r}_i(t)$ denote $\overline{r}_i(\mathbf{s}_t,\mathbf{a}_t)$ for notational simplicity, $i \in \{1,2\}$.
Note that $r_i(t)$ is equal to $D_i(t)$, since $D_i(t)$ is the \emph{service opportunity} given to link $i$ at time slot $t$. Now let $\mathbf{r}^* = (r_i^*)_i$ be the infinite horizon
average rate associated with policy $\boldsymbol{\pi}^*$. Let $\mathbf{x}^*$ be
the optimal vector of state-action frequencies corresponding to
$\boldsymbol{\pi}^*$. Define the time-average empirical reward from queue $i$ in
the saturated system, $\hat{r}_{T,i}(t_k)$, 
$i \in \{1,2\}$ by\\
\beq \hat{r}_{ T,i }(t_k) \doteq \frac{1}{T}\sum_{\tau=0}^{T-1}
\overline{r}_i(t_k+\tau)\nonumber.  
\eeq
Similarly, define the time average empirical state-action frequency
vector $\hat{\mathbf{x}}_{ T}(t_k;\mathbf{s},\mathbf{a})$.
\begin{equation}\label{eq:empr_saf}
\displaystyle \hat{\mathbf{x}}_T(t_k;\mathbf{s},\mathbf{a}) \doteq \frac{1}{T}\sum_{\tau=t_k}^{t_k+T-1} I_{\{\mathbf{s}_{\tau}=\mathbf{s},\mathbf{a}_{\tau}=\mathbf{a}\}},
\end{equation}
where $I_E$ is the indicator function of an event $E$, i.e., $I_E=1$ if $E$ occurs and $I_E=0$ otherwise.
%
Using the definition of the reward functions 
in (\ref{eq:dept_rate1}) and (\ref{eq:dept_rate}), we have that 
\begin{equation*}
\hat{r}_{ T,i }(t_k) =  \sum_{\mathbf{s}\in \mathcal{S}}\sum_{\mathbf{a}}\overline{r}_i(\mathbf{s},\mathbf{a}) \hat{\mathbf{x}}_{ T}(t_k;\mathbf{s},\mathbf{a}),\, i \in \{1,2\},
\end{equation*}
and $\hat{\mathbf{r}}_{ T}(t_k) = (\hat{r}_{ T,1 }(t_k))_{i}$. Similarly, we have
\begin{equation*}
r_i^* =  \sum_{\mathbf{s}\in \mathcal{S}}\sum_{\mathbf{a}}\overline{r}_i(\mathbf{s},\mathbf{a}) \mathbf{x}^*(\mathbf{s},\mathbf{a}),\, i \in \{1,2\}.
\end{equation*}
%
%
Now we utilize the following key MDP theory
result in Lemma 4.1 \cite{shie05}, which states that as $T$ increases, $\hat{\mathbf{x}}_{ T}(t_k) = (\hat{\mathbf{x}}_{ T}(t_k;\mathbf{s},\mathbf{a}))_{\mathbf{s},\mathbf{a}}$
converges to $\mathbf{x}^*$. 
\begin{lemma}\label{lem:MDP_convergence}
\emph{
For every choice of initial state distribution, there exists
constants $c_1$ and $c_2$ such that}
\begin{equation*}
\mathbb{P}(||\hat{\mathbf{x}}_{ T}(t_k)-\mathbf{x}^*|| \ge \delta_0)
\le c_1 e^{-c_2\delta_0^2 T}, \;\;\; \forall T\ge1,\; \forall
\delta_0 >0.
\end{equation*}
\emph{Furthermore, convergence of $\hat{\mathbf{x}}_{ T}(t_k)$ to $\mathbf{x}^*$ is with probability (w.p.) 1.}
\end{lemma}
This result applies in our system because every extreme point $\mathbf{x}^*$ of $\mathbf{X}$
\emph{can be} attained by a stationary and deterministic policy that has a
single irreducible recurrent class in its underlying Markov chain \cite{shie05}, \cite{puterman05}\footnote{Note that in general multiple stationary-deterministic policies can yield the same optimal reward vector $\mathbf{r}^*$. Among these, we choose the one that forms a Markov chain with a single recurrent class.}.
Due to the linear mapping from the state-action frequencies to the rewards, by Schwartz inequality, each component of
$\mathbf{\hat{r}}_T(t_k)$ also converges to the corresponding component of $\mathbf{r}^*$. Therefore, we have that
for every choice of initial state
distribution, there exists constants $c_1$ and $c_2$ such that
\begin{equation}\label{eq:MDP_rate_conv}
\mathbb{P}(||\hat{\mathbf{r}}_{ T}(t_k)-\mathbf{r}^*|| \ge \delta_1)
\le c_1 e^{-c_2\delta_1^2 T}, \; \forall T\ge1, \forall \delta_1
>0.
\end{equation}
Furthermore, convergence of $\hat{\mathbf{r}}_{ T}(t_k)$ to $\mathbf{r}^*$ is w.p. 1.
%
%
%
%
%
%
%
%
%
Now let 
$R_T(t_k) \doteq\sum_{i} Q_i(t_k) \hat{r}_{ T,i }(t_k)$ and
$R^*(t_k)\doteq \sum_{i} Q_i(t_k) r_i^*$.  We rewrite the drift
expression:
\begin{eqnarray}
\frac{\Delta_T(t_k)}{2T} \!\!\!\!\!\!\!\!\!\!\!\!\!\!\!\!\!\!\!\!\!\!&\leq&\!\!\!\!\!\!\!\!\!\!\!\!\!\!\!\!\!\! \!\!\!BT \!+
\sum_{i}Q_i(t_k) \lambda_i -
\mathbb{ E}\left[ R_T(t_k) \big|\mathbf{Q}(t_k) \right] \nonumber\\
&\!=&\!\!\!\!\!\!\!\!\!\!\!\!\!\!\!\!\!\!\!\!\!BT \!+ \sum_{i} Q_i(t_k)
\lambda_i -\sum_{i} Q_i(t_k) r_i^* \nonumber\\
&\;\;\;\;\;\;\;\;\;\;\;\;\;\;\;\;\;\;\;\;\;+& \mathbb{E}\left[ R^*(t_k)-\!R_T(t_k)\big|\mathbf{Q}(t_k) \right].
\label{eq:app-d_drift-mid}
\end{eqnarray}
Now we bound the last term. For all $\delta_2 >0$ we have
\begin{eqnarray} &\!\!\!\!\!\!\!\!\!\!\!\!\!\!\!\!\mathbb{ E}&
\!\!\!\!\!\!\!\!\!\!\!\!\!\!\!\!\left[
R^*(t_k)-R_T(t_k)\big|\mathbf{Q}(t_k) \right] = \nonumber \\
&\!\!\!\!\!\!\!\!\!\!\!\!\!\!\!\!=&\!\!\!\!\!\!\!\!\!\!\!\!\!\!\!\!
\mathbb{ E}\!\left[ R^*(t_k)\!-\!R_T(t_k)\big|\mathbf{Q}(t_k), R^*(t_k)\!-\!R_T(t_k) \ge \delta_2||\mathbf{Q}(t_k)|| \right] \nonumber \\
&\;\;\;\;\;.& \!\!\!\!\mathbb{P}\left(R^*(t_k)-R_T(t_k) \ge \delta_2||\mathbf{Q(t_k)}||\; \big|\mathbf{Q}(t_k)\right) \nonumber \\
&\;\!\!\!\!\!\!\!\!\!\!\!\!+& \!\!\!\!\!\!\!\!\!\!\!\!
\!\!\mathbb{E} \!\left[ R^*(t_k)\!-\!R_T(t_k\big)\big|\mathbf{Q}(t_k), R^*(t_k)\!-\!R_T(t_k) \!< \!\delta_2||\mathbf{Q}(t_k)|| \right] \nonumber\\
&\;\;\;\;\;.& \!\!\!\!\mathbb{P}\left(R^*(t_k)-R_T(t_k) < \delta_2||\mathbf{Q(t_k)}|| \; \big| \mathbf{Q}(t_k)\right)\nonumber \\
&\!\!\!\!\!\!\!\!\!\!\!\!\!\!\!\!\le&\!\!\!\!\!\!\!\!\!\!\!\!\!\!\!
\big(\!\sum_{i}Q_i(t_k) \big) \mathbb{P}\!\left
(|R^*\!(t_k)\!-\!\!R_T(t_k)| \!\ge\! \delta_2||\mathbf{Q}(t_k)|| \;
\!\big| \mathbf{Q}(t_k)\right) \nonumber\\
&+&\!\!\!\!\!\delta_2 ||\mathbf{Q}(t_k)||, \label{eq:app-d_expec_bound}
\end{eqnarray}
where we bound the first expectation by $\sum_{i}Q_i(t_k)$ by using $||\mathbf{r}^*|| <1$, the second expectation by $\delta_2||\mathbf{Q}(t_k)||$ and the second probability by 1. By Schwartz inequality we have
\begin{eqnarray}
&\!\!\!\!\!\!\!\!\!\!\!\!\mathbb{P}&\!\!\!\!\!\!\!\!\!\!\!\!\!\!\!\!\!\!\!\!\left(|R^*(t_k)-R_T(t_k)|\ge\delta_2||\mathbf{Q(t_k)}||\;\big| \mathbf{Q}(t_k)\right) \nonumber \\
&\;\;\;\;\;\;\;\;\;\;\;\le&
\mathbb{P}\left(||\mathbf{r}^*-\hat{\mathbf{r}}_T(t_k)|| \ge
\delta_2 \big| \mathbf{Q}(t_k)\right).\label{eq:app-d_prob_bound}
\end{eqnarray} 
Using (\ref{eq:MDP_rate_conv}) and (\ref{eq:app-d_prob_bound}) 
in (\ref{eq:app-d_expec_bound}), we have
\beqn
\mathbb{E}[R^*(t_k)\!-\!R_T\!(t_k\!)\big|\mathbf{Q}(t_k)] \!\le\!\! \Big(\!\!\sum_{i}\!Q_i(t_k\!)\!\Big)c_1e^{-c_2\delta_2 ^2T}\!\! +\!\delta_2 ||\mathbf{Q}(t_k)||.\nonumber
\eeqn
Hence, using $||\mathbf{Q}(t_k)|| \le \sum_i Q_i(t_k)$, we bound (\ref{eq:app-d_drift-mid}) as
\begin{eqnarray*}
\frac{\Delta_T(t_k)}{2T} \! \leq\! BT\!\!\!\!&+& \!\!\!\!\!
\sum_{i} Q_i(t_k) \lambda_i\! -\! \sum_{i} Q_i(t_k) r_i^* \\
\!\!\!\!\!\!\!\!&+& \!\!\!\!\!\Big(\sum_{i}Q_i(t_k) \!\Big)\left(c_1e^{-c_3\delta_2 ^2T} \!\!+\! \delta_2\right).\nonumber
\end{eqnarray*}
Therefore, calling $\delta \doteq c_1e^{-c_3\delta_2 ^2T} +
\delta_2$, we have
\begin{equation}
\frac{\Delta_T(t_k)}{2T} \!\leq \!\!BT \!+\!\! \sum_{i} Q_i(t_k)
\lambda_i\! -\!\! \sum_{i} Q_i(t_k) r_i^* \!\!+\!
\delta\!\sum_{i}Q_i(t_k). \label{eq:drift_comparison}
\end{equation}
Now for $\boldsymbol{\lambda}$ strictly inside the
$\mathbf{\delta}$-\emph{stripped} stability region
$\mathbf{\Lambda}_s^{\delta}$,
there exist a small $\xi >0$ such that $\boldsymbol{\lambda} +
\xi.\mathbf{1} = \mathbf{r} -\delta \mathbf{1}$, for some
$\mathbf{r} \in \mathbf{\Lambda}_s$. Utilizing this and the fact that
$\sum_{i} Q_i(t)(r_i- r_i^*) \le 0$ by definition of the FBDC policy in Algorithm~\ref{alg:FBDC}, we have,
\beqn \frac{\Delta_T(t_k)}{2T} \leq BT - \Big(\sum_{i}Q_i(t_k)
\Big)\xi. \label{eq:FBDC_proof_final} \eeqn
Therefore, the queue sizes
have negative drift when $\sum_{i} Q_i(t_k)$ is larger than $BT/\xi$.
\emph{This establishes stability of the queue sizes at the frame boundaries $t=kT,\, k=\{0,1,2,...\}$}
for $\boldsymbol \lambda$ within the
$\delta$-\emph{stripped} stability region
$\mathbf{\Lambda}_s ^{\delta}$. To see this, 
take expectations with respect to $\mathbf{Q}(t_k)$ to have
\begin{eqnarray*}
\mathbb{E}[L(\mathbf{Q}(t_{k+1}))] -
\mathbb{E}[L(\mathbf{Q}(t_k))] &\le&  2BT^2 - 2\xi T\mathbb{E}\left[\sum_{i}Q_i(t_k)\right]. 
\end{eqnarray*}
Writing a similar expression over the frame boundaries $t_k, k \in
\{0,1,2,...,K \}$, summing them and telescoping these expressions leads to
\begin{equation*}
L(\mathbf{Q}(t_{K})) -
L(\mathbf{Q}(0)) \le 2KBT^2  - 2\xi T \!\!\sum_{k=0}^{K-1} \!\!\mathbb{E}\left[\sum_{i}Q_i(t_k)\right]. 
\end{equation*}
Using $L(\mathbf{Q}(t_K))\ge 0$ and $L(\mathbf{Q}(0))=0$, we have 
\begin{equation*}
\frac{1}{K} \sum_{k=0}^{K-1}\mathbb{E}\left[\sum_{i}Q_i(t_k)\right] \le \frac{BT}{\xi } < \infty. 
\end{equation*}
This implies that
\begin{equation*}\label{eq:VFMW_final_drift8}
\displaystyle \limsup_{K\rightarrow \infty} \frac{1}{K} \sum_{k=0}^{K-1}\sum_i \mathbb{E}[Q_i(t_{k})] \le \frac{BT}{\xi} < \infty.
\end{equation*}
This establishes stability (as defined in Definition \ref{def:stab}) at the frame boundaries $t_k, k \in \{0,1,2,... \}$.

\emph{For any given time $t\in(t_k,t_{k+1})$ we have}
$Q_i(t) \le Q_i(t_k) + \sum_{\tau=0}^{T-1}A_i(t_k+\tau)$. Therefore,
$\mathbb{E}[Q_i(t)] \le \mathbb{E}[Q_i(t_k)] + T\lambda_i \le \mathbb{E}[Q_i(t_k)] + TA_{\max}$. Hence
stability at the frame boundaries implies the overall stability of the system. 
Finally, $\delta = c_1e^{-c_3 \delta_2^2T} + \delta_2$ for any $\delta_2 >0$. 
Therefore choosing $\delta_2$ appropriately (for example, $\delta_2=T^{-0.5+\delta_3}$ for some small $\delta_3>0$), we have that $\delta(T)$ is a decreasing function of $T$. Therefore, for any $\delta >0$, one can find $T$ such that the hypothesis of the theorem holds.

\section*{Appendix E - Proof of Lemma~\ref{lem:My_FBDC_comp} }
We first establish that
\begin{equation*}
\Psi=\frac{\sum_{i} Q_i(t) \hat{r}_i}{\sum_{i} Q_i(t) r_i^{*}} \ge 0.90.
\end{equation*}
Considering the mappings in tables \ref{Fig:myopic_table_eps_small} and
\ref{Fig:myopic_optimal_table}, for the regions of $\epsilon$ where the OLM policy and the
optimal policy ``choose'' the same corner point, we have $\Psi=1$.
In the following we analyze the ratio $\Psi$ in the regions where the two
policies choose different corner points, which we call 
``discrepant'' regions. We will use $Q_1$ and $Q_2$ instead of
$Q_1(t)$ and $Q_2(t)$ for notational simplicity. We first consider the case $Q_2>Q_1$, and divide the proof into separate cases for different regions of $\epsilon$
values.
\vspace{2mm}
\\
\textbf{Weighted Departure-Rate Ratio Analysis, Case 1: $\epsilon <\epsilon_c$ }\\
 Note that the following inequality always holds: $\frac{2-\epsilon}{1-\epsilon}
> \frac{(1+\epsilon-\epsilon^2)}{(1-\epsilon)}$. 
However, we have $\frac{2-\epsilon}{1-\epsilon} = \frac{(1-\epsilon)^2}{\epsilon}$ for $\epsilon = \epsilon_t \doteq 0.245$ for the
case of $\epsilon < \epsilon_c = 0.293$.\\
\\
\textbf{Case 1.1: $\epsilon < \epsilon_t$ }\\
For this case we have $\frac{2-\epsilon}{1-\epsilon} <
\frac{(1-\epsilon)^2}{\epsilon}$.\\
\textit{Discrepant Region 1: $\frac{(1-\epsilon)^2}{\epsilon} < \frac{Q_2}{Q_1}< \frac{1-\epsilon}{\epsilon}$}\\
In this case the OLM policy chooses the corner point $b_1$
whereas the optimal policy chooses the corner point $b_0$.
Therefore,
\begin{eqnarray*}
\Psi \!\!\!\!\!&=&\!\!\!\! \frac{Q_1\big(
\frac{(1-\epsilon)^2}{4}\big)+ Q_2 \big(
\frac{1}{2}-\frac{\epsilon}{4}\big)}  {Q_2\frac{1}{2}} \ge 1-
\frac{\epsilon}{2}+
\frac{(1-\epsilon)^2}{2}\frac{\epsilon}{1-\epsilon}\nonumber\\
&=& 1-\frac{\epsilon^2}{2} \ge 0.9700.
\end{eqnarray*}
\textit{Discrepant Region 2: $\frac{(1+\epsilon-\epsilon^2)}{1-\epsilon} < \frac{Q_2}{Q_1}< \frac{2-\epsilon}{1-\epsilon}$}\\
In this case the OLM policy chooses the corner point $b_2$
whereas the optimal policy chooses the corner point $b_1$.
Therefore,
\begin{eqnarray*}
\Psi \!\!\!\!\!&=&\!\!\!\! \frac{Q_1\big(
\frac{3}{8}-\frac{\epsilon}{2}+\frac{\epsilon}{8(2-\epsilon)}\big)+
Q_2 \big( \frac{3}{8}-\frac{\epsilon}{8(2-\epsilon)}\big)} {Q_1\big(
\frac{(1-\epsilon)^2}{4}\big)+ Q_2 \big(
\frac{1}{2}-\frac{\epsilon}{4}\big)} \\
&=& \frac{
\frac{3}{8}-\frac{\epsilon}{2}+\frac{\epsilon}{8(2-\epsilon)}+
\frac{Q_2}{Q_1} \big(
\frac{3}{8}-\frac{\epsilon}{8(2-\epsilon)}\big)} {
\frac{(1-\epsilon)^2}{4}+ \frac{Q_2}{Q_1} \big(
\frac{1}{2}-\frac{\epsilon}{4}\big)} \ge 0.9002.
\end{eqnarray*}
This is a minimization of a function of two variables for all
possible $\epsilon$ values in the interval $0 \le \epsilon \le
\epsilon_t$, and the ratio $\frac{Q_2}{Q_1}$ in the interval
$\frac{(1+\epsilon-\epsilon^2)}{1-\epsilon} < \frac{Q_2}{Q_1}<
\frac{2-\epsilon}{1-\epsilon}$.
\\
\\
\textbf{CASE 1.2: $\epsilon_t < \epsilon < \epsilon_c $}\\
For this case we have $\frac{2-\epsilon}{1-\epsilon} \ge
\frac{(1-\epsilon)^2}{\epsilon}$\\
\textit{Discrepant Region 1: $\frac{(2-\epsilon)}{(1-\epsilon)} < \frac{Q_2}{Q_1}< \frac{1-\epsilon}{\epsilon}$}\\
In this case the OLM policy chooses the corner point $b_1$
whereas the optimal policy chooses the corner point $b_0$.
Therefore,
\begin{eqnarray*}
\Psi\!\!\!\!\!&=&\!\!\!\! \frac{Q_1\big(
\frac{(1-\epsilon)^2}{4}\big)+ Q_2 \big(
\frac{1}{2}-\frac{\epsilon}{4}\big)}  {Q_2\frac{1}{2}} \ge 1-
\frac{\epsilon}{2}+
\frac{(1-\epsilon)^2}{2}\frac{\epsilon}{1-\epsilon}\\
&=& 1-\frac{\epsilon^2}{2} \ge 0.9500.
\end{eqnarray*}
\textit{Discrepant Region 2: $\frac{(1-\epsilon)^2}{\epsilon} < \frac{Q_2}{Q_1}< \frac{2-\epsilon}{1-\epsilon}$}\\
In this case the OLM policy chooses the corner point $b_2$
whereas the optimal policy chooses the corner point $b_0$.
Therefore,
\begin{eqnarray*}
\Psi\!\!\!\!\!&=&\!\!\!\! \frac{Q_1\big(
\frac{3}{8}-\frac{\epsilon}{2}+\frac{\epsilon}{8(2-\epsilon)}\big)+
Q_2 \big( \frac{3}{8}-\frac{\epsilon}{8(2-\epsilon)}\big)} { Q_2
\frac{1}{2}} \\
&\ge&
\!\!\!\!\!\big(\frac{1-\epsilon}{2-\epsilon}\big)\big(\frac{3}{4}-\epsilon+\frac{\epsilon}{4(2-\epsilon)}\big)
+ \frac{3}{4}-\frac{\epsilon}{4(2-\epsilon)} \ge 0.9150.
\end{eqnarray*}
\textit{Discrepant Region 3: $\frac{(1+\epsilon-\epsilon^2)}{1-\epsilon} < \frac{Q_2}{Q_1}< \frac{(1-\epsilon)^2}{\epsilon}$}\\
In this case the OLM policy chooses the corner point $b_2$
whereas the optimal policy chooses the corner point $b_1$.
Therefore,
\begin{eqnarray*}
\Psi\!\!\!\!\!&=&\!\!\!\! \frac{Q_1\big(
\frac{3}{8}-\frac{\epsilon}{2}+\frac{\epsilon}{8(2-\epsilon)}\big)+
Q_2 \big( \frac{3}{8}-\frac{\epsilon}{8(2-\epsilon)}\big)} {Q_1\big(
\frac{(1-\epsilon)^2}{4}\big)+ Q_2 \big(
\frac{1}{2}-\frac{\epsilon}{4}\big)} \\
&\ge& \frac{
\frac{3}{8}-\frac{\epsilon}{2}+\frac{\epsilon}{8(2-\epsilon)}+
\frac{Q_2}{Q_1} \big(
\frac{3}{8}-\frac{\epsilon}{8(2-\epsilon)}\big)} {
\frac{(1-\epsilon)^2}{4}+ \frac{Q_2}{Q_1} \big(
\frac{1}{2}-\frac{\epsilon}{4}\big)} \ge 0.9474.
\end{eqnarray*}
Due to symmetry, the same bounds on $\Psi$ applies for
$Q_2<Q_1$.\\
\textbf{Weighted Departure-Rate Ratio Analysis, Case 2: $\epsilon \ge
\epsilon_c$ }\\
For the case where $\epsilon \ge \epsilon_c$, we have $(1-\epsilon)(3-2\epsilon) \le (1-\epsilon)/\epsilon$ and $\frac{1-\epsilon}{\epsilon} < \frac{2-\epsilon}{1-\epsilon}$.
Therefore, the only discrepant region between the FBDC and the OLM policies for
$Q_2>Q_1$ is given by
%
%
$(1-\epsilon)(3-2\epsilon) \le \frac{Q_2}{Q_1}<
\frac{1-\epsilon}{\epsilon}$, where for this interval the OLM policy chooses the corner
point $b_1$, whereas the FBDC policy chooses the corner point $b_0$.\\
\textit{Discrepant Region 1: $(1-\epsilon)(3-2\epsilon) < \frac{Q_2}{Q_1}< \frac{1-\epsilon}{\epsilon}$}\\
In this case the OLM policy chooses the corner point $b_1$
whereas the optimal policy chooses the corner point $b_0$.
Therefore,
\begin{eqnarray*}\label{eq:Psi_6}
\Psi\!\!\!\!\!&=&\!\!\!\! \frac{Q_1\big(
\frac{3}{8}-\frac{\epsilon}{2}+\frac{\epsilon}{8(2-\epsilon)}\big)+
Q_2 \big( \frac{3}{8}-\frac{\epsilon}{8(2-\epsilon)}\big)}
{Q_2\frac{1}{2}} \\
&\ge&
\!\!\!\!\big(\frac{\epsilon}{1-\epsilon}\big)\big(\frac{3}{4}-\epsilon+\frac{\epsilon}{4(2-\epsilon)}\big)
+ \frac{3}{4}-\frac{\epsilon}{4(2-\epsilon)} \ge 0.914.
\end{eqnarray*}\\
Due to symmetry, the same bound on $\Psi$ applies for
$Q_2<Q_1$.
Combining all the cases, for all $\epsilon \in [0,0.5]$, we have
that $\Psi \ge 0.90$ for all possible $Q_1$ and $Q_2$.

Now the following drift expression for the
OLM policy can be derived similarly to the derivation of (\ref{eq:drift_comparison}) used in the proof
of Theorem~\ref{thm:FBDC}:
\begin{equation*}
\frac{\Delta_T(t_k)}{2T} \!\leq \!BT \!+\!\! \sum_{i} Q_i(t_k)
\lambda_i\! -\!\! \sum_{i} Q_i(t_k) \hat{r}_i \!\!+\!
\delta_4\!\sum_{i}Q_i(t_k),
\end{equation*}
where $\delta_4(T)$ 
is a decreasing function of $T$. 
Using (\ref{eq:Psi_bound}) 
\begin{equation*} \frac{\Delta_T(t_k)}{2T} \!\!\leq\!\! B T\! +\!\!
\sum_{i} \!Q_i(t_k) \lambda_i\! - 0.9\!\!\sum_{i}\! Q_i(t_k) r_i^{*}
\!+\! \delta_4 \!\!\sum_{i}\!Q_i(t_k).
\end{equation*}
Using an argument similar to that for (\ref{eq:FBDC_proof_final}) 
we have that
for $(\lambda_1,\lambda_2)$ strictly inside the 0.9 fraction
of the $\delta_4$-\emph{stripped} stability region,
there exist a small $\xi > 0$ such that $(\lambda_1,\lambda_2)
+ (\xi,\xi) = 0.9(r_1,r_2) -(\delta_4,\delta_4)$, for some
$\mathbf{r}=(r_1,r_2) \in \mathbf{\Lambda_s}$. Substituting this
expression for $(\lambda_1,\lambda_2)$ and using $\sum_{i} Q_i(t)(r-
r_i^*) \le 0$ we have,
\begin{eqnarray*}
\frac{\Delta_T(t_k)}{T} \leq (B+ K)T \!\!\!\!\! &+&\!\!\!\!\!
0.9\sum_{i} Q_i(t_k)(r-
r_i^*) \nonumber\\
&\!\!\!\!\!\!\!\!\!\!\!\!\!\!\!\!\!\!\!\!\!\!\!\!\!\!\!\!\!\!\!\!\!\!\!\!\!\!\!\!\!\!\!\!\!\!\!\!\!\!\!\!\!\!\!\!\!\!\!-&\!\!\!\!\!\!\!\!\!\!\!\!\!\!\!\!
\!\!\!\!\!\!\!\!\!\!\!\!\!\!\!\!\!\!\!\!
\bigg(\!\!\sum_{i}\!Q_i(t_k) \!\bigg)\delta_4 -\! \bigg(\!\!\sum_{i}\!Q_i(t_k)
\!\bigg)\xi\! +\!\! \bigg(\!\!\sum_{i}\!Q_i(t_k) \!\bigg)\delta_4.\nonumber
\end{eqnarray*}
After cancelations we have,
\beqn
\frac{\Delta_T(t_k)}{2T} \leq BT - \Big(\sum_{i}Q_i(t_k)
\Big)\xi. \nonumber
\eeqn
Therefore, using an argument similar to the proof of Theorem~\ref{thm:FBDC} in Appendix~C,
the system is stable for arrival rates within at least the 0.9
fraction of $\delta_4$-\emph{stripped} stability region, 
where
$\delta_4(T)$ is a decreasing function of $T$.
%
%
%

\section*{Appendix F - Proof of Lemma \ref{lemma:sum_tp_Nqueues}}

We follow similar steps to the proof of the sum-throughput upper bound for the case of two queues in Appendix A. In order to obtain an expression for $r_i, i \in \{1,...,N\}$, we sum the $2^{N-1}$ equations in (\ref{eq:saf_balance_N}) for which the server location $m$ is $i$ and the channel process of queue $i$, $C_i$, is 1. This gives for all $i \in \{1,...,N \}$
\begin{eqnarray*}
p_{10}r_i = &-& \!\!\!\sum_{\substack{s:\,m=i\\\;\;\;\;C_i=1}}\sum_{a\neq i} \mathbf{x}(s,a) \;\;+\;\; p_{01}\!\!\! \sum_{\substack{s:\,m=i\\\;\;\;\;C_i=0}} \mathbf{x}(s;i)\\
&+& \!\!\!(1-p_{10})\!\!\sum_{\substack{s:\,m \neq i\\\;\;\;\;C_i=1}}\!\!\mathbf{x}(s;i) \;+\; p_{01}\!\!\!\sum_{\substack{s:\,m\neq i\\\;\;\;\;C_i=0}}\mathbf{x}(s;i).
\end{eqnarray*}
Summing $r_i$ over all queues and using the normalization condition $\sum_s\sum_a \mathbf{x}(s,a)=1$, we have
\begin{eqnarray*}
(p_{10}+p_{01})\sum_{i=1}^Nr_i &=& \!\!p_{01} \;-\; \sum_{i=1}^N\sum_{j\neq i}\sum_{\substack{s:\,m=i\\C_i=1,C_j=0}}\!\!\mathbf{x}(s;j)\\
&-&\!\!(p_{01}+p_{10})\sum_{i=1}^N\sum_{j\neq i}\sum_{\substack{s:\,m=i\\C_i=1,C_j=1}}\mathbf{x}(s;j)\\
&+&\!\!\!(1\!-\!p_{01}\!-\!p_{10})\sum_{i=1}^N\sum_{j\neq i}\!\!\sum_{\substack{s:\,m=i\\C_i=0,C_j=1}}\!\!\mathbf{x}(s;j).
\end{eqnarray*}
From Corollary \ref{cor:X_vortex}, there exists a stationary-deterministic policy $\boldsymbol{\pi}$ that solves this LP of maximizing $\sum_i r_i(\mathbf{x})$ over the state-action polytope $\mathbf{X}$. Therefore, under this policy $\boldsymbol{\pi}$, at each state, at least one of the actions must have $0$ state-action frequency. Therefore, in order to maximize the sum-rate, the terms that have negative contribution to the sum-rate must be zero:
\begin{equation}\label{eq:sum_tp_Nqueues_2}
(p_{10}+p_{01})\!\sum_{i=1}^N \!r_i = p_{01} + (1-p_{01}-p_{10})\sum_{i=1}^N\sum_{j\neq i}\!\!\!\sum_{\substack{s:\,m=i\\C_i=0,C_j=1}}\mathbf{x}(s;j).
\end{equation}
Similar to the two-queue case in Appendix A, we utilize the expressions resulting from the fact that the steady state probability of each channel state vector is known. For instance, for $C_{0}^{(N)}\doteq \mathbb{P}\big((C_1,...,C_N)=(0,...,0)\big)=\frac{p_{10}^N}{(p_{10}+p_{01})^N}$, we have
\begin{equation*}
\sum_{i=1}^N\sum_{\substack{s:\,m=i\\C_j=0,\forall j}} \sum_{a\in A} \mathbf{x}(i,a) = C_{0}^{(N)}.
\end{equation*}
Summing these expressions we obtain
\begin{equation*}
\sum_{i=1}^N\sum_{j\neq i}\!\!\!\sum_{\substack{s:\,m=i\\C_i=0,C_j=1}}\mathbf{x}(s;j) = 1 - C_0^{(N)} - \sum_{i=1}^N r_i.
\end{equation*}
Combining this expression with (\ref{eq:sum_tp_Nqueues_2}) we obtain
\begin{equation*}
\sum_{i=1}^N r_i = 1 - C_{0}^{(N)} - \big( p_{10}(1-C_0^{(N)}) - p_{01} C_0^{(N)} \big).
\end{equation*}

\section*{Appendix G - Proof of Theorem~\ref{thm:FBDC_N}}
This proof follows very similar lines to the proof of Theorem~\ref{thm:FBDC} in Appendix~D.
The following $T$-step conditional Lyapunov drift expression can easily be derived
similar to (\ref{eq:drift_intermediate}).
\begin{eqnarray}
\Delta_{T}(t_k) &\!\!\!\!\!\le& \!\!\!\!\!NBT^2 + 2T \sum_{i=1}^N 
Q_i(t_k)\lambda_i \nonumber\\
&\!\!\!\!\!\!\!\!\!\!\!\!\!\!\!\!\!\!\!\!\!\!\!\!\!\!\!\!\!\!\!\!\!\!\!\!\!\!\!\!\!\!-&\!\!\!\!\!\!\!\!\!\!\!\!\!\!\!\!\!\!\!\!\!\!\!\!2\sum_{\ell}
Q_i(t_k) \mathbb{E}\left[ \sum_{\tau=0}^{T-1} D_i(t_k+\tau)
\big|\mathbf{Q}(t_k)\! \right]\!\!, \nonumber 
\end{eqnarray}
where
$B\doteq 1+A_{\max}^2$.
%
%
%
Recall the definition of the reward functions $\overline{r}_i(\mathbf{s}_t,\mathbf{a}_t)$, $i \in \{1,...,N\}$, 
in (\ref{eq:r_i-N})
and let
$\overline{r}_i(\mathbf{s}_t,\mathbf{a}_t)$ be the reward function associated with applying
policy $\boldsymbol{\pi}^*$ given in the definition of the FBDC policy in Algorithm~\ref{alg:FBDC_Nqueues} to the saturated system.
Let $\overline{r}_i(t)$ denote $\overline{r}_i(\mathbf{s}_t,\mathbf{a}_t)$ for notational simplicity, $i \in \{1,...,N\}$.
Note again that $r_i(t)$ is equal to $D_i(t)$, since $D_i(t)$ is the \emph{service opportunity} given to link $i$ at time slot $t$. Now let $\mathbf{r}^* = (r_i^*)_i$ be the infinite horizon
average rate associated with policy $\boldsymbol{\pi}^*$. Let $\mathbf{x}^*$ be
the optimal vector of state-action frequencies corresponding to
$\boldsymbol{\pi}^*$. Define the time-average empirical reward from queue $i$ in
the saturated system, $\hat{r}_{T,i}(t_k)$, 
$i \in \{1,...,N\}$ by\\
\beq \hat{r}_{ T,i }(t_k) \doteq \frac{1}{T}\sum_{\tau=0}^{T-1}
\overline{r}_i(t_k+\tau)\nonumber.  
\eeq
Similarly, define the time average empirical state-action frequency
vector $\hat{\mathbf{x}}_{ T}(t_k;\mathbf{s},\mathbf{a})$.
\begin{equation*}
\displaystyle \hat{\mathbf{x}}_T(t_k;\mathbf{s},\mathbf{a}) \doteq \frac{1}{T}\sum_{\tau=t_k}^{t_k+T-1} I_{\{\mathbf{s}_{\tau}=\mathbf{s},\mathbf{a}_{\tau}=\mathbf{a}\}},
\end{equation*}
where $I_E$ is the indicator function of an event $E$, i.e., $I_E=1$ if $E$ occurs and $I_E=0$ otherwise.
%
Using the definition of the reward functions 
in (\ref{eq:r_i-N}), we have that 
\begin{equation*}
\hat{r}_{ T,i }(t_k) =  \sum_{\mathbf{s}\in \mathcal{S}}\sum_{\mathbf{a}\in\mathcal{A}}\overline{r}_i(\mathbf{s},\mathbf{a}) \hat{\mathbf{x}}_{ T}(t_k;\mathbf{s},\mathbf{a}),\, i \in \{1,...,N\},
\end{equation*}
and $\hat{\mathbf{r}}_{ T}(t_k) = (\hat{r}_{ T,1 }(t_k))_{i}$. Similarly, we have
\begin{equation*}
r_i^* =  \sum_{\mathbf{s}\in \mathcal{S}}\sum_{\mathbf{a}\in\mathcal{A}}\overline{r}_i(\mathbf{s},\mathbf{a}) \mathbf{x}^*(\mathbf{s},\mathbf{a}),\, i \in \{1,...,N\}.
\end{equation*}
%
%
Again utilizing Lemma 4.1 in \cite{shie05}, we have that
for every choice of initial state
distribution, there exists constants $c_1$ and $c_2$ such that
\begin{equation}\label{eq:MDP_rate_conv_N}
\mathbb{P}(||\hat{\mathbf{r}}_{ T}(t_k)-\mathbf{r}^*|| \ge \delta_1)
\le c_1 e^{-c_2\delta_1^2 T}, \; \forall T\ge1, \forall \delta_1
>0.
\end{equation}
Furthermore, convergence of $\hat{\mathbf{r}}_{ T}(t_k)$ to $\mathbf{r}^*$ is w.p. 1.
%
%
%
%
%
%
%
%
%
Now let 
$R_T(t_k) \doteq\sum_{i} Q_i(t_k) \hat{r}_{ T,i }(t_k)$ and
$R^*(t_k)\doteq \sum_{i} Q_i(t_k) r_i^*$.  We rewrite the drift
expression:
\begin{eqnarray}
\frac{\Delta_T(t_k)}{2T} \!\!\!\!\!\!\!\!\!\!\!\!\!\!\!\!\!\!\!\!\!\!&\leq&\!\!\!\!\!\!\!\!\!\!\!\!\!\!\!\!\!\! \!\!\!\frac{NBT}{2} \!+
\sum_{i}Q_i(t_k) \lambda_i -
\mathbb{ E}\left[ R_T(t_k) \big|\mathbf{Q}(t_k) \right] \nonumber\\
&\!=&\!\!\!\!\!\!\!\!\!\!\!\!\!\!\!\!\!\!\!\!\!\frac{NBT}{2} \!+ \sum_{i} Q_i(t_k)
\lambda_i -\sum_{i} Q_i(t_k) r_i^* \nonumber\\
&\;\;\;\;\;\;\;\;\;\;\;\;\;\;\;\;\;\;\;\;\;+& \mathbb{E}\left[ R^*(t_k)-\!R_T(t_k)\big|\mathbf{Q}(t_k) \right].
\label{eq:app-d_drift-mid_N}
\end{eqnarray}
The last term can be bounded similarly to (\ref{eq:app-d_expec_bound}) to have for all $\delta_2 >0$
\begin{equation}
\mathbb{ E}\left[
R^*(t_k)-R_T(t_k)\big|\mathbf{Q}(t_k) \right]
\le \big(\!\sum_{i}Q_i(t_k) \big) \mathbb{P}\!\left
(|R^*\!(t_k)\!-\!\!R_T(t_k)| \!\ge\! \delta_2||\mathbf{Q}(t_k)|| \;
\!\big| \mathbf{Q}(t_k)\right)
+\delta_2 ||\mathbf{Q}(t_k)||. \label{eq:app-d_expec_bound_N}
\end{equation}
By Schwartz inequality we have
\begin{eqnarray}
&\!\!\!\!\!\!\!\!\!\!\!\!\mathbb{P}&\!\!\!\!\!\!\!\!\!\!\!\!\!\!\!\!\!\!\!\!\left(|R^*(t_k)-R_T(t_k)|\ge\delta_2||\mathbf{Q(t_k)}||\;\big| \mathbf{Q}(t_k)\right) \nonumber \\
&\;\;\;\;\;\;\;\;\;\;\;\le&
\mathbb{P}\left(||\mathbf{r}^*-\hat{\mathbf{r}}_T(t_k)|| \ge
\delta_2 \big| \mathbf{Q}(t_k)\right).\label{eq:app-d_prob_bound_N}
\end{eqnarray} 
Using (\ref{eq:MDP_rate_conv_N}) and (\ref{eq:app-d_prob_bound_N}) 
in (\ref{eq:app-d_expec_bound_N}), we have
\beqn
\mathbb{E}[R^*(t_k)\!-\!R_T\!(t_k\!)\big|\mathbf{Q}(t_k)] \!\le\!\! \Big(\!\!\sum_{i}\!Q_i(t_k\!)\!\Big)c_1e^{-c_2\delta_2 ^2T}\!\! +\!\delta_2 ||\mathbf{Q}(t_k)||.\nonumber
\eeqn
Hence, using $||\mathbf{Q}(t_k)|| \le \sum_i Q_i(t_k)$, we bound (\ref{eq:app-d_drift-mid_N}) as
\begin{eqnarray*}
\frac{\Delta_T(t_k)}{2T} \! \leq\! \frac{NBT}{2}\!\!\!\!&+& \!\!\!\!\!
\sum_{i} Q_i(t_k) \lambda_i\! -\! \sum_{i} Q_i(t_k) r_i^* \\
\!\!\!\!\!\!\!\!&+& \!\!\!\!\!\Big(\sum_{i}Q_i(t_k) \!\Big)\left(c_1e^{-c_3\delta_2 ^2T} \!\!+\! \delta_2\right).\nonumber
\end{eqnarray*}
Therefore, calling $\delta \doteq c_1e^{-c_3\delta_2 ^2T} +
\delta_2$, we have
\begin{equation*}
\frac{\Delta_T(t_k)}{2T} \!\leq \!\!\frac{NBT}{2} \!+\!\! \sum_{i} Q_i(t_k)
\lambda_i\! -\!\! \sum_{i} Q_i(t_k) r_i^* \!\!+\!
\delta\!\sum_{i}Q_i(t_k). 
\end{equation*}
Now for $\boldsymbol{\lambda}$ strictly inside the
$\mathbf{\delta}$-\emph{stripped} stability region
$\mathbf{\Lambda}_s^{\delta}$,
there exist a small $\xi >0$ such that $\boldsymbol{\lambda} +
\xi.\mathbf{1} = \mathbf{r} -\delta \mathbf{1}$, for some
$\mathbf{r} \in \mathbf{\Lambda}_s$. Utilizing this and the fact that
$\sum_{i} Q_i(t)(r_i- r_i^*) \le 0$ by definition of the FBDC policy in Algorithm~\ref{alg:FBDC}, we have,
\beqn \frac{\Delta_T(t_k)}{2T} \leq \frac{NBT}{2} - \Big(\sum_{i}Q_i(t_k)
\Big)\xi. \label{eq:FBDC_proof_final} \eeqn
Therefore, the queue sizes
have negative drift when $\sum_{i} Q_i(t_k)$ is larger than $\frac{NBT}{2\xi}$.
\emph{This establishes stability of the queue sizes at the frame boundaries $t=kT,\, k=\{0,1,2,...\}$}
for $\boldsymbol \lambda$ within the
$\delta$-\emph{stripped} stability region
$\mathbf{\Lambda}_s ^{\delta}$ (see e.g., \cite[Theorem 3]{neely05}).

\emph{For any given time $t\in(t_k,t_{k+1})$ we have}
$Q_i(t) \le Q_i(t_k) + \sum_{\tau=0}^{T-1}A_i(t_k+\tau)$. Therefore,
$\mathbb{E}[Q_i(t)] \le \mathbb{E}[Q_i(t_k)] + T\lambda_i \le \mathbb{E}[Q_i(t_k)] + TA_{\max}$. Hence
stability at the frame boundaries implies the overall stability of the system since the frame length is constant.
Finally, $\delta = c_1e^{-c_3 \delta_2^2T} + \delta_2$ for any $\delta_2 >0$. 
Therefore choosing $\delta_2$ appropriately (for example, $\delta_2=T^{-0.5+\delta_3}$ for some small $\delta_3>0$), we have that $\delta(T)$ is a decreasing function of $T$. Therefore, for any $\delta >0$, one can find $T$ such that the hypothesis of the theorem holds.
%


\bibliographystyle{IEEE}

\begin{thebibliography}{1}

\bibitem{AhmadTaraKrish09} S. Ahmad, L. Mingyan, T. Javidi, Q. Zhao, and B. Krishnamachari,
``Optimality of myopic sensing in multichannel opportunistic
access,'' \emph{IEEE Trans. Infor. Theory}, vol.~55, no. 9,
pp.~4040-4050, Sept. 2009.


\bibitem{AhmadLiu} S. Ahmad and M. Liu, ``Multi-channel opportunistic access: A case of
restless bandits with multiple plays,'' In \emph{Proc. Allerton'09}, Oct. 2009.

%
%
\bibitem{AkyilWang09} I. F. Akyildiz and X. Wang, ``Wireless mesh networks'', \emph{Wiley}, 2009.

\bibitem{altman99} E. Altman, ``Constrained Markov decision processes'', \emph{Chapman \& Hall}, London, 1999.

\bibitem{AltmanLevy94} E. Altman and H. Levy, ``Queueing in space'', \emph{Adv. Appl. Prob.}, vol. 26, no. 4, pp. 1095-1116, Dec. 1994.

\bibitem{AltSch91} E. Altman and A. Shwartz, ``Markov decision problems and state-action frequencies,'' \emph{SIAM J.
on Control and Opt.}, vol.~29, no~4, pp.~786-809, Jul. 1991.

\bibitem{AltKonsLiu92} E. Altman, P. Konstantopoulos, and Z. Liu, ``Stability, monotonicity and invariant quantities in
general polling systems,'' \emph{Queuing Sys.}, vol.~11, pp.~35-57,
Mar. 1992.

\bibitem{AltKush00} E. Altman and H. J. Kushner, ``Control of polling in presense of vacations in heavy traffic with
applications to satellite and mobile radio systems,'' \emph{SIAM J.
on Control and Opt.}, vol.~41, pp.~217-252, 2002.

\bibitem{BlakeLong09} L. Blake and M. Long, ''Antennas: Fundamentals, Design, Measurement,'' \emph{SciTech}, 2009.

\bibitem{BertTsit97} D. Bertsimas and J. Tsitsiklis, ``Introduction to Linear Optimization'', \emph{Athena Scientific}, 1997.

\bibitem{Berz} A.~Brzezinski and E.~Modiano, ``Dynamic Reconfiguration and Routing Algorithms for IP-over-WDM networks with Stochastic Traffic,'' \emph{IEEE Journal of Lightwave Tech.}, vol. 23, no. 10, pp. 3188-3205, Oct. 2005


\bibitem{Boxma90} O.J. Boxma, W.P. Groenendijk, and J.A. Weststrate, ``A Pseudoconservation law for service systems
with a polling table'', \emph{IEEE Trans. Commun.}, vol. 38, no. 10, pp. 1865–1870, 1990.


\bibitem{ChapKarSar05} P. Chaporkar, K. Kar, and S. Sarkar, ``Throughput guarantees through maximal scheduling
in wireless networks,'' In \emph{Proc. Allerton'05}, Sept. 2005.



\bibitem{Celik_Infocom11} G. D. \c{C}elik, L. B. Le, and E. Modiano, ``Scheduling in parallel queues with randomly varying connectivity and switchover delay,'' In \emph{Proc. IEEE INFOCOM'11 (Mini Conference)}, Apr. 2011.


\bibitem{Celik_ISIT11} G. \c{C}elik, S. Borst, P. Whiting, and E. Modiano, ``Variable frame based Max-Weight algorithms for networks with switchover delay,'' In \emph{Proc. IEEE ISIT'11}, Aug. 2011.

\bibitem{ChenZhaoSwami08}
Y. Chen, Q. Zhao, and A. Swami, ``Joint design and separation principle
for opportunistic spectrum access in the presence of sensing errors,''
\emph{IEEE Trans. Infor. Theory}, vol. 54, no. 5, pp. 2053–2071, May 2008.

\bibitem{EryilOzMod07} A. Eryilmaz, A. Ozdaglar, and E. Modiano, ``Polynomial complexity algorithms for full utilization of multi-hop
wireless networks,'' In \emph{Proc. IEEE INFOCOM'07}, May. 2007.

\bibitem{Gallager} R. G. Gallager, ``Discrete Stochastic Processes,'' Kluwer, 1996,
2nd edition online: http://www.rle.mit.edu/rgallager/notes.htm.

\bibitem{Geor_Neely_Tass06} L. Georgiadis, M. Neely, and L. Tassiulas,
``Resource Allocation and Cross-Layer Control in Wireless Networks,'' Now Publishers, 2006.


\bibitem{Gilbert60} E. N. Gilbert, ``Capacity of burst-noise channels,'' \emph{Bell Syst. Tech. J.},
vol. 39, pp. 1253–1265, Sept. 1960.

\bibitem{KarLuo07} K. Kar, X. Luo, and S. Sarkar, ``Throughput-optimal scheduling in
multichannel access point networks under infrequent channel
measurements,'' In \emph{Proc. IEEE Infocom'07}, May. 2007.

\bibitem{Long10}
L. B. Le, E. Modiano, C. Joo, and N. B. Shroff,
``Longest-queue-first scheduling under SINR interference model,'' In
\emph{Proc. ACM MobiHoc'10}, Sept. 2010.

\bibitem{Levy} H. Levy, M. Sidi, and O.L. Boxma, ``Dominance relations in polling systems,''
\emph{Queueing Systems}, vol. 6, pp. 155-172, Apr. 1990.

\bibitem{LiNeely10} C. Li and M. Neely, ``On achievable network capacity and throughput-achieving policies over Markov
ON/OFF channels,'' In \emph{Proc. WiOpt'10}, Jun. 2010.

\bibitem{LinShr05} X. Lin and N. B. Shroff, ``The impact of imperfect scheduling on
cross-layer rate control in wireless networks,'' In \emph{Proc. IEEE
Infocom'05}, Mar. 2005.

\bibitem{LiuNainTow} Z. Lui, P. Nain, and D. Towsley, ``On optimal polling policies,'' \emph{Queuing Sys.},
vol. 11, pp.~59-83, Jul. 1992.

\bibitem{shie05}
S.~Mannor and J.~N.~Tsitsiklis, ``On the emperical state-action
frequencies in Markov Decision Processes under general policies,''
{\em Mathematics of Operation Research}, vol. 30, no. 3, Aug. 2005.

\bibitem{ModBarry00}
E. Modiano and R. Barry, ``A novel medium access control protocol
for WDM-based LAN's and access networks using a Master/Slave
scheduler,'' {\em IEEE J. Lightwave Tech.}, vol. 18, no. 4, pp.
461--468, Apr. 2000.

\bibitem{ModShahZuss06} E. Modiano, D. Shah and G. Zussman, ``Maximizing throughput in
wireless networks via Gossip,'' In \emph{Proc. ACM
SIGMETRICS/Performance'06}, June 2006.



\bibitem{Neely10} M. J. Neely, ``Stochastic Network Optimization with Application to Communication and Queueing Systems,'' Morgan and Claypool, 2010.

\bibitem{neely09}
M. J. Neely, ``Stochastic optimization for Markov modulated networks with
application to delay constrained wireless scheduling,'' In {\em Proc.
IEEE CDC'09}, Dec. 2009.

\bibitem{neely05}
M. J. Neely, E. Modiano, and C. E. Rohrs, ``Dynamic power allocation
and routing for time varying wireless networks,'' {\em IEEE J. Sel.
Areas Commun.}, vol. 23, no. 1, pp. 89--103, Jan. 2005.

\bibitem{neely03}
M. J. Neely, E. Modiano, and C. E. Rohrs, ``Power allocation and
routing in multi-beam satellites with time varying channels,'' {\em
IEEE Trans. Netw.}, vol. 11, no. 1, pp. 138--152, Feb. 2003.

%

\bibitem{neely02}
M. J. Neely, E. Modiano, and C. E. Rohrs, ``Tradeoffs in delay guarantees
and computation complexity in nxn packet switches,'' I n \emph{Proc. CISS'02}, Mar. 2002.

%

\bibitem{puterman05} M. Puterman, ''Markov Decision Processes: Discrete Stochastic Dynamic Programming,'' \emph{Wiley}, 2005.


\bibitem{ShahWis06} D. Shah and D. J. Wischik, ``Optimal scheduling algorithms for input-queued switches,'' In \emph{Proc. IEEE
Infocom'06}, Mar. 2006.

\bibitem{AnnaEphr09}
A. Pantelidou, A. Ephremides, and A. Tits ``A cross-layer approach
for stable throughput maximization under channel state
uncertainty,'' {\em Wireless Networks}, vol. 15, no.5, pp. 555-569,
Jul. 2009.

\bibitem{Stolyar04}
A. L. Stolyar, ``Maxweight scheduling in a generalized switch: State
space collapse and workload minimization in heavy traffic,'' {\em
Annals of Appl. Prob.}, vol. 14, no. 1, pp. 1–-53, 2004.
%

\bibitem{tass92}
L.~Tassiulas and A.~Ephremides, ``Stability properties of
constrained queueing systems and scheduling policies for maximum
throughput in multihop radio networks,'' {\em IEEE Trans. Auto.
Control}, vol. 37, no. 12, pp. 1936-1948, Dec. 1992.

\bibitem{tass93}
L.~Tassiulas and A.~Ephremides, ``Dynamic server allocation to
parallel queues with randomly varying connectivity,'' {\em IEEE
Trans. Infor. Theory}, vol. 39, no. 2, pp. 466-478, Mar. 1993.

\bibitem{TolShu96} A. Tolkachev, V. Denisenko, A. Shishlov, and A. Shubov, ``High gain antenna systems for millimeter wave radars
with combined electronical and mechanical beam steering,'' In
\emph{Proc. IEEE Symp. Phased Array Sys. Tech.}, Oct. 2006.

\bibitem{vish06}
V. M. Vishnevskii and O. V. Semenova, ``Mathematical methods to
study the polling systems,'' {\em Auto. and Rem. Cont.}, vol. 67,
no. 2, pp. 173-220, Feb. 2006.
%

\bibitem{Walrand} J. Walrand, ``Queuing Networks,'' Englewood Cliffs, NJ:Prentice Hall, 1988.


\bibitem{WangChang96}
H.~Wang and P.~Chang, ``On verifying the first-order Markovian
assumption for a Rayleigh fading channel model,'' {\em IEEE Trans.
Veh. Tech.}, vol. 45, no. 2, pp. 353-357, May 1996.


\bibitem{WuSri06} X. Wu and R. Srikant, ``Bounds on the capacity
region of multi-hop wireless networks under distributed greedy
scheduling,'' In \emph{Proc. IEEE Infocom'06}, Mar. 2006.

\bibitem{shakk08}
L. Ying, and S. Shakkottai, ``On throughput-optimality with delayed
network-state information,'' In {\em Proc. Inform. Theory and
Applic. Workshop}, Jan. 2008.

\bibitem{ZhaoKrishLiu08} Q. Zhao, B. Krishnamachari, and K. Liu, ``On myopic sensing for multichannel
opportunistic access: Structure, optimality, and performance,'' \emph{ IEEE Trans. Wireless Commun.}, vol. 7, no. 12, pp. 5431–5440, Dec.
2008.

\bibitem{ZorziRao95}
M. Zorzi, R. Rao, and L. Milstein, ``On the accuracy of a
first-order Markov model for data transmission on fading channels,''
In {\em Proc. ICUPC'95}, 1995.
%

\end{thebibliography}


\begin{biography}[{\includegraphics[width=1in,height=1.25in,clip]{celik.eps}}]{G\"{u}ner D. \c{C}elik}
received his B.S. degree in Electrical and Electronics Engineering from the Middle East Technical University (METU), Ankara, Turkey with second highest ranking in 2005, and his S.M. degree in Electrical Engineering and Computer
Science Department (EECS) from the Massachusetts Institute of Technology (MIT)
in 2007. In August 2010, he was a visiting researcher at Bell Labs in Mathematics of Networks and Complex Systems Research Group in NJ and in summer 2007, he joined Microsoft Research (MSR) Systems and Networking Group in Camrbidge, UK as a research intern. He is
currently pursuing a Ph.D. degree in Laboratory For Information and Decision Systems (LIDS) in the EECS department at MIT under the supervision of Prof. Eytan Modiano. His research interests are in wireless network control, queuing theory and communication theory. In particular, he is currently working on scheduling, routing and resource allocation for network control and optimization, dynamic server scheduling in queuing systems and trajectory control for data gathering in mobile wireless networks.
\end{biography}


\begin{biography}[{\includegraphics[width=1in,height=1.25in,clip]{le.eps}}]{Long B. Le}
received the B.Eng. degree with highest distinction from Ho Chi Minh City University of Technology in 1999, the M.Eng. degree from Asian Institute of Technology (AIT) in 2002 and the PhD degree from University of Manitoba in 2007. He is currently an assistant professor at INRS-EMT, University of Quebec, Montreal, Quebec, Canada. Before that, he was a postdoctoral researcher first at University of Waterloo then at Massachusetts Institute of Technology. His current research interests include wireless protocol engineering and resource allocation, stochastic control and cross-layer design for communication networks. He is a Member of the IEEE.
\end{biography}


\begin{biography}[{\includegraphics[width=1in,height=1.25in,clip]{modiano.eps}}]{Eytan Modiano}
received his B.S. degree in Electrical Engineering and Computer Science from the University of Connecticut at Storrs in 1986 and his M.S. and PhD degrees, both in Electrical Engineering, from the University of Maryland, College Park, MD, in 1989 and 1992 respectively.  He was a Naval Research Laboratory Fellow between 1987 and 1992 and a National Research Council Post Doctoral Fellow during 1992-1993.  Between 1993 and 1999 he was with MIT Lincoln Laboratory where he was a project leader for MIT Lincoln Laboratory's Next Generation Internet (NGI) project.  Since 1999 he has been on the faculty at MIT; where he is a Professor in the Department of Aeronautics and Astronautics and the Laboratory for Information and Decision Systems (LIDS).  His research is on communication networks and protocols with emphasis on satellite, wireless, and optical networks. He is the co-recipient of the Sigmetrics 2006 Best paper award for the paper "Maximizing Throughput in Wireless Networks via Gossiping," and the Wiopt 2005 best student paper award for the paper "Minimum Energy Transmission Scheduling Subject to Deadline Constraints."

He is currently an Associate Editor for IEEE/ACM Transactions on Networking.  He had served as Associate Editor for IEEE Transactions on Information Theory, and as guest editor for IEEE JSAC special issue on WDM network architectures; the Computer Networks Journal special issue on Broadband Internet Access; the Journal of Communications and Networks special issue on Wireless Ad-Hoc Networks; and for IEEE Journal of Lightwave Technology special issue on Optical Networks.  He was the Technical Program co-chair for  IEEE Wiopt 2006, IEEE Infocom 2007, and ACM MobiHoc 2007.

\end{biography}

\end{document}